\documentclass{article}

\usepackage{amsmath}
\usepackage{amssymb}
\usepackage{amsthm}
\usepackage[inline]{enumitem}
\usepackage{imakeidx}
\usepackage{hyperref}

\usepackage{mathrsfs}           
\usepackage{stmaryrd}           

\usepackage{xcolor}
\usepackage{framed}     
\usepackage[most]{tcolorbox}
\usepackage{xparse}

\usepackage{footnote}

\usepackage{tikz}
\usetikzlibrary{calc} 

\usepackage{bm}                 
\usepackage{bbm}                

\usepackage{geometry}           
\usepackage{xparse}             

\usepackage[framemethod=TikZ]{mdframed}

\def\dif{\mathop{}\!\mathrm{d}}         
\DeclareMathOperator{\tr}{tr}           
\newcommand{\Inner}[2]{\left\langle {#1}, {#2} \right\rangle}            
\newcommand{\norm}[1]{\left\Vert {#1} \right\Vert}

\newcommand{\laplacian}{\Delta}

\DeclareMathOperator{\divr}{div}

\newcommand{\vol}{\mathrm{vol}}

\newcommand{\Rm}{\mathrm{Rm}}           

\DeclareMathOperator{\second}{I\!I}

\usepackage{amsthm}
\usepackage{tcolorbox}
\usepackage{boites}

\NewCommandCopy\originalproof\proof
\NewCommandCopy\endoriginalproof\endproof

\makeatletter
\renewenvironment{proof}[1][\proofname]
{%
    \def\bkvz@left{\vrule width 0.4pt\hskip 8pt}%
    \def\bkvz@right{}%
    \def\bkvz@top{}%
    \def\bkvz@bottom{}%
    \def\bkvz@set@linewidth{%
        \advance\linewidth -8.4pt\relax
    }%
    \bkcountfalse
    \begin{breakbox}%
        \originalproof[#1]%
}
{%
        \endoriginalproof
    \end{breakbox}%
}
\makeatother
\setlist[enumerate]{itemsep=0pt, parsep=0pt, topsep=0pt, partopsep=0pt} 

\newcommand{\regs}{\mathit{\Omega}}
\newcommand{\regt}{\mathit{\Lambda}}

\usepackage{float}
\usepackage{caption} 

\makeatletter
\def\namedlabel#1#2{\begingroup
#2%
\def\@currentlabel{#2}%
\phantomsection\label{#1}\endgroup
}
\makeatother

\DeclareMathOperator{\cl}{cl}
\DeclareMathOperator{\interior}{int}
\DeclareMathOperator{\supp}{supp}
\DeclareMathOperator{\diam}{diam}
\DeclareMathOperator{\dist}{dist}

\newcommand{\matr}[1]{\boldsymbol{#1}} 

\makeindex

\newcommand{\point}[1]{\bm{\mathrm{#1}}}

\renewcommand{\vec}[1]{\bm{#1}}
\renewcommand{\matr}[1]{\bm{\mathsf{#1}}}

\usepackage{aliascnt}
\newtheorem{theorem}{Theorem}[section]
\newaliascnt{lemma}{theorem}
\newtheorem{lemma}[lemma]{Lemma}
\aliascntresetthe{lemma}
\newaliascnt{proposition}{theorem}
\newtheorem{proposition}[proposition]{Proposition}
\aliascntresetthe{proposition}
\newaliascnt{corollary}{theorem}
\newtheorem{corollary}[corollary]{Corollary}
\aliascntresetthe{corollary}
\newtheorem{definition}{Definition}[section]

\newtheorem{remark}{Remark}[section]
\newtheorem*{claim}{Claim}

\newcommand{\seccur}{sec}
\newcommand{\cexp}{c\text{-}\exp}

\usepackage{cleveref}
\crefname{lemma}{Lemma}{Lemmata}
\Crefname{lemma}{Lemma}{Lemmata}
\crefname{proposition}{Proposition}{Propositions}
\Crefname{proposition}{Proposition}{Propositions}
\crefname{corollary}{Corollary}{Corollaries}
\Crefname{corollary}{Corollary}{Corollaries}

\usepackage[backend=biber, style=alphabetic]{biblatex}
\defbibheading{bibliography}[References]{\section*{#1}}

\title{Boundary $ C^{2, \gamma} $ regularity for optimal transport potentials under the A3s MTW condition}
\author{Yizhong Hong \and Feng Wang}

\begin{document}

\maketitle

\begin{abstract}
    We establish boundary regularity for optimal transport potentials under the A3s Ma--Trudinger--Wang condition. Let the source and target be bounded, mutually $ c $-convex domains in $ \mathbb{R}^n $, $ n\geqslant2 $, with $ C^{1,\beta} $ boundaries, and suppose that the cost is smooth, bi-twisted, and nondegenerate. For positive $ C^{0,\alpha} $ densities bounded away from zero and having equal total mass, we prove that both Kantorovich potentials belong to $ C^{2,\gamma} $ up to the boundary, where $ \gamma = \min\{\alpha,\beta\} $. Consequently, the optimal transport map is $ C^{1,\gamma} $ up to the boundary. The proof extends the monotonicity formula method of Collins and Tong to general costs through the geometry of the Kim--McCann metric. The resulting formula yields homogeneity of blow-up limits and uniform estimates for sections. These estimates give boundary obliqueness and allow comparison with solutions of the constant-density Monge--Ampère equation, leading to the claimed regularity.
\end{abstract}

\section{Introduction}\label{section::thesis::introduction}

This thesis studies boundary regularity for optimal transport with a general cost satisfying the strong Ma-Trudinger-Wang condition. The aim is to extend the monotonicity formula method of Collins and Tong \cite{collins2025boundary}, developed for the cost $ c(\point{x}, \overline{\point{x}}) = -\Inner{\point{x}}{\overline{\point{x}}} $, to general costs and obtain regularity under $ C^{1, \beta} $ assumptions on the boundaries of the source and target domains.

Let $ \regs, \regt \subset \mathbb{R}^n $ be bounded domains, and let $ \rho, \overline{\rho} $ be positive densities of equal total mass. We use the cost convention $ b = -c $. The Kantorovich potentials $ u, v $ satisfy $ u(\point{x}) + v(\overline{\point{x}}) \geqslant b(\point{x}, \overline{\point{x}}) $, with equality on the graph of the optimal transport map. The main result is the following.

\begin{theorem}\label{theorem::thesis::maintheorem}
    Let $ n\geqslant2 $ and $ 0 < \alpha, \beta < 1 $. Suppose that $ \regs, \regt \subset \mathbb{R}^n $ are bounded connected domains with $ C^{1, \beta} $ boundaries and are mutually $ c $-convex. Let the cost $ c $ be smooth on a neighborhood of $ \cl \regs \times \cl \regt $ and satisfy \ref{condition::thesis::A1}, \ref{condition::thesis::A2}, and \ref{condition::thesis::A3s} there. Suppose that $ \rho \in C^{0, \alpha}(\cl \regs) $ and $ \overline{\rho} \in C^{0, \alpha}(\cl \regt) $ have equal total mass and satisfy
    \[
        0 < \lambda^{-1} \leqslant \rho, \overline{\rho} \leqslant \lambda < \infty.
    \]
    Then the Kantorovich potentials satisfy
    \[
        u \in C^{2, \gamma}(\cl \regs), \qquad v \in C^{2, \gamma}(\cl \regt), \qquad \gamma := \min\{ \alpha, \beta \}.
    \]
    After fixing their additive constants, their $ C^{2, \gamma} $ norms are bounded in terms of the domains, the cost, the constants in the stated conditions, the density bounds, and the H\"older norms of the densities. In particular, if $ \beta \geqslant \alpha $, both potentials are $ C^{2, \alpha} $ up to the boundary.
\end{theorem}

The proof begins with the geometry of the optimal transport graph. The Kim-McCann metric expresses the MTW condition as a curvature condition, while the result of Kim, McCann, and Warren identifies the graph as a calibrated submanifold after a conformal change involving the densities. We use these identities to derive a monotonicity formula and estimate its error terms.

The monotonicity formula is then applied to blow-up limits of the potentials. We prove that the limits are homogeneous of degree two and use convergence of their sections and John ellipsoids to obtain uniform section estimates. These estimates give $ C^{1, 1 - \epsilon} $ regularity and are used to prove obliqueness at the boundary. Finally, comparison with solutions of the constant-density Monge-Amp\`ere equation and an iteration of these comparisons give the boundary $ C^{2, \gamma} $ estimates.

We assume smoothness of the cost to simplify the calculations and do not attempt to optimize the number of derivatives required. The conditions and notation used in the theorem are introduced in the following sections.

\noindent \textbf{Use of AI.} The initial proof was generated entirely by OpenAI's GPT-6 Astra through the Danus project\footnote{Danus project: \url{https://github.com/frenzymath/Danus}.}. The first author supplied a handwritten mathematical note deriving a monotonicity formula using $ z = \tilde{b}(\point{x}, \overline{\point{x}}) $ (see \cref{definition::thesis::normalization-at-a-contact,sec::thesis::insert-c} for the definition and motivation), with several estimates left unproved, and asked GPT-6 Astra to develop a complete proof of $ C^{2, \gamma} $ regularity. Astra first produced a complete draft on 19 September 2026 after about four days. The first author subsequently audited the proof, and further discussion with Astra led to corrections and simpler arguments for some parts.

The manuscript was written collaboratively by the authors with AI assistance, with AI also assisting in English translation and editing.

\section{The optimal transport problem}
Let $ M $, $ \overline{M} $ be manifolds. Assume that they are compact, or that they admit continuous embeddings into separable metric spaces in which their closures $ \cl M $, $ \cl \overline{M} $ are compact. Let $ \rho$ , $ \overline{\rho} $ be Borel probability measures on the two manifolds, respectively, and let a lower semicontinuous function $ c : \cl (M \times \overline{M}) \to \mathbb{R} \cup \{ +\infty \} $ be the transport cost.

The Kantorovich optimal transport problem seeks a measure $ \gamma \geqslant 0 $ on $ M \times \overline{M} $, whose marginals on $ M $ and $ \overline{M} $ are $ \rho $ and $ \overline{\rho} $, respectively, that minimizes the total transport cost
\[
    \int_{M \times \overline{M}} c(\point{x}, \overline{\point{x}}) \dif \gamma(\point{x}, \overline{\point{x}}).
\]
A transport plan $ \gamma $ minimizing this cost is called an optimal transport plan, which is the object sought in the optimal transport problem.

The Monge problem minimizes the transport cost over measurable maps $ \point{T}:M\to\overline{M} $ that transport $ \rho $ to $ \overline{\rho} $. If a Kantorovich minimizer is concentrated on the graph of such a map, that map also solves the Monge problem. Under suitable regularity assumptions on the cost and absolute continuity of the source measure, the twist condition \ref{condition::thesis::A1}, stated below, guarantees existence and uniqueness of the optimal map up to a set of source measure zero; see \cite{kim2010continuity}. This has motivated extensive study of its regularity.

For domains in Euclidean space and the cost $ c(\point{x},\overline{\point{x}})=-\Inner{\point{x}}{\overline{\point{x}}} $, Caffarelli \cite{caffarelli1992boundary,Ca1996}, Delano\"e \cite{delanoe1991classical}, and Urbas \cite{urbas1997second} established regularity results using geometric and partial differential equation methods, under suitable assumptions on the densities and domains. Convexity of the target plays a central role in these results.

After Caffarelli's extensive study of this standard case, Ma, Trudinger, and Wang \cite{ma2005regularity} formulated the MTW condition, motivated by the design of antenna reflectors. This condition involves fourth derivatives of the cost $ c $ and was used to establish interior estimates and regularity of optimal transport potentials for general costs under suitable assumptions on the data. Although the condition initially appears unusual, G. Loeper \cite{loeper2009regularity} subsequently proved its necessity for regularity: if the cost does not satisfy the weak MTW condition, there are examples of discontinuous maps even with smooth densities. In the same paper, Loeper obtained an interior H\"older exponent for optimal transport maps under the strong MTW condition. Later, A. Figalli and collaborators \cite{figalli2013holder} proved interior H\"older continuity under the weak MTW condition.

Loeper described the MTW expression as a cost-sectional curvature. Kim and McCann \cite{kim2010continuity} developed its geometric interpretation by constructing a pseudo-Riemannian metric on $ M\times\overline{M} $ and identifying the MTW expression with its unnormalized cross-curvature. The weak MTW condition corresponds to nonnegativity under the relevant null condition. They also showed that $ c $-segments are geodesics of this metric with one coordinate fixed, giving a geometric interpretation of the $ c $-exponential map.

Building on this, Kim, McCann, and M. Warren \cite{kim2010pseudo} proved that, after modifying the pseudo-metric by the densities, the graph of the optimal transport map is a calibrated submanifold and hence a stable maximal submanifold. Y. Yuan \cite{yuan2022monotonicity} developed a monotonicity argument on the maximal gradient graph for $ \det D^2u=1 $, following a mean-value inequality argument in \cite{gilbarg1998elliptic}, and used it to obtain Pogorelov\textquotesingle s Hessian estimate. T. Collins and F. Tong \cite{collins2025boundary} extended this method to general densities for $ c = -\Inner{\point{x}}{\overline{\point{x}}} $, thereby weakening the boundary assumptions in earlier $ C^{2, \alpha} $ estimates. 

The ideas of Y. Yuan and Collins--Tong are central to this proof. Their approach allows us to exploit the geometric structure associated with the cost function, avoiding the need to extend Caffarelli's section-based arguments directly to every stage of the general-cost proof.

\section{The Ma-Trudinger-Wang conditions and the Kim-McCann metric}\label{sec::thesis::MTW}
\subsection{The Kim-McCann metric}
Following the notation and theory of \cite{kim2010continuity}, this section introduces the Ma-Trudinger-Wang conditions and their reformulation by Kim and McCann using a pseudo-Riemannian metric.

Throughout, $ \point{x}, \point{y} $ denote points of $ M $, and barred symbols $ \overline{\point{x}}, \overline{\point{y}} $ denote points of $ \overline{M} $. Indices and functions are distinguished in the same way. When all indices belong to one space, as in $ c_{ij}(\point{x}, \overline{\point{x}}) = \partial_{x^i x^j} c(\point{x}, \overline{\point{x}}) $, we call this a derivative \textbf{on one side}, or \textbf{on the $ M $ side} when the space needs to be specified. If indices from both spaces occur, as in $ c_{i \overline{k}} = \partial_{x^i \overline{x}^{\overline{k}}} c $, we call it a \textbf{mixed derivative}. Bold symbols denote vectors, while nonbold symbols denote numbers or functions.

Take a subdomain $ N \subset M \times \overline{M} $ of the product manifold. It can be viewed as a multivalued map from $ M $ to $ \overline{M} $, with inverse $ N^{\ast} := \left\{ (\overline{\point{x}}, \point{x}) : (\point{x}, \overline{\point{x}}) \in N \right\} $. Also write $ c^{\ast}(\overline{\point{x}}, \point{x}) := c(\point{x}, \overline{\point{x}}) $ for the reversed cost.

Write $ T_{\point{x}} M $ for the tangent space of $ M $ at $ \point{x} $, and $ T^{\ast}_{\point{x}} M $ for the cotangent space. The tangent and cotangent spaces of the product manifold are direct sums of the corresponding spaces of the two factors. At $ (\point{x}, \overline{\point{x}}) $, the mixed derivative $ \bar{D}Dc $ induces a linear map from $ T_{\overline{\point{x}}}\overline{M} $ to $ T^{\ast}_{\point{x}}M $: for $ \vec{v}\in T_{\point{x}}M $ and $ \overline{\vec{w}}\in T_{\overline{\point{x}}}\overline{M} $, we have $ \Inner{\bar{D}Dc(\overline{\vec{w}})}{\vec{v}}=\bar{D}_{\overline{\vec{w}}}D_{\vec{v}}c $. Its transpose $ D\bar{D}c $ maps $ T_{\point{x}}M $ to $ T^{\ast}_{\overline{\point{x}}}\overline{M} $. Define
\[
    \mathfrak{h} := \frac{1}{2} \begin{pmatrix}
        0            & -\bar{D} D c \\
        -D \bar{D} c & 0
    \end{pmatrix},
\]
which is a symmetric bilinear form on the tangent space $ T_{(\point{x}, \overline{\point{x}})} N $.

\begin{definition}[\namedlabel{condition::thesis::A0}{\textbf{A0}}]
    The cost satisfies $ c \in C^4(N) $.
\end{definition}

Actually in this article we assume the cost is smooth, so \ref{condition::thesis::A0} automatically holds.

\begin{definition}[\namedlabel{condition::thesis::A1}{\textbf{A1}}, also called the bi-twist condition]
    A cost $ c \in C^1(N) $ satisfies the twist condition if, for every $ \point{x}_0 \in M $, the map $ \overline{\point{x}} \to -D_{\point{x}} c(\point{x}_0, \overline{\point{x}}) $ is injective.

    If both $ c $ and $ c^\ast $ satisfy this condition, the cost satisfies the bi-twist condition. Condition \textbf{A1} means the bi-twist condition.
\end{definition}

\begin{definition}[\namedlabel{condition::thesis::A2}{\textbf{A2}}, also called the nondegeneracy condition]
    A cost $ c \in C^2(N) $ satisfies \textbf{A2} if $ \bar{D} Dc: T_{\overline{\point{x}}} \overline{M} \to T_{\point{x}}^{\ast} M $ is bijective at every $ (\point{x}, \overline{\point{x}}) \in N $.
\end{definition}

Condition \ref{condition::thesis::A2} guarantees that $ \overline{\point{x}} \to -D_{\point{x}} c(\point{x}_0, \overline{\point{x}}) $ is locally a diffeomorphism. Its inverse is called the ``$ c $-exponential map'' in \cite{ma2005regularity}. This condition also ensures that, although $ \mathfrak{h} $ is not Riemannian, it is a pseudo-Riemannian metric. Nondegeneracy excludes zero eigenvalues, so the metric has signature $ (n, n) $. A tangent vector $ \vec{p} \oplus \overline{\vec{p}} \in T_{\point{x}} M \oplus T_{\overline{\point{x}}} \overline{M} $ is called a null vector if $ \mathfrak{h}(\vec{p} \oplus \overline{\vec{p}}, \vec{p} \oplus \overline{\vec{p}}) = 0 $.

Using the Riemann curvature tensor $ \Rm_{\alpha \beta \gamma \delta} $ induced by $ \mathfrak{h} $, define the sectional curvature by
\[
    \seccur_{(\point{x}, \overline{\point{x}})} \vec{P} \wedge \vec{Q} = \seccur_{(\point{x}, \overline{\point{x}})}^{N, \mathfrak{h}} \vec{P} \wedge \vec{Q} = \sum_{\alpha, \beta, \gamma, \delta = 1}^{2n} \Rm_{\alpha \beta \gamma \delta} P^{\alpha} Q^{\beta} P^{\gamma} Q^{\delta}.
\]
Unlike the usual definition on a Riemannian manifold, this definition does not divide by $ \mathfrak{h}(\vec{P}, \vec{P})\mathfrak{h}(\vec{Q}, \vec{Q}) - \mathfrak{h}(\vec{P}, \vec{Q})^2 $, since this factor vanishes for the null vectors considered in the MTW condition.

\begin{definition}[\namedlabel{condition::thesis::A3w}{\textbf{A3w}}]
    For every null vector $ \vec{p} \oplus \overline{\vec{p}} $,
    \begin{equation}\label{equation::thesis::seccur}
        \seccur_{(\point{x}, \overline{\point{x}})} (\vec{p} \oplus \vec{0}) \wedge (\vec{0} \oplus \overline{\vec{p}}) \geqslant 0.
    \end{equation}
\end{definition}

\begin{definition}[\namedlabel{condition::thesis::A3s}{\textbf{A3s}}]
    For every null vector $ \vec{p} \oplus \overline{\vec{p}} $, (\ref{equation::thesis::seccur}) holds, with equality only if $ \vec{p} = \vec{0} $ or $ \overline{\vec{p}} = \vec{0} $.
\end{definition}

In the terminology of the Kim-McCann metric, a cost satisfying \ref{condition::thesis::A0}, \ref{condition::thesis::A2}, and \ref{condition::thesis::A3w} is called \textbf{weakly regular}; it is called \textbf{regular} if \ref{condition::thesis::A3w} is replaced by \ref{condition::thesis::A3s}. If \ref{equation::thesis::seccur} also holds when $ \vec{p} \oplus \overline{\vec{p}} $ is not null, the pseudo-metric induced by $ \mathfrak{h} $ has \textbf{nonnegative cross-curvature}, denoted by \namedlabel{condition::thesis::NNCC}{\textbf{NNCC}}. If equality in \ref{equation::thesis::seccur} further implies $ \vec{p} = \vec{0} $ or $ \overline{\vec{p}} = \vec{0} $, the pseudo-metric has \textbf{positive cross-curvature}.

\begin{definition}[$ c $-convex domains]\label{definition::thesis::c-convex-domains}
    The domain $ \regs $ is $ c $-convex with respect to $ \overline{\point{x}}_0 \in \cl \regt $ if
    \[
        -D_{\overline{\point{x}}} c(\regs, \overline{\point{x}}_0) := \{ -D_{\overline{\point{x}}} c(\point{x}, \overline{\point{x}}_0) \colon \point{x} \in \regs \}
    \]
    is a convex subset of $ T^{\ast}_{\overline{\point{x}}_0} \overline{M} $. Similarly, $ \regt $ is $ c^{\ast} $-convex with respect to $ \point{x}_0 \in \cl \regs $ if $ -D_{\point{x}} c(\point{x}_0, \regt) $ is a convex subset of $ T^{\ast}_{\point{x}_0} M $.

    We say that $ \regs $ and $ \regt $ are \textbf{mutually $ c $-convex} if these two conditions hold for every $ \overline{\point{x}}_0 \in \cl \regt $ and every $ \point{x}_0 \in \cl \regs $, respectively.
\end{definition}

\subsection{Important results}
This subsection collects several results about the Kim-McCann metric from \cite{kim2010continuity}.

We begin with the Christoffel symbols and Riemann curvature tensor of the pseudo-metric $ \mathfrak{h} $:
\begin{lemma}\label{lemma::thesis::Christoffel-symbol-and-curvature}
    The only nonzero Christoffel symbols of $ \mathfrak{h} $ are
    \[
        {\Gamma_{ij}}^{k} = c^{k \overline{m}} c_{\overline{m} ij} \quad \text{and} \quad {\Gamma_{\overline{i}\overline{j}}}^{\overline{k}} = c^{\overline{k} m} c_{m \overline{i} \overline{j}}.
    \]
    A component of the Riemann curvature tensor vanishes if the numbers of indices from the two sides, namely from $ M $ and $ \overline{M} $, differ. Among the components with equal numbers of indices from the two sides, $ \Rm_{ij \overline{k} \overline{l}} = 0 $, while
    \[
        2 \Rm_{i \overline{k} j \overline{l}} = c_{ij \overline{a}} c^{\overline{a} b} c_{\overline{k} \overline{l} b} - c_{ij \overline{k} \overline{l}}.
    \]
\end{lemma}

We next consider geodesics for this metric. First recall the $ c $-exponential map defined in \cite{ma2005regularity}:
\begin{definition}[$ c $-exponential map]\label{definition::thesis::c-exp-map}
    Under condition \ref{condition::thesis::A1}, define the $ c $-exponential map
    \[
        \cexp_{\point{x}}: -D_{\point{x}}c(\point{x},\overline{N}(\point{x})) \to \overline{N}(\point{x})
    \]
    as follows: if
    \[
        \vec{p}^{\ast} = -Dc(\point{x}, \overline{\point{x}}) \in T^{\ast}_{\point{x}} M,
    \]
    then
    \[
        \cexp_{\point{x}} \vec{p}^{\ast} := \overline{\point{x}}.
    \]

    {Here $ \overline{N}(\point{x}):=\{\overline{\point{x}}:(\point{x},\overline{\point{x}})\in N\} $. Under \ref{condition::thesis::A1} and \ref{condition::thesis::A2}, this map is a diffeomorphism onto $ \overline{N}(\point{x}) $.}
\end{definition}

The image under the $ c $-exponential map of a line segment in $ T^{\ast}_{\point{x}} M $ is called a $ c $-segment in $ \overline{M} $. The following lemma shows that these segments do indeed correspond to geodesics:
\begin{lemma}
    Let $ \mathfrak{h} $ be the pseudo-metric on $ N \subset M \times \overline{M} $ defined by a cost satisfying \ref{condition::thesis::A2}. Fix $ \point{x}_0 \in M $. Take a line segment $ (1 - s) \vec{p}^{\ast} + s \vec{q}^{\ast} $ in $ T^{\ast}_{\point{x}_0} M $, and assume that the $ c $-exponential map $ \cexp $ is defined on it
    \footnote{The original statement does not assume \ref{condition::thesis::A1}, so a branch of a multivalued map should be considered. We assume this condition throughout and therefore omit this distinction.}.

    Consider the corresponding $ c $-segment
    \[
        \point{\sigma}(s) := (\point{x}_0, \cexp_{\point{x}_0}\left[(1 - s) \vec{p}^{\ast} + s \vec{q}^{\ast}\right]),
    \]
    Then $ \point{\sigma}(s) $ is a null geodesic in $ (N, \mathfrak{h}) $. Moreover, every geodesic in the totally geodesic submanifold $ \{\point{x}_0\} \times \overline{N}(\point{x}_0) $ can be represented in this way.

    In particular, in local coordinates on $ M \times \overline{M} $, if $ \point{\sigma}(s) $ is represented by $ (x^1, \cdots, x^n, \overline{x}^1(s), \cdots, \overline{x}^n(s)) $, it satisfies the geodesic equation
    \[
        \ddot{\overline{x}}^{\overline{m}} + c^{\overline{m} i} c_{i \overline{k} \overline{l}} \dot{\overline{x}}^{\overline{k}} \dot{\overline{x}}^{\overline{l}} = 0.
    \]
\end{lemma}

The following lemma shows that the MTW tensor appears in the fourth mixed derivative of the cost along two curves.
\begin{lemma}[Curvature and fourth derivatives, {\cite[Lemma 4.5]{kim2010continuity}}]\label{lemma::thesis::fourtime-derivative}
    Suppose the surface $ (s, t) \to (\point{x}(s), \overline{\point{x}}(t)) $ contains the curves $ \point{\sigma}(s) = (\point{x}(s), \overline{\point{x}}(0)) $ and $ \overline{\point{\tau}}(t) = (\point{x}(0), \overline{\point{x}}(t)) $. Then the vector field $ \vec{0}_{\point{x}(s)} \oplus \vec{\dot{\overline{x}}}(0) $ is parallel along $ \point{\sigma}(s) $.

    If $ \sigma(s) $ is additionally assumed to be a geodesic, then {with the curvature convention above,}
    \[
        \left. {-\frac{1}{2}} \frac{\partial^4}{\partial s^2 \partial t^2} \right|_{s = t = 0} c(\point{x}(s), \overline{\point{x}}(t)) = \seccur_{(\point{x}(0), \overline{\point{x}}(0))} \vec{\dot{\sigma}}(0) \wedge \vec{\dot{\overline{\tau}}}(0).
    \]
\end{lemma}

Finally, we consider the property called Double Mountain Above Sliding Mountain (DASM), which plays a central role in Loeper's regularity theory and subsequent work by Kim, McCann, Figalli, and others. Roughly, the ``double mountain'' is the maximum of two mountains of shape $ -c $, which are supporting surfaces of a $ c $-convex function: $ \max \{ -c(\point{x}, \overline{\point{x}}_1), -c(\point{x}, \overline{\point{x}}_2) \} $. The ``sliding mountain'' is a mountain of shape $ -c $ moving along the geodesic from $ \overline{\point{x}}_1 $ to $ \overline{\point{x}}_2 $. In the language of convex functions, this says that if
\[
    u(\point{z}) \geqslant u(\point{x}) + \Inner{\vec{p}_1}{\point{z} - \point{x}},
\]
\[
    u(\point{z}) \geqslant u(\point{x}) + \Inner{\vec{p}_2}{\point{z} - \point{x}},
\]
then, for every $ (1 - s)\vec{p}_1 + s \vec{p}_2 $, $ s \in [0, 1] $, we have
\[
    u(\point{z}) \geqslant u(\point{x}) + \Inner{(1 - s)\vec{p}_1 + s \vec{p}_2}{\point{z} - \point{x}}.
\]
\begin{theorem}[DASM]\label{theorem::thesis::DASM}
    Suppose $ c $ satisfies \ref{condition::thesis::A0} - \ref{condition::thesis::A2} and \ref{condition::thesis::A3w} on the product of mutually $ c $-convex domains $ \regs $ and $ \regt $. Take a $ c $-geodesic $ (\point{x}_0, \overline{\point{x}}(t)) $, $ 0\leqslant t\leqslant1 $, with $ \point{x}_0\in\regs $ and $ \overline{\point{x}}(t)\in\regt $. Then, for every $ \point{x}\in\regs $,
    \[
        f(t) := -c(\point{x}, \overline{\point{x}}(t)) + c(\point{x}_0, \overline{\point{x}}(t))
    \]
    the so-called ``sliding mountain,'' satisfies
    \[
        f(t) \leqslant \max\{ f(0), f(1) \}.
    \]
\end{theorem}
Under the stated domain hypotheses, DASM is equivalent to the weak MTW condition, \ref{condition::thesis::A3w}. DASM ensures that the $ c $-gradient set of a $ c $-convex function is connected.

\section{Optimal transport potentials and normalization}\label{section::thesis::potentials-and-normalization}

Write $ b := -c $ since the cost function always appears in its negative form. We can call this $ b $ a \textbf{benefit} function, since in Kantorovich duality it measures the profit of a company taking care of the transport problem as in the example of \cite[Section 5]{villani2008optimal}. But we still use the terms $ c $-convexity and $ c $-segments, following the established terminology in optimal transport.

\subsection{Optimal transport potentials}
We first recall the definition of $ c $-convexity:
\begin{definition}[$ c $-convex functions]
    A function $ u: \cl \regs \to \mathbb{R} \cup \{ +\infty \} $ is $ c $-convex if there exists a function $ v: \cl\regt\to\mathbb{R}\cup\{+\infty\} $ such that
    \[
        u(\point{x}) = v^c(\point{x}) = \sup_{\overline{\point{x}} \in \cl \regt} \{b(\point{x}, \overline{\point{x}}) - v(\overline{\point{x}})\},
    \]
    where $ v^c $ denotes the $ c $-transform of $ v $.
\end{definition}

For a function $ u $ on the source, define
\[
    u^c(\overline{\point{x}}):=\sup_{\point{x}\in\cl\regs}\{b(\point{x},\overline{\point{x}})-u(\point{x})\}.
\]
Its $ c $-subdifferential is
\[
    \partial^c u:=\{(\point{x},\overline{\point{x}})\in\cl\regs\times\cl\regt:
    u(\point{x})+u^c(\overline{\point{x}})=b(\point{x},\overline{\point{x}})\}.
\]

Kantorovich duality guarantees the existence of a $ c $-convex function $ u $ such that the measure induced by the optimal transport map vanishes outside $ \partial^c u $, see \cite[Proof of Theorem D.1]{kim2010continuity}. That is, the optimal transport map satisfies
\[
    u(\point{x}) + u^c(\point{T}(\point{x})) - b(\point{x}, \point{T}(\point{x})) = 0.
\]
Moreover, by the definition of $ c $-convexity, for every $ \point{z} \in \regs $,
\[
    u(\point{z}) + u^c(\point{T}(\point{x})) - b(\point{z}, \point{T}(\point{x})) \geqslant 0.
\]

We now assume that $ u, \point{T} $ are sufficiently differentiable. The preceding two formulas show that $ u(\point{z}) + u^c(\point{T}(\point{x})) - b(\point{z}, \point{T}(\point{x})) $ attains its minimum at $ \point{z} = \point{x} $, so
\begin{equation}\label{equation::thesis::formula-first-derivative}
    Du(\point{x}) - D_{\point{x}} b(\point{x}, \point{T}(\point{x})) = \vec{0},
\end{equation}
and
\begin{equation}\label{equation::thesis::positive-definite-second-derivative}
    D^2 u(\point{x}) - D^2_{\point{x}\point{x}} b(\point{x}, \point{T}(\point{x})) \geqslant 0.
\end{equation}
Using (\ref{equation::thesis::formula-first-derivative}), we also obtain
\begin{equation}\label{equation::thesis::formula-second-derivative}
    D^2 u(\point{x}) - D^2_{\point{x}\point{x}} b(\point{x}, \point{T}(\point{x})) - D^2_{\point{x} \overline{\point{x}}} b(\point{x}, \point{T}(\point{x})) D \point{T}(\point{x}) = \matr{0}.
\end{equation}

Equation (\ref{equation::thesis::formula-first-derivative}) gives
\[
    \point{T}(\point{x}) = (D_{\point{x}} b(\point{x}, \cdot))^{-1}\left( Du(\point{x}) \right) = \cexp_{\point{x}}\left(Du(\point{x})\right),
\]
Writing (\ref{equation::thesis::formula-second-derivative}) in local coordinates gives
\[
    (T^{\overline{k}})_{, i}(\point{x}) = b^{\overline{k} m}(\point{x}, \point{T}(\point{x})) (u_{mi}(\point{x}) - b_{mi}(\point{x}, \point{T}(\point{x}))).
\]
In what follows, as in \cite{ma2005regularity}, write
\[
    h_{ij} := u_{ij} - b_{ij},
\]
Then $ (T^{\overline{k}})_{, i} = b^{\overline{k} m} h_{mi} $.

Finally, using (\ref{equation::thesis::formula-second-derivative}), the Monge-Amp\'ere equation satisfied by $ u $ is
\[
    \det \left[ D^2 u - D^2_{\point{x}\point{x}} b \right]_{(\point{x}, \point{T}(\point{x}))} = \left|\det [D^2_{\point{x}\overline{\point{x}}} b]_{(\point{x}, \point{T}(\point{x}))}\right| \frac{\rho(\point{x})}{\overline{\rho}(\point{T}(\point{x}))}.
\]

\subsection{Normalization at a transport pair}

We use the following normalization throughout the boundary regularity, monotonicity, and blow-up arguments.

Take the continuous Kantorovich potentials in their canonical form
\[
    u(\point{x}) = \max_{\overline{\point{x}} \in \cl \regt} \{ b(\point{x}, \overline{\point{x}}) - v(\overline{\point{x}}) \}, \qquad v(\overline{\point{x}}) = \max_{\point{x} \in \cl \regs} \{ b(\point{x}, \overline{\point{x}}) - u(\point{x}) \},
\]
and fix their additive constants by $ u(\point{x}_{\mathrm{ref}}) = 0 $ at a fixed source point. Before uniqueness of the paired source and target points has been proved, write
\[
    \cl \Sigma := \{ (\point{x}, \overline{\point{x}}) \in \cl \regs \times \cl \regt \colon u(\point{x}) + v(\overline{\point{x}}) - b(\point{x}, \overline{\point{x}}) = 0 \}
\]
for the closed set of transport pairs. The optimal plan $ \gamma $ is supported on this set and has marginals $ \rho \dif \point{x} $ and $ \overline{\rho} \dif \overline{\point{x}} $, extended by zero to the boundaries. The canonical formulas and smoothness of the cost give uniform Lipschitz bounds for both potentials. For nearby points this follows by comparing the costs at the same opposite point, and for separated points it follows from boundedness of the cost. After normalization, the potentials also have common value bounds.

\begin{definition}\label{definition::thesis::normalization-at-a-contact}
Fix a transport pair $ (\point{x}_0, \overline{\point{x}}_0) \in \cl \Sigma $ and define
\[
    \begin{aligned}
        \tilde{u}(\point{x}) :={}& u(\point{x}) - u(\point{x}_0) - b(\point{x}, \overline{\point{x}}_0) + b(\point{x}_0, \overline{\point{x}}_0), \\
        \tilde{v}(\overline{\point{x}}) :={}& v(\overline{\point{x}}) - v(\overline{\point{x}}_0) - b(\point{x}_0, \overline{\point{x}}) + b(\point{x}_0, \overline{\point{x}}_0), \\
        \tilde{b}(\point{x}, \overline{\point{x}}) :={}& b(\point{x}, \overline{\point{x}}) - b(\point{x}, \overline{\point{x}}_0) - b(\point{x}_0, \overline{\point{x}}) + b(\point{x}_0, \overline{\point{x}}_0).
    \end{aligned}
\]
Thus $ \tilde{u}, \tilde{v} \geqslant 0 $, $ \tilde{u}(\point{x}_0) = \tilde{v}(\overline{\point{x}}_0) = 0 $, and
\[
    \tilde{u}(\point{x}) + \tilde{v}(\overline{\point{x}}) \geqslant \tilde{b}(\point{x}, \overline{\point{x}}),
\]
with equality at transport pairs. In particular, $ \tilde{b} $ is nonnegative on $ \cl\Sigma $.

\end{definition}

Subtracting functions of $ \point{x} $ alone and $ \overline{\point{x}} $ alone leaves the mixed derivatives of $ b $ unchanged. Thus the Kim-McCann metric is unchanged. Moreover,
\[
    \tilde{u}(\point{x}) + \tilde{v}(\overline{\point{x}}) - \tilde{b}(\point{x}, \overline{\point{x}}) = u(\point{x}) + v(\overline{\point{x}}) - b(\point{x}, \overline{\point{x}}),
\]
so the set of transport pairs and the optimal transport problem are unchanged. By abuse of notation, we will continue to write $ u, v, b $ for the normalized functions. Unless explicitly stated otherwise, this convention is applied at each fixed base transport pair, and any subsequent use of $ c $ means $ c = -b $ for the same normalization. In particular, $ b(\point{x}, \overline{\point{x}}_0) = b(\point{x}_0, \overline{\point{x}}) = 0 $, so $ b_i(\point{x}, \overline{\point{x}}_0) = b_{\overline{k}}(\point{x}_0, \overline{\point{x}}) = 0 $.

\section{Boundary $ C^{1, \alpha} $ regularity}\label{section::thesis::boundary-C1a}

We prove the following boundary regularity result, following the geometric method of \cite{loeper2009regularity}.

\begin{theorem}[Boundary H\"older regularity]\label{theorem::thesis::C1a-estimate-under-A3s}
    Let $ n \geqslant 2 $. Assume that $ \regs, \regt $ are bounded connected domains with $ C^{1,\beta} $ boundaries and that $ b $ is smooth on a neighborhood of $ \cl \regs \times \cl \regt $. Assume that the cost $ c = -b $ satisfies \ref{condition::thesis::A1}, \ref{condition::thesis::A2}, and \ref{condition::thesis::A3s} on this neighborhood. Assume that $ \regs $ and $ \regt $ are mutually $ c $-convex in the sense of \cref{definition::thesis::c-convex-domains}. The densities have equal total mass and satisfy
\[
    0 < \lambda^{-1} \leqslant \rho, \overline{\rho} \leqslant \lambda < \infty.
\]

    The set $ \cl\Sigma $ of transport pairs is the graph of a homeomorphism $ \vec{T}: \cl \regs \to \cl \regt $. Put $ \theta := 1/(4n-1) $. There is a constant $ C $, depending only on the domains, the cost, and the density bounds, such that
    \[
        |\vec{T}(\point{x}) - \vec{T}(\point{z})| \leqslant C|\point{x} - \point{z}|^\theta, \qquad |\vec{T}^{-1}(\overline{\point{x}}) - \vec{T}^{-1}(\overline{\point{z}})| \leqslant C|\overline{\point{x}} - \overline{\point{z}}|^\theta.
    \]
    The normalized potentials satisfy
    \[
        \norm{u}_{C^{1,\theta}(\cl \regs)} + \norm{v}_{C^{1,\theta}(\cl \regt)} \leqslant C.
    \]
\end{theorem}
\begin{remark}
    Loeper obtained the exponent $ 1 / (4n - 1) $, while Liu improved it to $ 1 / (2n - 1) $ for bounded densities; see \cite[Theorem 1 and Remark 4.1]{liu2009holder}. We make no claim of novelty for the regularity result in this section. We include the proof to present the argument in the setting and notation used below. For convenience in the subsequent proofs, we shall use the weaker exponent $ 1 / (4n) $.
\end{remark}

\subsection{Supporting functions and $ c $-segments}

Fix a transport pair $ (\point{x}_0, \overline{\point{x}}_0) \in \cl \Sigma $ and use the normalization in \cref{definition::thesis::normalization-at-a-contact}.

The maps $ D_{\overline{\point{x}}} b(\cdot, \overline{\point{x}}_0) $ and their inverses have common Lipschitz bounds as the focus varies. Indeed, nondegeneracy gives uniform local inverse bounds, and injectivity on the compact closures gives a uniform separation of images for pairs of points outside those local neighborhoods. Smoothness bounds the derivatives of the maps and their inverses. The same assertions hold with source and target interchanged. Their convex images have a common diameter bound and contain balls of a common positive radius: apply the uniform local inverse bounds to a fixed ball compactly contained in the original domain. Their Jacobian determinants are bounded above and away from zero.

For $ \point{x}_1 \ne \point{x}_0 $, let $ \point{x}(s) $, $ 0 \leqslant s \leqslant 1 $, be the source $ c $-segment with focus $ \overline{\point{x}}_0 $, defined by
\[
    D_{\overline{\point{x}}} b(\point{x}(s), \overline{\point{x}}_0) = s D_{\overline{\point{x}}} b(\point{x}_1, \overline{\point{x}}_0).
\]
It lies in $ \cl \regs $. Set
\[
    \vec{a} := D_{\overline{\point{x}}} b(\point{x}_1, \overline{\point{x}}_0), \qquad \ell := |\vec{a}|.
\]
The preceding bounds give
\[
    C^{-1}|\point{x}_1 - \point{x}_0| \leqslant \ell \leqslant C|\point{x}_1 - \point{x}_0|, \qquad C^{-1}\ell \leqslant |\dot{\point{x}}(s)| \leqslant C\ell.
\]
Also,
\[
    D_{\overline{\point{x}}} b(\point{x}(s), \overline{\point{x}}_0) = s\vec{a}, \qquad b(\point{x}(s), \overline{\point{x}}_0) = 0.
\]

The following estimate is related to the quantitative double-mountain estimate of \cite[Proposition 5.1]{loeper2009regularity} and \cite[Proposition B.1]{kim2010continuity}. Those estimates contain an additional cubic remainder in the distance from the base point.

\begin{lemma}\label{lemma::thesis::endpoint-separation}
    There exist $ c_0 > 0 $ and $ 0 < \varrho \leqslant 1 $, independent of the base transport pair, such that for every source $ c $-segment above with $ 0 < \ell \leqslant 1 $, every $ 1/3 \leqslant s \leqslant 2/3 $, and every $ \overline{\point{x}} \in \cl \regt $ with $ |\overline{\point{x}} - \overline{\point{x}}_0| \leqslant \varrho $,
    \[
        \max\{ -b(\point{x}(s), \overline{\point{x}}), b(\point{x}_1, \overline{\point{x}}) - b(\point{x}(s), \overline{\point{x}}) \} \geqslant c_0\ell^2 |\overline{\point{x}} - \overline{\point{x}}_0|^2.
    \]
\end{lemma}
\begin{proof}
    Fix $ s $ and let $ \overline{\point{x}}_s(t) $ be the target $ c $-segment from $ \overline{\point{x}}_0 $ to $ \overline{\point{x}} $, with focus $ \point{x}(s) $. Write $ \overline{\vec{v}} := \overline{\point{x}} - \overline{\point{x}}_0 $. Uniform inverse bounds give, for $ |\overline{\vec{v}}| $ small,
    \[
        \overline{\point{x}}_s(t) = \overline{\point{x}}_0 + t\overline{\vec{v}} + O(|\overline{\vec{v}}|^2), \qquad \dot{\overline{\point{x}}}_s(t) = \overline{\vec{v}} + O(|\overline{\vec{v}}|^2).
    \]
    Define
    \[
        f_s(t) := \left. \frac{\dif^2}{\dif r^2} b(\point{x}(r), \overline{\point{x}}_s(t)) \right|_{r=s}.
    \]
    The two identities preceding the lemma give $ f_s(0) = f_s'(0) = 0 $. Since $ D_{\point{x}} b(\point{x}(s), \overline{\point{x}}_s(t)) $ is affine in $ t $, the terms containing $ \ddot{\point{x}}(s) $ vanish after differentiating twice in $ t $. By \cref{lemma::thesis::Christoffel-symbol-and-curvature}, the remaining terms are
    \[
        \begin{aligned}
            f_s''(t) ={}& \left( b_{ij\overline{k}\overline{l}} - b_{ij\overline{a}} b^{\overline{a}m} b_{m\overline{k}\overline{l}} \right) \dot{x}^i(s)\dot{x}^j(s)\dot{\overline{x}}_s^{\overline{k}}(t)\dot{\overline{x}}_s^{\overline{l}}(t) \\
            ={}& 2\Rm_{i\overline{k}j\overline{l}} \dot{x}^i(s)\dot{x}^j(s)\dot{\overline{x}}_s^{\overline{k}}(t)\dot{\overline{x}}_s^{\overline{l}}(t),
        \end{aligned}
    \]
    where the derivatives of $ b $ and the curvature are evaluated at $ (\point{x}(s), \overline{\point{x}}_s(t)) $.

    Strict MTW and compactness give a positive curvature bound on null pairs of unit vectors. By continuity, the same bound, reduced by a fixed factor, holds for pairs sufficiently close to being null. Differentiating the source-segment equation and using the preceding estimates gives
    \[
        b_{i\overline{k}}(\point{x}(s), \overline{\point{x}}_s(t))\dot{x}^i(s)\dot{\overline{x}}_s^{\overline{k}}(t) = \vec{a}\cdot\overline{\vec{v}} + O(\ell|\overline{\vec{v}}|^2).
    \]
    Consequently, there is a common $ \vartheta > 0 $ such that, if $ |\vec{a}\cdot\overline{\vec{v}}| \leqslant \vartheta\ell|\overline{\vec{v}}| $, reducing $ \varrho $ gives $ f_s''(t) \geqslant c_0\ell^2|\overline{\vec{v}}|^2 $. Integrating twice in $ t $ yields
    \[
        \frac{\dif^2}{\dif s^2} b(\point{x}(s), \overline{\point{x}}) \geqslant c_0\ell^2|\overline{\vec{v}}|^2.
    \]
    Hence the difference between the affine interpolation of the endpoint values and the value at $ s \in [1/3,2/3] $ is at least a fixed positive multiple of $ \ell^2|\overline{\vec{v}}|^2 $. The maximum of the two endpoint differences is at least this interpolation difference.

    In the remaining case, $ |\vec{a}\cdot\overline{\vec{v}}| > \vartheta\ell|\overline{\vec{v}}| $, Taylor expansion in the target variable gives
    \[
        \frac{\dif}{\dif s} b(\point{x}(s), \overline{\point{x}}) = \vec{a}\cdot\overline{\vec{v}} + O(\ell|\overline{\vec{v}}|^2).
    \]
    After reducing $ \varrho $, this derivative has constant sign and absolute value at least $ \vartheta\ell|\overline{\vec{v}}|/2 $. Integrating toward the appropriate endpoint over an interval of length at least $ 1/3 $ proves a lower bound $ \vartheta\ell|\overline{\vec{v}}|/6 $. Since $ \ell, |\overline{\vec{v}}| \leqslant 1 $, this implies the stated quadratic lower bound. The case $ \overline{\vec{v}} = \vec{0} $ is immediate.
\end{proof}

We also record the geometric volume estimate needed in the mass comparisons. It concerns subsets of the original source domain, including when the $ c $-segment lies on its boundary.

\begin{lemma}\label{lemma::thesis::volume-near-c-segment}
    For a source $ c $-segment as above with $ 0 < \ell \leqslant 1 $, and for $ 0 < \delta \leqslant \ell $ sufficiently small, there is a measurable set $ A_\delta \subset \regs $ such that every point of $ A_\delta $ is within distance $ \delta $ of some $ \point{x}(s) $ with $ 1/3 < s < 2/3 $, and
    \[
        |A_\delta| \geqslant c_0\ell\delta^{n-1}.
    \]
    The constants are independent of the base transport pair and of the segment.
\end{lemma}
\begin{proof}
    The convex set $ D_{\overline{\point{x}}} b(\regs, \overline{\point{x}}_0) $ contains a ball $ B_d(\vec{a}_0) $ with a common $ d > 0 $ and has uniformly bounded diameter. Choose a common sufficiently large $ K_0 $ and set $ t := \delta/K_0 < 1/2 $. Define $ A_\delta $ as the inverse image under $ D_{\overline{\point{x}}} b(\cdot, \overline{\point{x}}_0) $ of
    \[
        \left\{ (1-t)D_{\overline{\point{x}}} b(\point{x}(s), \overline{\point{x}}_0) + t(\vec{a}_0 + \vec{z}) \colon \frac{1}{3} < s < \frac{2}{3}, \quad \vec{z}\cdot\vec{a} = 0, \quad |\vec{z}| < \frac{d}{2} \right\}.
    \]
    Convexity puts this set in the open image domain, even if the original segment is in its boundary. The uniform inverse Lipschitz bound and the choice of $ K_0 $ put each inverse image within distance $ \delta $ of the corresponding $ \point{x}(s) $. Its image volume is
    \[
        (1-t)\frac{\ell}{3}\omega_{n-1}\left( \frac{td}{2} \right)^{n-1}.
    \]
    The uniform Jacobian bounds therefore give the desired lower bound for the original source volume.
\end{proof}

\begin{remark}
    Loeper's proof of H\"older regularity uses a localized volume comparison; see \cite[Section 5]{loeper2009regularity}. Kim and McCann describe this argument as ``fitting sausages into balls'' in \cite[Appendix C]{kim2010continuity}. Here we use the interior ball $ B_d(\vec{a}_0) $ in the convex image domain to construct an interior sausage near the middle part of the $ c $-segment. Taking convex combinations with points in this ball keeps the sausage inside the domain, including when the segment lies on the boundary, while preserving the volume lower bound $ c_0\ell\delta^{n-1} $.
\end{remark}

\subsection{Uniqueness of the transport pairs on the closures}

\begin{lemma}\label{lemma::thesis::unique-closure-contacts}
    Every source point in $ \cl \regs $ is paired with exactly one target point, and every target point in $ \cl \regt $ is paired with exactly one source point. Consequently, the set $ \cl\Sigma $ of transport pairs is the graph of a homeomorphism $ \vec{T}: \cl \regs \to \cl \regt $.
\end{lemma}
\begin{proof}
    Suppose that $ \overline{\point{x}}_0 $ is paired with two distinct source points. By \cref{theorem::thesis::DASM}, with source and target interchanged, $ u $ vanishes on the source $ c $-segment joining them. Indeed, each supporting function along that segment is bounded above by the larger of its endpoint values, which are nonpositive; taking their supremum gives $ u \leqslant 0 $, while $ u \geqslant 0 $. Taking a shorter segment and choosing its first endpoint as $ \point{x}_0 $ if necessary, assume $ 0 < \ell \leqslant 1 $.

    For $ s \in [1/3,2/3] $, any target point $ \overline{\point{x}} $ paired with $ \point{x}(s) $ makes $ b(\point{x}(r), \overline{\point{x}}) $ attain a maximum at $ r=s $. Thus
    \[
        \left. \frac{\dif}{\dif r} b(\point{x}(r), \overline{\point{x}}) \right|_{r=s} = 0, \qquad \left. \frac{\dif^2}{\dif r^2} b(\point{x}(r), \overline{\point{x}}) \right|_{r=s} \leqslant 0.
    \]
    If $ \overline{\point{x}} \ne \overline{\point{x}}_0 $, take the target $ c $-segment $ \overline{\point{x}}_s(t) $ at $ \point{x}(s) $. Its defining equation and the first equality imply
    \[
        b_{i\overline{k}}(\point{x}(s), \overline{\point{x}}_s(t))\dot{x}^i(s)\dot{\overline{x}}_s^{\overline{k}}(t) = 0.
    \]
    Both velocities are nonzero by twist and nondegeneracy. The exact derivative identity in the proof of \cref{lemma::thesis::endpoint-separation} now gives $ f_s''(t) > 0 $ by strict MTW, with $ f_s(0) = f_s'(0) = 0 $. Therefore $ f_s(1) > 0 $, contradicting the second inequality. Thus every middle point is paired only with $ \overline{\point{x}}_0 $.

    For a middle point define the nonnegative difference
    \[
        w_s(\overline{\point{x}}) := v(\overline{\point{x}}) - b(\point{x}(s), \overline{\point{x}}).
    \]
    By the endpoint supporting inequalities and \cref{lemma::thesis::endpoint-separation},
    \[
        w_s(\overline{\point{x}}) \geqslant c_0\ell^2|\overline{\point{x}} - \overline{\point{x}}_0|^2 \quad \text{if } |\overline{\point{x}} - \overline{\point{x}}_0| \leqslant \varrho.
    \]
    On the compact set where $ s \in [1/3,2/3] $ and $ |\overline{\point{x}} - \overline{\point{x}}_0| \geqslant \varrho $, it has a positive minimum, unless that set is empty. This follows from the uniqueness just proved at the middle points.

    If $ (\point{x}, \overline{\point{x}}) $ is a transport pair with $ |\point{x} - \point{x}(s)| < \delta $, then the defining equality for transport pairs and $ u \geqslant 0 $ give
    \[
        w_s(\overline{\point{x}}) = b(\point{x}, \overline{\point{x}}) - u(\point{x}) - b(\point{x}(s), \overline{\point{x}}) \leqslant C\delta|\overline{\point{x}} - \overline{\point{x}}_0|.
    \]
    Here we used $ D_{\point{x}}b(\point{x}, \overline{\point{x}}_0) = 0 $ and smoothness of the cost on a neighborhood of the compact product. The positive minimum first excludes the targets at distance at least $ \varrho $ for small $ \delta $. The quadratic estimate then gives $ |\overline{\point{x}} - \overline{\point{x}}_0| \leqslant C_\ell\delta $ for all such transport pairs.

    Apply \cref{lemma::thesis::volume-near-c-segment}. Every transport pair whose source lies in $ A_\delta $ has its target in $ B_{C_\ell\delta}(\overline{\point{x}}_0) $. The marginal identities and the density bounds imply
    \[
        c_\ell\delta^{n-1} \leqslant \gamma(A_\delta \times \cl \regt) \leqslant \gamma\left( \cl \regs \times \left[\cl \regt \cap B_{C_\ell\delta}(\overline{\point{x}}_0)\right] \right) \leqslant C_\ell\delta^n,
    \]
    a contradiction as $ \delta \to 0 $. Hence each target is paired with at most one source point. Interchanging source and target proves the same assertion for source points. The canonical formulas ensure that every source and target point belongs to a transport pair. The set $ \cl\Sigma $ of transport pairs is closed in a compact product, so the resulting bijection and its inverse are continuous.
\end{proof}

\subsection{A uniform positive separation from a supporting mountain}

\begin{lemma}\label{lemma::thesis::uniform-positive-contact-gap}
    For every $ d > 0 $, there exists $ g(d) > 0 $, depending only on $ d $, the domains, the cost, and the density bounds, such that for every base transport pair,
    \[
        u(\point{x}) \geqslant g(d) \quad \text{if } |\point{x} - \point{x}_0| \geqslant d, \qquad v(\overline{\point{x}}) \geqslant g(d) \quad \text{if } |\overline{\point{x}} - \overline{\point{x}}_0| \geqslant d.
    \]
\end{lemma}
\begin{proof}
    In this proof only, $ b $ denotes the original benefit before normalization at a transport pair, and the canonical potentials are normalized by $ u_j(\point{x}_{\mathrm{ref}}) = 0 $. Suppose the source assertion fails. Take a sequence of density pairs satisfying the common bounds, canonical potentials $ u_j,v_j $, transport pairs $ (\point{x}_j,\overline{\point{x}}_j) $, and source points $ \point{z}_j $ with $ |\point{z}_j - \point{x}_j| \geqslant d $ such that
    \[
        u_j(\point{z}_j) + v_j(\overline{\point{x}}_j) - b(\point{z}_j,\overline{\point{x}}_j) \longrightarrow 0.
    \]
    The common value and Lipschitz bounds give uniformly convergent subsequences of both potentials on the closures. Take weak-* convergent subsequences of the densities in $ L^\infty $ on the fixed domains, weakly convergent subsequences of the optimal plans on the compact product, and convergent subsequences of the selected points. The limiting densities retain the same positive bounds. Testing the marginal identities shows that the limiting plan has those densities as marginals.

    Uniform convergence preserves both canonical formulas and the condition that the transport plan is supported on $ \{u+v=b\} $: the nonnegative functions $ u_j+v_j-b $ converge uniformly, and their integrals against the corresponding plans are zero. Thus the limit satisfies the same hypotheses. In the limit, two distinct source points are paired with the same target, because the source separation remains at least $ d $. This contradicts \cref{lemma::thesis::unique-closure-contacts}. The target assertion follows by interchanging source and target, and we take the smaller of the two constants.
\end{proof}

\subsection{A lower bound for the potentials}

\begin{lemma}[Lower bound for the potentials]\label{lemma::thesis::4n-lowerbound-for-u-v}
    There exist $ c_0,r_0 > 0 $, depending only on the domains, the cost, and the density bounds, such that for every base transport pair and every $ \point{x} \in \cl \regs $ with $ |\point{x} - \point{x}_0| \leqslant r_0 $,
    \[
        u(\point{x}) \geqslant c_0|\point{x} - \point{x}_0|^{4n}.
    \]
    The corresponding assertion holds for $ v $ with source and target interchanged.
\end{lemma}
\begin{proof}
    Smoothness on a neighborhood of the compact product and the vanishing of $ b $ on both axes give common estimates
    \[
        |D_{\point{x}}b(\point{x}, \overline{\point{x}})| \leqslant C|\overline{\point{x}} - \overline{\point{x}}_0|, \qquad |b(\point{x}, \overline{\point{x}})| \leqslant C|\point{x} - \point{x}_0|\,|\overline{\point{x}} - \overline{\point{x}}_0|.
    \]
    For nearby points these follow by Taylor expansion in the ambient neighborhood; when one separation is bounded below, the corresponding uniform bounds give the same estimates after increasing $ C $.

    Let $ \varrho $ be the radius in \cref{lemma::thesis::endpoint-separation}. There is a common $ r_* > 0 $ such that every transport pair with $ |\point{x} - \point{x}_0| < r_* $ satisfies $ |\overline{\point{x}} - \overline{\point{x}}_0| < \varrho $. Indeed, at a transport pair,
    \[
        0 \leqslant v(\overline{\point{x}}) = b(\point{x}, \overline{\point{x}}) - u(\point{x}) \leqslant C|\point{x} - \point{x}_0|,
    \]
    whereas targets at distance at least $ \varrho $ have $ v \geqslant g(\varrho) $ by \cref{lemma::thesis::uniform-positive-contact-gap}.

    Take $ \point{x}_1 \ne \point{x}_0 $ close to the base and the source $ c $-segment $ \point{x}(s) $ defined above. Put $ E := u(\point{x}_1) > 0 $, where positivity follows from uniqueness of the source point paired with the base target. The uniform Lipschitz bound gives $ E \leqslant C\ell $. Choose a common $ \ell_0 \leqslant 1 $ so small that $ E \leqslant 1 $ and every point at distance at most $ \ell $ from the segment has distance less than $ r_* $ from $ \point{x}_0 $, whenever $ \ell \leqslant \ell_0 $.

    For $ s \in [1/3,2/3] $, the two endpoint supporting inequalities give
    \[
        v(\overline{\point{x}}) - b(\point{x}(s), \overline{\point{x}}) \geqslant \max\{ -b(\point{x}(s), \overline{\point{x}}), b(\point{x}_1, \overline{\point{x}}) - b(\point{x}(s), \overline{\point{x}}) \} - E.
    \]
    If $ (\point{x}, \overline{\point{x}}) $ is a transport pair with $ |\point{x} - \point{x}(s)| < \delta \leqslant \ell $, the defining equality for transport pairs gives the upper bound
    \[
        v(\overline{\point{x}}) - b(\point{x}(s), \overline{\point{x}}) \leqslant C\delta|\overline{\point{x}} - \overline{\point{x}}_0|.
    \]
    The preceding choice puts this target within distance $ \varrho $ of $ \overline{\point{x}}_0 $. Applying \cref{lemma::thesis::endpoint-separation} yields
    \[
        c_0\ell^2|\overline{\point{x}} - \overline{\point{x}}_0|^2 \leqslant E + C\delta|\overline{\point{x}} - \overline{\point{x}}_0|,
    \]
    and consequently
    \[
        |\overline{\point{x}} - \overline{\point{x}}_0| \leqslant C\left( \frac{\sqrt{E}}{\ell} + \frac{\delta}{\ell^2} \right).
    \]

    Set $ \delta := \ell\sqrt{E} \leqslant \ell $, reducing $ \ell_0 $ if needed so that \cref{lemma::thesis::volume-near-c-segment} applies. Every transport pair whose source lies in the set $ A_\delta $ in that lemma has its target in $ B_{C\sqrt{E}/\ell}(\overline{\point{x}}_0) $. The original density bounds and the marginal identities give
    \[
        c_0\ell\delta^{n-1} \leqslant \gamma(A_\delta \times \cl \regt) \leqslant \gamma\left( \cl \regs \times \left[\cl \regt \cap B_{C\sqrt{E}/\ell}(\overline{\point{x}}_0)\right] \right) \leqslant C\left( \frac{\sqrt{E}}{\ell} \right)^n.
    \]
    Substituting the choice of $ \delta $ gives
    \[
        c_0\ell^n E^{\frac{n-1}{2}} \leqslant C E^{\frac{n}{2}}\ell^{-n}.
    \]
    Dividing by $ E^{(n-1)/2} > 0 $ proves $ E \geqslant c_0\ell^{4n} $. Since $ \ell $ is uniformly comparable to $ |\point{x}_1 - \point{x}_0| $, this is the asserted inequality. The target assertion follows by interchanging source and target.
\end{proof}

\subsection{H\"older continuity of the maps and gradients}

\begin{proof}[Proof of \cref{theorem::thesis::C1a-estimate-under-A3s}]
    By \cref{lemma::thesis::unique-closure-contacts}, the set $ \cl\Sigma $ of transport pairs is the graph of a homeomorphism between the closures. Take two transport pairs $ (\point{x}_0, \overline{\point{x}}_0) $ and $ (\point{x}_1, \overline{\point{x}}_1) $. If $ 0 < |\point{x}_1 - \point{x}_0| \leqslant r_0 $, the preceding lower bound and the estimate for $ b $ give
    \[
        c_0|\point{x}_1 - \point{x}_0|^{4n} \leqslant u(\point{x}_1) = b(\point{x}_1, \overline{\point{x}}_1) - v(\overline{\point{x}}_1) \leqslant C|\point{x}_1 - \point{x}_0|\,|\overline{\point{x}}_1 - \overline{\point{x}}_0|.
    \]
    Dividing and taking the power $ \theta = 1/(4n-1) $ gives
    \[
        |\point{x}_1 - \point{x}_0| \leqslant C|\overline{\point{x}}_1 - \overline{\point{x}}_0|^\theta.
    \]
    If $ |\point{x}_1 - \point{x}_0| \geqslant r_0 $, the uniform positive separation gives $ u(\point{x}_1) \geqslant g(r_0) $, while the defining equality for transport pairs gives $ u(\point{x}_1) \leqslant C|\overline{\point{x}}_1 - \overline{\point{x}}_0| $. Thus the target separation has a common positive lower bound, and the same estimate follows by boundedness of the domains. The zero-separation case is immediate. Interchanging source and target proves the other estimate.

    Set $ \vec{G}(\point{x}) := D_{\point{x}} b(\point{x}, \vec{T}(\point{x})) $. Smoothness of $ b $ and the map estimate give
    \[
        |\vec{G}(\point{x}) - \vec{G}(\point{z})| \leqslant C|\point{x} - \point{z}|^\theta.
    \]
    The supporting inequalities at $ \point{x} $ and $ \point{z} $ give
    \[
        b(\point{z}, \vec{T}(\point{x})) - b(\point{x}, \vec{T}(\point{x})) \leqslant u(\point{z}) - u(\point{x}) \leqslant b(\point{z}, \vec{T}(\point{z})) - b(\point{x}, \vec{T}(\point{z})).
    \]
    For nearby points, Taylor expansion of the cost in its ambient neighborhood yields
    \[
        -C|\point{z} - \point{x}|^2 \leqslant u(\point{z}) - u(\point{x}) - \vec{G}(\point{x})\cdot(\point{z} - \point{x}) \leqslant C|\point{z} - \point{x}|^{1+\theta} + C|\point{z} - \point{x}|^2.
    \]
    For separated points the uniform value and gradient-field bounds give the same first-order remainder estimate after increasing $ C $. Therefore
    \[
        |u(\point{z}) - u(\point{x}) - \vec{G}(\point{x})\cdot(\point{z} - \point{x})| \leqslant C|\point{z} - \point{x}|^{1+\theta}.
    \]
    This proves differentiability in the interior, with $ Du = \vec{G} $, and gives its H\"older extension and first-order remainder on the closure. Interchanging source and target gives
    \[
        Dv(\overline{\point{x}}) = D_{\overline{\point{x}}} b(\vec{T}^{-1}(\overline{\point{x}}), \overline{\point{x}})
    \]
    and the corresponding estimates. Together with the value bounds from the normalization, these prove the stated $ C^{1,\theta} $ estimates.
\end{proof}

\section{Step 1 for the monotonicity formula: basic formulae}\label{sec::thesis::frame}
In this section, we calculate some general identities using an arbitrary function $ z $. The choice of $ z $ will be determined in the next section. We use the potential identities and normalization from \cref{section::thesis::potentials-and-normalization}, and assume enough differentiability for the potentials, transport maps and densities.

\subsection{Geometry of the submanifold}
We next analyze the tangent and normal bundles of the graph $ \Sigma = (\point{x}, \point{T}(\point{x})) $ of the optimal transport map with respect to $ \mathfrak{h} $.

Using the tangent map induced by the embedding $ \iota: \regs \to \Sigma \subset M \times \overline{M} $, denote the tangent vectors of the submanifold by $ \iota_\ast(\partial_{x^i}) =: \vec{s}_i \in T_{\point{p}} \Sigma \subset T_{\point{p}} (M \times \overline{M}) $. We have $ \vec{s}_i := \vec{e}_i + b^{j \overline{l}} h_{ji} \vec{e}_{\overline{l}} $. 

The manifold $ \Sigma $ can be expressed as the level set $ \{ \vec{\Phi}(\point{x}, \overline{\point{x}}) = \vec{0} \} $ of the function $ \vec{\Phi}(\point{x}, \overline{\point{x}}) := \overline{\point{x}} - \point{T}(\point{x}) $, so the covectors $ \dif (\overline{x}^{\overline{k}} - T^{\overline{k}}(\point{x})) $ span its conormal bundle. Taking one half of their gradients with respect to $ \mathfrak{h} $ gives the normal frame
\[
    \vec{s}_{\overline{k}} := \frac{1}{2}\nabla (\overline{x}^{\overline{k}} - T^{\overline{k}}(\point{x})) = b^{\overline{k} i}(\vec{e}_i - b^{j \overline{l}} h_{ij} \vec{e}_{\overline{l}}).
\]

The inner products of the tangent vectors are
\[
    \begin{aligned}
        \Inner{\vec{s}_i}{\vec{s}_j} ={} & \Inner{\vec{e}_i + h_{is} b^{s \overline{k}} \vec{e}_{\overline{k}}}{\vec{e}_j + h_{jt} b^{t \overline{l}} \vec{e}_{\overline{l}}} \\
        ={}                                & \frac{1}{2}\left[b_{i \overline{l}} h_{jt} b^{t \overline{l}} + b_{j \overline{k}} h_{is} b^{s \overline{k}}\right]          \\
        ={}                                & h_{ij}.
    \end{aligned}
\]
and,
\[
    \begin{aligned}
        \Inner{\vec{s}_{\overline{k}}}{\vec{s}_{\overline{l}}} ={} & \Inner{b^{\overline{k} i}(\vec{e}_i - h_{is} b^{s \overline{a}} \vec{e}_{\overline{a}})}{b^{\overline{l} j}(\vec{e}_j - h_{jt}b^{t \overline{b}} \vec{e}_{\overline{b}})} \\
        ={}                                                          & b^{\overline{k} i} b^{\overline{l} j} \Inner{\vec{e}_i - h_{is} b^{s \overline{a}} \vec{e}_{\overline{a}}}{\vec{e}_j - h_{jt} b^{t \overline{b}} \vec{e}_{\overline{b}}}    \\
        ={}                                                          & b^{\overline{k} i} b^{\overline{l} j} \frac{1}{2} \left[ -b_{i \overline{b}} h_{jt} b^{t \overline{b}} - b_{j \overline{a}} h_{is} b^{s \overline{a}} \right]                   \\
        ={}                                                          & -b^{\overline{k} i} h_{ij} b^{j \overline{l}}.
    \end{aligned}
\]
We use $ |\cdot|_h $ to denote the norm induced by the positive definite metric $ h $ on $ T\Sigma $. For covectors on $ \Sigma $, the corresponding norm is defined using the inverse matrix $ (h^{ij}) $.

\paragraph{Conformal change}
By the result of Kim, McCann, and Warren \cite{kim2010pseudo}, $ \Sigma $ is a maximal submanifold for the following pseudo-metric, which is one half of the metric in their convention. Multiplication by a positive constant preserves maximality. We write
\[
    \mathfrak{h}^{\rho, \overline{\rho}} = \frac{1}{2} \left( \frac{\rho(\point{x})\overline{\rho}(\overline{\point{x}})}{|\det D \bar{D} b(\point{x}, \overline{\point{x}})|} \right)^{\frac{1}{n}} \begin{pmatrix}
        \matr{0}        & D\bar{D}b \\
        (D\bar{D}b)^\top & 0
    \end{pmatrix}.
\]

The metric $ \mathfrak{h}^{\rho, \overline{\rho}} $ is obtained from $ \mathfrak{h} $ by a conformal change. Define
\[
    f := \frac{1}{2n} \log \frac{\rho(\point{x})\overline{\rho}(\overline{\point{x}})}{|\det D \bar{D} b(\point{x}, \overline{\point{x}})|}.
\]
Then $ \mathfrak{h}^{\rho, \overline{\rho}} = e^{2f} \mathfrak{h} $.

\paragraph{The mean curvature vector}

Denote the mean curvature vector with respect to $ \mathfrak{h}^{\rho, \overline{\rho}} $ by $ \tilde{\vec{H}} = \tr_{\mathfrak{h}^{\rho, \overline{\rho}}}\tilde{\second} $. Since $ \Sigma $ is a maximal submanifold for this pseudo-metric, $ \tilde{\vec{H}} = \vec{0} $. By basic geometry calculations, we have
\begin{equation}\label{equation::thesis::conformal-mean-curvature-vector}
    \tilde{\vec{H}} = e^{-2f} \left[ \vec{H} - n \left(\nabla f\right)^{\perp} \right].
\end{equation}

\paragraph{The Hessian and Laplace operator}

We next consider the Hessians and Laplace operators for the two pseudo-metrics. For any function $ z $ on the manifold, by basic conformal calculation, we have
\begin{equation}\label{equation::thesis::conformal-hessian}
    \tilde{\nabla}^2 z - \nabla^2 z = -\dif f \otimes \dif z - \dif z \otimes \dif f + \Inner{\nabla f}{\nabla z}\mathfrak{h}.
\end{equation}

The Laplace operator is the trace of the Hessian. Taking traces gives
\[
    \tilde{\laplacian} z = e^{-2f} \left[ \laplacian z + (2n - 2) \Inner{\nabla f}{\nabla z} \right].
\]

Then we compare the induced Hessians and Laplace operators on the submanifold.

Denote the Levi-Civita connection on the submanifold $ \Sigma $ by $ \nabla^{\Sigma} $, and the corresponding Riemannian metric by $ g_{ij} = h_{ij} = \Inner{\vec{s}_i}{\vec{s}_j} $. Then $ \nabla_{\vec{X}} \vec{Y} = \nabla^{\Sigma}_{\vec{X}} \vec{Y} + \second_{\vec{X}\vec{Y}} $. Hence
\[
    (\nabla^{\Sigma})^2_{\vec{X}, \vec{Y}} z = \nabla^2_{\vec{X}, \vec{Y}} z + \second_{\vec{X}\vec{Y}}(z).
\]
The same identity holds for $ \tilde{\nabla} $. Taking traces gives
\[
    \begin{aligned}
        \laplacian_{\Sigma} z ={}         & h^{ij} \nabla^2_{\vec{s}_i, \vec{s}_j} z + \vec{H}(z),                       \\
        \tilde{\laplacian}_{\Sigma} z ={} & e^{-2f} h^{ij} \tilde{\nabla}^2_{\vec{s}_i, \vec{s}_j} z + \tilde{\vec{H}}(z).
    \end{aligned}
\]
Using (\ref{equation::thesis::conformal-hessian}) and (\ref{equation::thesis::conformal-mean-curvature-vector}), subtraction gives
\[
    e^{2f} \tilde{\laplacian}_{\Sigma} z - \laplacian_{\Sigma} z = -2\Inner{\nabla^{\Sigma} f}{\nabla^{\Sigma} z} + n\Inner{\nabla f}{\nabla z} - n\Inner{\left(\nabla f\right)^{\perp}}{\nabla z}.
\]
Since $ \Inner{\left(\nabla f\right)^{\perp}}{\nabla z} = \Inner{\left(\nabla f\right)^{\perp}}{\left(\nabla z\right)^{\perp}} = \Inner{\nabla f}{\nabla z} - \Inner{\nabla^{\Sigma} f}{\nabla^{\Sigma} z} $, we obtain
\begin{equation}\label{equation::thesis::difference-de-two-submanifold-laplacian}
    e^{2f} \tilde{\laplacian}_{\Sigma} z - \laplacian_{\Sigma} z = (n - 2) \Inner{\nabla^{\Sigma} f}{\nabla^{\Sigma} z}.
\end{equation}
Finally, using $ \tilde{\vec{H}} = \vec{0} $, which follows from maximality of the submanifold,
\begin{equation}\label{equation::thesis::expression-de-submanifold-laplacian-et-base-laplacian}
    \laplacian_{\Sigma} z = \tr_{h} \nabla^2 z + n \left[\Inner{\nabla z}{\nabla f} - \Inner{\nabla^{\Sigma}z}{\nabla^{\Sigma} f}\right].
\end{equation}

All terms associated with $ \mathfrak{h}^{\rho, \overline{\rho}} $ have now been expressed in terms of $ \mathfrak{h} $. The preceding formula suggests defining the weighted Laplace operator with respect to $ e^{nf} \dif \vol_{h} $ by
\[
    \laplacian_f z := \laplacian_{\Sigma} z + \Inner{\nabla^{\Sigma} (nf)}{\nabla^{\Sigma} z},
\]
Then
\begin{equation}\label{equation::thesis::expression-weighted-laplacian}
    \laplacian_f z = \tr_{h} \nabla^2 z + \Inner{\nabla z}{\nabla nf}.
\end{equation}

The weighted Laplace operator defined in this way is precisely $ \laplacian_f $ in \cite{collins2025boundary}. The preceding formula gives \cite[equation (3.3)]{collins2025boundary}. It is easy to verify that
\[
    \laplacian_f z = \frac{1}{\sqrt{\det h_{ij}} e^{nf}} \partial_i \left[ h^{ij} \sqrt{\det h_{ij}} e^{nf} \partial_{x^j} z(\point{x}, \point{T}(\point{x})) \right].
\]

\subsection{The trace of the Hessian on the submanifold}
We next analyze the first term on the right-hand side of (\ref{equation::thesis::expression-weighted-laplacian}).

In this subsection, let $ h_{\perp}^{\overline{k}\overline{l}} := -b_{\overline{k}a}h^{ab}b_{b\overline{l}} $ denote the inverse of the normal Gram matrix $ (\Inner{\vec{s}_{\overline{k}}}{\vec{s}_{\overline{l}}}) $. Since $ \vec{s}_i $ and $ \vec{s}_{\overline{k}} $ are orthogonal,
\[
    \laplacian z = \tr_{\mathfrak{h}} \nabla^2 z = h^{ij} \nabla^2_{\vec{s}_i, \vec{s}_j} z + h_{\perp}^{\overline{k} \overline{l}} \nabla^2_{\vec{s}_{\overline{k}}, \vec{s}_{\overline{l}}} z.
\]
Thus the Laplace operator on $ M \times \overline{M} $ decomposes along $ \Sigma $ into a tangential part and a normal part.

For the pseudo-Riemannian metric, only $ \Gamma^{i}_{jk} $ and $ \Gamma^{\overline{i}}_{\overline{j} \overline{k}} $ can be nonzero, so the covariant derivatives satisfy $ \nabla_{\vec{e}_i} \vec{e}_{\overline{k}} = \nabla_{\vec{e}_{\overline{k}}} \vec{e}_i = \vec{0} $. For convenience in the calculation, write $ \vec{f}_j := -h_{jt} b^{t \overline{l}} \vec{e}_{\overline{l}} $. Then $ \vec{s}_i = \vec{e}_i - \vec{f}_i $ and $ \vec{s}_{\overline{k}} = b^{\overline{k} i} \left( \vec{e}_i + \vec{f}_i \right) $. Therefore
\[
    \nabla_{\vec{e}_i} \vec{f}_j = -\vec{e}_i \left(h_{jt} b^{t \overline{l}}\right) \vec{e}_{\overline{l}} - h_{jt} b^{t \overline{l}} \nabla_{\vec{e}_i} \vec{e}_{\overline{l}} = -\vec{e}_i \left(h_{jt} b^{t \overline{l}}\right) \vec{e}_{\overline{l}},
\]
and
\[
    \nabla_{\vec{f}_j} \vec{e}_i = \vec{0}.
\]

Using the bilinearity of the Hessian as a tensor, we can simplify this to
\[
    \nabla^2_{\vec{s}_i, \vec{s}_j} z = \nabla^2_{\vec{e}_i, \vec{e}_j} z + \nabla^2_{\vec{f}_i, \vec{f}_j} z - \nabla^2_{\vec{e}_i, \vec{f}_j} z - \nabla^2_{\vec{f}_i, \vec{e}_j} z,
\]
and
\[
    \nabla^2_{\vec{s}_{\overline{k}}, \vec{s}_{\overline{l}}} z = b^{\overline{k} i} b^{\overline{l} j} \left[\nabla^2_{\vec{e}_i, \vec{e}_j} z + \nabla^2_{\vec{f}_i, \vec{f}_j} z + \nabla^2_{\vec{e}_i, \vec{f}_j} z + \nabla^2_{\vec{f}_i, \vec{e}_j} z\right].
\]

Thus, for the Laplacian of the function, taking traces gives
\[
    h^{ij} \nabla^2_{\vec{s}_i, \vec{s}_j} z = h^{ij} \left[ \nabla^2_{\vec{e}_i, \vec{e}_j} z + \nabla^2_{\vec{f}_i, \vec{f}_j} z - \nabla^2_{\vec{e}_i, \vec{f}_j} z - \nabla^2_{\vec{f}_i, \vec{e}_j} z \right],
\]
and
\[
    \begin{aligned}
        h_{\perp}^{\overline{k} \overline{l}} \nabla^2_{\vec{s}_{\overline{k}}, \vec{s}_{\overline{l}}} z ={} & -b_{\overline{k} a} h^{ab} b_{b \overline{l}} b^{\overline{k} i} b^{\overline{l} j} \left[\nabla^2_{\vec{e}_i, \vec{e}_j} z + \nabla^2_{\vec{f}_i, \vec{f}_j} z + \nabla^2_{\vec{e}_i, \vec{f}_j} z + \nabla^2_{\vec{f}_i, \vec{e}_j} z\right] \\
        ={}                                                                                                        & -h^{ij}\left[\nabla^2_{\vec{e}_i, \vec{e}_j} z + \nabla^2_{\vec{f}_i, \vec{f}_j} z + \nabla^2_{\vec{e}_i, \vec{f}_j} z + \nabla^2_{\vec{f}_i, \vec{e}_j} z \right].
    \end{aligned}
\]

Write
\begin{equation}\label{equation::thesis::formula-of-I-1}
    I_1 := h^{ij} \left[ \nabla^2_{\vec{e}_i, \vec{e}_j} z + \nabla^2_{\vec{f}_i, \vec{f}_j} z \right],
\end{equation}
and
\begin{equation}\label{equation::thesis::formula-of-I-2}
    I_2 := h^{ij} \left[ \nabla^2_{\vec{e}_i, \vec{f}_j} z + \nabla^2_{\vec{f}_i, \vec{e}_j} z \right],
\end{equation}
Then
\[
    h^{ij} \nabla^2_{\vec{s}_i, \vec{s}_j} z = I_1 - I_2, \quad h_{\perp}^{\overline{k} \overline{l}} \nabla^2_{\vec{s}_{\overline{k}}, \vec{s}_{\overline{l}}} z = -I_1 - I_2, \quad \laplacian z = -2I_2.
\]
We calculate explicit expressions for $ I_1 $ and $ I_2 $. First,
\[
    \begin{aligned}
        \nabla^2_{\vec{e}_i, \vec{e}_j} z ={} & z_{ij} - \Gamma^{k}_{ij} z_k = z_{ij} - b^{k \overline{m}} b_{\overline{m} ij} z_k,                                                                                                                                                                                                         \\
        \nabla^2_{\vec{f}_i, \vec{f}_j} z ={} & h_{is} b^{s \overline{k}} h_{jt} b^{t \overline{l}} \nabla^2_{\vec{e}_{\overline{k}}, \vec{e}_{\overline{l}}} z = h_{is} b^{s \overline{k}} h_{jt} b^{t \overline{l}} \left[ z_{\overline{k} \overline{l}} - b^{\overline{p} m} b_{m \overline{k} \overline{l}} z_{\overline{p}} \right].
    \end{aligned}
\]
Therefore,
\[
    I_1 = h^{ij} z_{ij} - b^{q \overline{m}} b_{ij \overline{m}} h^{ij} z_q + b^{i \overline{k}} h_{ij} b^{j \overline{l}} z_{\overline{k} \overline{l}} - b^{\overline{k} i} h_{ij} b^{j \overline{l}} b_{\overline{k} \overline{l} a} b^{a \overline{p}} z_{\overline{p}}, \quad I_2 = -2 b^{i \overline{k}} z_{i \overline{k}}.
\]

\subsubsection{The inner product term}
We now calculate the remaining term $ \Inner{\nabla nf}{\nabla z} $. First, the gradient with respect to the pseudo-metric $ \mathfrak{h} $ has the following expression in local coordinates:
\[
    \nabla z = 2b^{i \overline{l}} z_{\overline{l}} \vec{e}_i + 2b^{\overline{k} j} z_j \vec{e}_{\overline{k}},
\]
Note the coefficients here: the derivative with respect to $ \overline{x}^{\overline{k}} $ instead corresponds to the coefficient of $ \vec{e}_i $. Although length is not appropriate for a pseudo-Riemannian metric, we can still calculate the ``squared norm''
\[
    \left|\nabla z\right|^2 = \Inner{\nabla z}{\nabla z} = 4z_{\overline{l}} b^{\overline{l} j} z_j.
\]

By the definition of the gradient, substituting $ nf = \frac{1}{2} \left[ \log \rho(\point{x}) + \log \overline{\rho}(\overline{\point{x}}) - \log \left| \det \nabla^2_{\point{x} \overline{\point{x}}} b \right| \right] $ gives
\begin{equation}\label{equation::thesis::formula-de-inner-product}
    \Inner{\nabla nf}{\nabla z} = z_{\overline{l}} b^{\overline{l} i} \frac{\rho_i(\point{x})}{\rho(\point{x})} + z_{j} b^{j \overline{k}} \frac{\overline{\rho}_{\overline{k}}(\overline{\point{x}})}{\overline{\rho}(\overline{\point{x}})} - z_{\overline{l}} b^{\overline{l} i} b_{a \overline{b} i} b^{a \overline{b}} - z_{j} b^{j \overline{k}} b_{a \overline{b} \overline{k}} b^{a \overline{b}}.
\end{equation}

\subsection{Summary}
Collecting the preceding results,
\[
    \laplacian_f z = I_1 - I_2 + z_{\overline{l}} b^{\overline{l} i} \frac{\rho_i(\point{x})}{\rho(\point{x})} + z_{j} b^{j \overline{k}} \frac{\overline{\rho}_{\overline{k}}(\overline{\point{x}})}{\overline{\rho}(\overline{\point{x}})} - z_{\overline{l}} b^{\overline{l} i} b_{a \overline{b} i} b^{a \overline{b}} - z_{j} b^{j \overline{k}} b_{a \overline{b} \overline{k}} b^{a \overline{b}}.
\]
We aim to choose $ z $ so that $ I_2 = -2n $, as in \cite{collins2025boundary}, or, if this is not possible, so that the error can be controlled.

\section{Step 2 for the monotonicity formula: choice of the function $ z $}\label{sec::thesis::insert-c}
In \cite{collins2025boundary}, the extrinsic distance $ \point{x} \cdot \nabla u(\point{x}) $ is considered. This is precisely the restriction of $ z = \Inner{\point{x}}{\overline{\point{x}}} $ to $ \Sigma $.

Formally, there are two ways to generalize this function. The first retains the ``extrinsic distance'', namely the ``distance'' function defined by geodesics on $ M \times \overline{M} $; the second uses the benefit function $ b(\point{x}, \overline{\point{x}}) $. In the quadratic-distance case, the two choices give the same function, whereas in our setting each has advantages and disadvantages. Choosing the distance function agrees more closely with the origin of the monotonicity formula: recall the example in \cite{gilbarg1998elliptic}, which considers an integral over the part of a minimal surface cut out by a ball in the ambient space. However, \citeauthor{collins2025boundary} also use positivity of the ``distance'', which is not apparent on a pseudo-Riemannian manifold. Choosing the benefit function makes positivity a consequence of optimality of the transport plan, but the benefit function presents some sign difficulties. The author initially chose the benefit function because it was easier to calculate with, but eventually abandoned this choice because it was unclear how to control the density integral. We therefore turn to the ``extrinsic distance''.

More precisely, fix a base point $ \point{o} := (\point{x}_0, \overline{\point{x}}_0) $ on $ \Sigma $ and a normal neighborhood $ U $ of $ \point{o} $ for the ambient metric $ \mathfrak{h} $. For every $ \point{p} \in U $, write $ \vec{v} := \exp_{\point{o}}^{-1}(\point{p}) $ and let $ \gamma(t) := \exp_{\point{o}}(t\vec{v}) $, $ 0 \leqslant t \leqslant 1 $, be the radial geodesic from $ \point{o} $ to $ \point{p} $. Define
\[
    \sigma(\point{p}) := \int_{0}^{1} \mathfrak{h}(\dot{\gamma}(t), \dot{\gamma}(t)) \dif t.
\]

We now calculate the terms in $ \laplacian_f z $ for $ z = \sigma $.

\subsection{First and second derivatives}

In this section we calculate the first and second derivatives of $ \sigma $.

\paragraph{First derivatives} Take a tangent vector $ \vec{X} \in T_{\point{p}} (M \times \overline{M}) $. Construct a variation $ \alpha(t, s) \colon [0, 1] \times (-\varepsilon, \varepsilon) \to M \times \overline{M} $ through geodesics $ t \mapsto \alpha(t,s) $, with $ \alpha(t,0) = \gamma(t) $ and $ \alpha(0,s) = \point{o} $. Write $ \point{T} := \alpha_{\ast}\partial_t $ for the tangent vector field, which equals $ \dot{\gamma} $ when $ s = 0 $ and $ \vec{J} = \alpha_{\ast} \partial_s $ for the variation vector field, with $ \vec{J}(1) = \vec{X} $. Then $ \vec{J}(0) = \vec{0} $ and $ [\point{T}, \vec{J}] = \vec{0} $. By the first variation formula,
\[
    \begin{aligned}
        \vec{X}(\sigma) ={}& \int_{0}^{1} 2\mathfrak{h}(\nabla_{\vec{J}} \point{T}, \point{T}) \dif t \\
        ={}& \int_{0}^{1} 2\mathfrak{h}(\nabla_{\point{T}} \vec{J}, \point{T}) \dif t \\
        ={}& 2\mathfrak{h}(\vec{X}, \point{T}(1)),
    \end{aligned}
\]
so $ \nabla \sigma = 2\dot{\gamma}(1) $.

Since $ \nabla \sigma = 2 b^{i \overline{l}} \sigma_{\overline{l}} \vec{e}_i + 2 b^{\overline{k} j} \sigma_j \vec{e}_{\overline{k}} = 2\dot{x}^i \vec{e}_i + 2\dot{\overline{x}}^{\overline{k}} \vec{e}_{\overline{k}} $, we also have $ \sigma_{\overline{k}} = b_{i \overline{k}}\dot{x}^i $ and $ \sigma_{i} = b_{\overline{k} i} \dot{\overline{x}}^{\overline{k}} $.

\paragraph{Second derivatives} For the second derivatives, the definition of the Hessian gives
\[
    \begin{aligned}
        \nabla^2_{\vec{X}\vec{Y}} \sigma ={}& \Inner{\nabla_{\vec{X}} \nabla \sigma}{\vec{Y}} \\
        ={}& 2\Inner{\nabla_{\point{T}} \vec{J}}{\vec{Y}} \\
    \end{aligned}
\]
Thus the key to calculating the Hessian is to calculate derivatives of Jacobi fields.

\paragraph{A parallel frame}

Choose a frame $ \{\partial_{x_1}, \partial_{x_2}, \cdots, \partial_{x_n}; \partial_{\overline{x}_1}, \partial_{\overline{x}_2}, \cdots, \partial_{\overline{x}_n}\} $ at the initial point $ \point{o} $, and parallel transport it along $ \gamma(t) $ to obtain a frame $ \{ \vec{E}_i(t); \vec{E}_{\overline{k}}(t) \} $ along the geodesic. The parallel transport equation is $ \nabla_{\dot{\gamma}} \vec{E}_i \equiv 0 $. Write
\[
    \vec{E}_i(t) = E_i^{a}(t) \partial_{x_a} + E_i^{\overline{b}}(t) \partial_{\overline{x}_{\overline{b}}},
\]
Using capital letters $ A $, $ B $, and so on for indices on either the $ M $ or $ \overline{M} $ side, we have
\[
    \frac{\dif}{\dif t}E_i^{a}(t) + \Gamma^{a}_{AB} \dot{\gamma}^{A}(t) E_i^{B}(t) = 0,
\]
Since the Christoffel symbols are nonzero only when their upper and lower indices are all on the same side, the preceding equation is equivalent to
\[
    \frac{\dif}{\dif t}E_i^{a}(t) + \Gamma^{a}_{pq} \dot{\gamma}^{p}(t) E_i^{q}(t) = 0.
\]
Similarly, on the other side,
\[
    \frac{\dif}{\dif t}E_i^{\overline{b}}(t) + \Gamma^{\overline{b}}_{\overline{p}\overline{q}} \dot{\gamma}^{\overline{p}}(t) E_i^{\overline{q}}(t) = 0.
\]
Since $ \vec{E}_i(0) = \partial_{x_i} $ and $ E_i^{\overline{b}}(0) = 0 $, the equation for its $ \overline{M} $ components has the zero solution. Uniqueness for linear ordinary differential equations implies that all its $ \overline{M} $ components vanish. Thus parallel transport preserves vectors on the $ M $ or $ \overline{M} $ side. The resulting $ \vec{E}_{A}(t) $ is therefore a frame along the geodesic that preserves the splitting into the two sides.

Since the tangent vector is parallel along the geodesic, there are constants $ T^A $ such that $ \point{T}(t) = T^i \vec{E}_i(t) + T^{\overline{k}} \vec{E}_{\overline{k}}(t) $. Write the metric matrix as
\[
    \matr{H}(t) := \left[ \mathfrak{h}(\vec{E}_A(t), \vec{E}_B(t)) \right]_{AB} =: \begin{bmatrix}
        \matr{0} & \matr{G} \\ \matr{G}^{\top} & \matr{0}
    \end{bmatrix},
\]
Then parallel transport implies that $ \matr{H} $ is constant in $ t $.

\paragraph{The coefficient matrix of Jacobi fields}

Decomposing the Jacobi field $ \vec{J} $ in the frame $ \vec{E}_{A}(t) $ gives $ \vec{J}(t) = J^{A}(t) \vec{E}_{A}(t) $. In the Jacobi equation $ \nabla_{\point{T}} \nabla_{\point{T}} \vec{J} = \Rm (\point{T}, \vec{J})\point{T} $, the covariant derivative of $ \vec{E}_{A} $ in the direction $ \point{T} $ vanishes, so the equation becomes
\[
    \ddot{J}^{F}(t) \vec{E}_{F}(t) = \Rm (\point{T}, \vec{J}) \point{T} = -{\Rm^{F}}_{BCD} J^{C} T^{D} T^{B} \vec{E}_{F}(t).
\]
Taking the inner product of both sides with $ \vec{E}_{A} $ gives
\[
    \ddot{J}^{F}(t) h_{FA} = -\Rm_{ABCD} T^{B} J^{C} T^{D}.
\]
By \ref{lemma::thesis::Christoffel-symbol-and-curvature}, the curvature tensor can be nonzero only when each of its two pairs of indices contains one index from each side. The metric matrix $ \matr{H} $ also has nonzero entries only between the two sides. Thus the preceding equation can be written as
\begin{equation}\label{equation::thesis::jacobi-equation-for-single-jacobi-field}
    \begin{aligned}
        \ddot{J}^{i}(t) h_{i \overline{k}}
        ={}& -\Rm_{\overline{k} i \overline{l} j} T^{i} J^{\overline{l}} T^{j}
             -\Rm_{\overline{k} i j \overline{l}} T^{i} J^{j} T^{\overline{l}}, \\
        \ddot{J}^{\overline{k}}(t) h_{\overline{k} i}
        ={}& -\Rm_{i \overline{k} j \overline{l}} T^{\overline{k}} J^{j} T^{\overline{l}}
             -\Rm_{i \overline{k} \overline{l} j} T^{\overline{k}} J^{\overline{l}} T^{j}.
    \end{aligned}
\end{equation}

The coefficients of Jacobi fields solve a second-order ordinary differential equation. Moreover, the Jacobi fields considered here arise from variations fixing the initial point $ \point{o} $, so $ \vec{J}(0) = \vec{0} $. These Jacobi fields therefore form an $ N := 2n $ dimensional linear space, and can be represented by their derivatives at the initial point and a matrix $ \matr{A}(t) \in \mathbb{R}^{2n \times 2n} $.

More precisely, $ \matr{A}(t)(v^A\vec{E}_{A}) $ represents the Jacobi field whose derivative at $ \point{o} $ is $ \nabla_{\point{T}} \vec{J}(0) = v^A\vec{e}_A $. Therefore
\[
    \matr{A}(0) = \matr{0}, \qquad \dot{\matr{A}}(0) = \matr{I},
\]
and it satisfies the Jacobi equation. Substituting \ref{equation::thesis::jacobi-equation-for-single-jacobi-field} and writing $ \matr{R}_{AC} := \Rm_{ABCD}T^{B}T^{D} $, we obtain
\[
    \matr{H} \ddot{\matr{A}}(t) = -\matr{R} \matr{A}(t).
\]

We return to the calculation of $ \nabla^2 \sigma $. For $ s > 0 $ in the normal neighborhood, identify vectors with their coefficients in the parallel frame. The Jacobi field with $ \vec{J}(0) = \vec{0} $ and $ \vec{J}(s) = \vec{X} $ is
\[
    \vec{J}(t) = \matr{A}(t)\matr{A}(s)^{-1}\vec{X}, \qquad
    \nabla_{\point{T}(t)}\vec{J}(t) = \dot{\matr{A}}(t)\matr{A}(s)^{-1}\vec{X}.
\]
Write
\[
    \matr{L}(t) := \dot{\matr{A}}(t)\matr{A}^{-1}(t).
\]
The first variation formula, applied to the radial geodesic reparametrized on $ [0,1] $, gives $ \nabla\sigma(\gamma(s)) = 2s\point{T}(s) $. Differentiating this identity along a geodesic variation with endpoint vector $ \vec{X} $ gives
\[
    \nabla^2_{\vec{X}\vec{Y}}\sigma\big|_{\gamma(s)}
    = 2s\Inner{\nabla_{\point{T}}\vec{J}(s)}{\vec{Y}}
    = 2s\vec{Y}^{\top}\matr{H}\matr{L}(s)\vec{X},
\]
and therefore
\begin{equation}\label{equation::thesis::formula-of-hessian-sigma-in-L}
    \nabla^2 \sigma = 2t \matr{H}\matr{L}(t).
\end{equation}

Moreover, the equation satisfied by $ \matr{L}(t) $ is
\begin{equation}\label{equation::thesis::equation-of-L}
    \matr{H} \matr{L}^\prime + \matr{H}\matr{L}^2 = -\matr{R}.
\end{equation}

Write $ \matr{H} \matr{L} $ in block form as
\[
    \matr{H} \matr{L} =: \begin{bmatrix}
        \matr{E} & \matr{B} \\
        \matr{B}^{\top} & \matr{D}
    \end{bmatrix},
\]
The traces of $ \matr{E} $, $ \matr{B} $, and $ \matr{D} $ are precisely what we need when calculating $ I_1 $ and $ I_2 $. Since $ \sigma $ has no explicit algebraic expression, we control $ I_1 $, $ I_2 $, and the other terms by studying the equations satisfied by these blocks.

\subsection{Analysis of \texorpdfstring{$I_1$}{I1}}\label{subsection::thesis::analyse-of-I_1}

Taking $ \matr{E} $ as an example, (\ref{equation::thesis::equation-of-L}) gives the equation for $ \matr{E} $
\begin{equation}\label{equation::thesis::equation-of-E}
    \matr{E}^{\prime} + \matr{E} \matr{N}^{\top} + \matr{N} \matr{E} = -\matr{R}_{xx},
\end{equation}
where $ \matr{N} := \matr{B}\matr{G}^{-1} $, and $ \matr{R}_{xx} $ denotes the upper-left block of $ \matr{R} $, so $ (\matr{R}_{xx})_{ij} = \Rm_{i \overline{k} j \overline{l}} T^{\overline{k}} T^{\overline{l}} $. For convenience, write $ \point{T} = (\vec{a}, \vec{b}) $, so $ (\matr{R}_{xx})_{ij} = \Rm_{i \overline{k} j \overline{l}} b^{\overline{k}} b^{\overline{l}} $. The Jacobi field parallel to $ \point{T} $ is $ \vec{J}(t) = t \point{T}(t) $, so $ \matr{A}(t) (\vec{a}, \vec{b})^{\top} = t(\vec{a}, \vec{b})^{\top} $ and $ \matr{A}^\prime(t) (\vec{a}, \vec{b})^{\top} = (\vec{a}, \vec{b})^{\top} $. Therefore
\begin{equation}\label{equation::thesis::equation-of-Ea}
    \begin{aligned}
        \matr{E} \vec{a} + \matr{B} \vec{b} ={}& \frac{1}{t} \matr{G}\vec{b}, \\
        \matr{B}^{\top} \vec{a} + \matr{D} \vec{b} ={}& \frac{1}{t} \matr{G}^\top\vec{a}.
    \end{aligned}
\end{equation}

Applying the fundamental theorem of calculus to $ \matr{A}(t) $ gives
\[
    \begin{aligned}
        \matr{A}(t) ={}& \int_0^{t} \dot{\matr{A}}(s) \dif s \Rightarrow \norm{\matr{A}(t)} \leqslant M t; \\
        \dot{\matr{A}}(t) ={}& \matr{I}_{2n} - \int_{0}^{t} \matr{H}^{-1} \matr{R}(s) \matr{A}(s) \dif s = \matr{I}_{2n} + O(t^2); \\
        \matr{A}(t) ={}& \int_{0}^{t} \dot{\matr{A}}(s) \dif s = t \matr{I}_{2n} + O(t^3),
    \end{aligned}
\]
Therefore
\[
    \matr{L}(t) = \frac{1}{t} \matr{I}_{2n} + O(t).
\]
Since $ \matr{H}\matr{L} = t^{-1}\matr{H} + O(t) $, we have $ \matr{B} = t^{-1}\matr{G} + O(t) $, so $ \matr{N} $ satisfies $ \matr{N} = \frac{1}{t} \matr{I}_{n} + O(t) $. Thus equation \ref{equation::thesis::equation-of-E} is approximately
\[
    \matr{E}^\prime + \frac{2}{t} \matr{E} = -\matr{R}_{xx},
\]
For this approximate equation, $ (t^2 \matr{E})^\prime = -t^2 \matr{R}_{xx} $, suggesting $ \matr{E} \leqslant 0 $ when $ \matr{R}_{xx} \geqslant 0 $. However, $ \matr{R}_{xx} = (\matr{R}_{xx})_{ij} = \Rm_{i \overline{k} j \overline{l}} T^{\overline{k}} T^{\overline{l}} $, and conditions (\ref{condition::thesis::A3w}) or (\ref{condition::thesis::A3s}) only give, when $ \vec{v}^{\top} \matr{G} \vec{b} = 0 $,
\[
    \Rm_{i \overline{k} j \overline{l}} v^i v^j T^{\overline{k}} T^{\overline{l}} \geqslant \begin{cases}
        0, & (\ref{condition::thesis::A3w}), \\
        C_0 \norm{\vec{v}}^2 \norm{\vec{b}}^2, & (\ref{condition::thesis::A3s}).
    \end{cases}
\]
If $ \vec{b} = \vec{0} $, then $ \matr{R}_{xx} = \matr{0} $ and the desired lower bound is immediate. Otherwise, $ \vec{l} := \matr{G}\vec{b} \ne \vec{0} $. Decompose any direction $ \vec{w} $ into components perpendicular and parallel to $ \vec{l} $, namely $ \vec{w} = \vec{w}_{\perp} + \vec{w}_{\parallel} $. Then
\[
    \matr{R}_{xx}(\vec{w}, \vec{w}) = \matr{R}_{xx}(\vec{w}_{\perp}, \vec{w}_{\perp}) + 2\matr{R}_{xx}(\vec{w}_{\perp}, \vec{w}_{\parallel}) + \matr{R}_{xx}(\vec{w}_{\parallel}, \vec{w}_{\parallel}),
\]
Since $ \Rm_{i \overline{k} j \overline{l}} = \frac{1}{2} \left( b_{ij \overline{k} \overline{l}} - b_{ij \overline{a}} b^{\overline{a} m} b_{\overline{k} \overline{l} m} \right) $, the curvature tensor is controlled by the norm of $ b \in C^{4} $. Since the geodesic $ \gamma(t) $ need not lie entirely in $ \cl \regs \times \cl \regt $, assume there is an open set $ N \supset \cl \regs \times \cl \regt $ on which conditions (\ref{condition::thesis::A0}) - (\ref{condition::thesis::A3s}) or (\ref{condition::thesis::A3w}) hold on $ N $. Then, using the local geometric neighborhoods described in \cref{remark::thesis::small-height-and-geometric-neighborhood}, there is $ C > 0 $ such that $ |\matr{R}_{xx}(\vec{v}, \vec{v})| \leqslant C \norm{\vec{v}}_{\mathbb{R}^n}^2 \norm{\vec{b}}_{\mathbb{R}^n}^2 $.

Now assume (\ref{condition::thesis::A3s}) holds. Then there is $ \lambda_{\mathrm{sec}} > 0 $ such that
\[
    \matr{R}_{xx}(\vec{w}_{\perp}, \vec{w}_{\perp}) \geqslant \lambda_{\mathrm{sec}} \norm{\vec{w}_{\perp}}_{\mathbb{R}^n}^2 \norm{\vec{b}}_{\mathbb{R}^n}^2,
\]
Hence
\[
    \begin{aligned}
        \matr{R}_{xx}(\vec{w}, \vec{w}) ={}& \matr{R}_{xx}(\vec{w}_{\perp}, \vec{w}_{\perp}) + 2\matr{R}_{xx}(\vec{w}_{\perp}, \vec{w}_{\parallel}) + \matr{R}_{xx}(\vec{w}_{\parallel}, \vec{w}_{\parallel}) \\
        \geqslant{}& \norm{\vec{b}}_{\mathbb{R}^n}^2 \left[ \lambda_{\mathrm{sec}} \norm{\vec{w}_{\perp}}_{\mathbb{R}^n}^2 - 2C \norm{\vec{w}_{\perp}}_{\mathbb{R}^n} \norm{\vec{w}_{\parallel}}_{\mathbb{R}^n} - C \norm{\vec{w}_{\parallel}}_{\mathbb{R}^n}^2 \right] \\
        ={}& \norm{\vec{b}}_{\mathbb{R}^n}^2 \left[ \lambda_{\mathrm{sec}} \left( \norm{\vec{w}_{\perp}}_{\mathbb{R}^n} - \frac{C}{\lambda_{\mathrm{sec}}}\norm{\vec{w}_{\parallel}}_{\mathbb{R}^n} \right)^2 - \left( C + \frac{C^2}{\lambda_{\mathrm{sec}}} \right) \norm{\vec{w}_{\parallel}}^2 \right] \\
        \geqslant{}& - \left( C + \frac{C^2}{\lambda_{\mathrm{sec}}} \right) \frac{\norm{\vec{b}}_{\mathbb{R}^n}^2}{\norm{\vec{l}}_{\mathbb{R}^n}^2} (\vec{w}^{\top} \vec{l})^2 \\
        \geqslant{}& -\frac{C + C^2 / \lambda_{\mathrm{sec}}}{m^2} (\vec{w}^{\top} \vec{l})^2.
    \end{aligned}
\]
Here $ m > 0 $ is the smallest singular value of the invertible matrix $ \matr{G} $, so that $ \norm{\vec{l}}_{\mathbb{R}^n} \geqslant m \norm{\vec{b}} $. Therefore there is a constant $ K(C, \lambda_{\mathrm{sec}}, m) > 0 $ such that
\begin{equation}\label{equation::thesis::lowerbound-of-Rm}
    \matr{R}_{xx} \geqslant -2K \vec{l} \vec{l}^{\top}.
\end{equation}

\begin{remark}
    The factor $ 2 $ is included only to simplify the coefficients below.

    The purpose of this step is to obtain an estimate of the form $ \vec{l} \vec{l}^{\top} $, and the constant $ \lambda_{\mathrm{sec}} $ supplied by (\ref{condition::thesis::A3s}) is essential here. Without this positive lower bound, as under condition (\ref{condition::thesis::A3w}), a lower bound of this form cannot be obtained.

    If condition (\textrm{NNCC}) holds on $ N $, then $ \matr{R}_{xx} \geqslant 0 $, and the matrix comparison argument below, with $ \kappa = 0 $, gives $ \matr{E} \leqslant 0 $.
\end{remark}

Consider $ \matr{F} := \matr{E} - \kappa t \vec{l} \vec{l}^{\top} $. Since $ \matr{G} $ and the components of $ \point{T} $ in the parallel frame are constant, $ \vec{l} $ is constant. The equation satisfied by $ \matr{F} $ is
\[
    \matr{F}^\prime + \kappa \vec{l} \vec{l}^\top + \matr{F} \matr{N}^{\top} + \matr{N} \matr{F} + \kappa t \vec{l} \vec{l}^{\top} \matr{N}^{\top} + \kappa t \matr{N} \vec{l} \vec{l}^{\top} = -\matr{R}_{xx}.
\]
Using \cref{equation::thesis::equation-of-Ea}, together with $ \matr{N} = \matr{B}\matr{G}^{-1} $ and $ \vec{l} = \matr{G}\vec{b} $, gives
\[
    \matr{N}\vec{l} = \matr{B}\vec{b} = \frac{1}{t}\vec{l} - \matr{E}\vec{a}
    = \frac{1}{t}\vec{l} - \matr{F}\vec{a} - \kappa t\vec{l}\vec{l}^{\top}\vec{a}.
\]
Write $ e := \vec{l}^{\top}\vec{a} = \frac{1}{2}\mathfrak{h}(\point{T},\point{T}) $. Since $ \matr{F} $ is symmetric, its equation becomes
\[
    \matr{F}^\prime + \matr{F}\matr{N}^{\top} + \matr{N}\matr{F}
    - \kappa t\matr{F}\vec{a}\vec{l}^{\top}
    - \kappa t\vec{l}\vec{a}^{\top}\matr{F}
    + (3\kappa - 2\kappa^2 e t^2)\vec{l}\vec{l}^{\top}
    = -\matr{R}_{xx}.
\]
Set $ \matr{N}_1 := \matr{N} - \kappa t\vec{l}\vec{a}^{\top} $. By \cref{equation::thesis::lowerbound-of-Rm},
\[
    \begin{aligned}
        \matr{F}^{\prime} + \matr{F}\matr{N}_1^{\top} + \matr{N}_1\matr{F}
        ={}& -\matr{R}_{xx} - 3\kappa\vec{l}\vec{l}^{\top}
             + 2\kappa^2 e t^2\vec{l}\vec{l}^{\top} \\
        \leqslant{}& (2K - 3\kappa + 2\kappa^2 e t^2)\vec{l}\vec{l}^{\top}.
    \end{aligned}
\]
Since $ \matr{N}_1 = t^{-1}\matr{I}_n + O(t) $, write $ \matr{N}_1 = t^{-1}\matr{I}_n + \matr{N}_2(t) $. Let $ \matr{W} $ solve $ \matr{W}^{\prime} = \matr{W}\matr{N}_2 $ with $ \matr{W}(0) = \matr{I}_n $. The matrix $ \matr{W} $ is invertible, and direct differentiation gives
\[
    (t^2\matr{W}\matr{F}\matr{W}^{\top})^{\prime}
    = t^2\matr{W}\left(
        \matr{F}^{\prime} + \matr{F}\matr{N}_1^{\top} + \matr{N}_1\matr{F}
      \right)\matr{W}^{\top}.
\]
Take $ \kappa = 2K $. Since
\[
    \sigma(\gamma(t)) = t^2\mathfrak{h}(\point{T},\point{T}) = 2et^2,
\]
we obtain, whenever $ K\sigma(\gamma(t)) \leqslant 1 $,
\[
    (t^2\matr{W}\matr{F}\matr{W}^{\top})^{\prime}
    \leqslant (-2\kappa + \kappa^2\sigma(\gamma(t)))t^2
        (\matr{W}\vec{l})(\matr{W}\vec{l})^{\top}
    \leqslant 0.
\]
If this condition holds at an endpoint of the radial segment, it holds at every preceding point because $ \sigma(\gamma(t)) = 2et^2 $. Moreover, $ \matr{H}\matr{L} = t^{-1}\matr{H} + O(t) $ and the upper-left block of $ \matr{H} $ vanishes, so $ \matr{E} = O(t) $ and $ \matr{F} = O(t) $. Consequently,
\[
    \lim_{t\to 0} t^2\matr{W}(t)\matr{F}(t)\matr{W}(t)^{\top} = \matr{0}.
\]
Integrating the preceding inequality and using invertibility of $ \matr{W} $ gives $ \matr{F} \leqslant 0 $, namely
\[
    \matr{E} \leqslant 2Kt\vec{l}\vec{l}^{\top}.
\]
Returning to $ \nabla^2\sigma $, its source-source block in the parallel frame is $ 2t\matr{E} $. In this frame,
\[
    d_{\point{x}}\sigma(\vec{p}\oplus\overline{\vec{0}})
    = \Inner{2t\point{T}}{\vec{p}\oplus\overline{\vec{0}}}
    = 2t\vec{l}^{\top}\vec{p}.
\]
Thus $ 2t\matr{E} \leqslant K(2t\vec{l})(2t\vec{l})^{\top} $. Parallel transport preserves the source-target splitting, so this quadratic-form inequality also holds in the coordinate frame:
\[
    (\nabla^2_{\vec{e}_i,\vec{e}_j}\sigma)_{ij}
    \leqslant K(\sigma_i\sigma_j)_{ij}.
\]
Interchanging source and target gives the same argument for $ \matr{D} $. Taking $ K $ large enough for both curvature bounds, we obtain, on the region $ K\sigma \leqslant 1 $,
\[
    (\nabla^2_{\vec{e}_{\overline{k}},\vec{e}_{\overline{l}}}\sigma)_{\overline{k}\overline{l}}
    \leqslant K(\sigma_{\overline{k}}\sigma_{\overline{l}})_{\overline{k}\overline{l}}.
\]

In the monotonicity formula we want $ I_1 \leqslant 0 $, so we may take
\[
    \tau = \frac{1 - e^{-K \sigma}}{K}, 
\]
If $ K = 0 $ (for example, in the \textrm{NNCC} case), take $ \tau = \sigma $ directly. Then
\[
    \nabla^2 \tau = e^{-K \sigma} \left( \nabla^{2} \sigma - K \dif \sigma \otimes \dif \sigma \right), 
\]
Therefore
\[
    I_1(\tau) \leqslant 0.
\]

\subsection{Analysis of \texorpdfstring{$I_2$}{I2}}

Having defined $ \tau $, we now consider $ I_2 $. At $ (\point{x},\overline{\point{x}})=\gamma(t) $, for $ \sigma $,
\[
    I_2(\sigma) = -\frac{1}{2} \tr_{\mathfrak{h}} \nabla^{2} \sigma = -t \tr \matr{L}(t) =: -2n + \epsilon(\point{x},\overline{\point{x}}),
\]
Here $ \epsilon(\gamma(t)) = -t\tr\matr{L}(t) + 2n = O(t^2) $
as $ t \to 0 $ along the fixed geodesic.

Hence
\[
    I_2(\tau) = e^{-K \sigma} \left( -2n + \epsilon(\gamma(t)) + \frac{K}{2}\norm{\nabla \sigma}^2_{\mathfrak{h}} \right) = e^{-K \sigma} \left( -2n + \epsilon(\point{x},\overline{\point{x}}) + 2K \sigma \right).
\]

For now, retain the term $ \epsilon(\point{x},\overline{\point{x}}) $ and use the resulting $ 2n $ to derive the monotonicity formula.

\subsection{Integration by parts}

For the chosen $ \tau $, we obtain
\[
    \laplacian_f \tau = I_1 - e^{-K \sigma} \left( -2n + \epsilon(\point{x},\overline{\point{x}}) + 2K \sigma \right) + \Inner{\nabla \tau}{\nabla nf},
\]
so
\begin{equation}\label{equation::thesis::expression-of-one}
    2n = \laplacian_f \tau - I_1 - 2n\left( e^{-K \sigma} - 1 \right) + e^{-K\sigma} \left( \epsilon(\point{x},\overline{\point{x}}) + 2K \sigma \right) - \Inner{\nabla \tau}{\nabla nf}.
\end{equation}

First consider the weighted volume element $ e^{nf} \dif \vol_{h} $. Substituting the expression for $ f $ gives
\[
    e^{nf} = \sqrt{\frac{\rho(\point{x}) \overline{\rho}(\overline{\point{x}})}{\left|\det \nabla^2_{\point{x}\overline{\point{x}}} b\right|}},
\]
while
\[
    \dif \vol_{h} = \sqrt{\det h_{ij}} \dif \point{x} = \sqrt{\frac{\rho(\point{x})}{\overline{\rho}(\overline{\point{x}})} \left|\det \nabla^2_{\point{x}\overline{\point{x}}} b\right|} \dif \point{x}.
\]
Therefore
\[
    e^{nf} \dif \vol_{h} = \rho(\point{x}) \dif \point{x}.
\]

\paragraph{The main part}
Take $ \kappa > 0 $ and consider the function analogous to that in \cite{collins2025boundary}:
\[
    \frac{\dif}{\dif r} \left[ r^{-\kappa} \int_{\{ \tau \leqslant r \}} e^{nf} \dif \vol_{h} \right] = -\kappa r^{-(\kappa + 1)} \int_{\{ \tau \leqslant r \}} \rho(\point{x}) \dif \point{x} + r^{-\kappa} \int_{\{ \tau = r \} \cap \regs} \frac{\rho(\point{x})}{\norm{\nabla^{\mathbb{R}^n} \tau}} \dif \mathcal{H}^{n - 1}(\point{x}).
\]
The last term uses the coarea formula, and all integrals over level sets are restricted to $ \regs $. For the following calculations, first assume that the functions involved are smooth.
Using (\ref{equation::thesis::expression-of-one}),
\begin{equation}\label{equation::thesis::inserting-laplacian-f-2n-others}
    \begin{aligned}
        \int_{\{ \tau \leqslant r \}} \rho(\point{x}) \dif \point{x} ={} & \frac{1}{2n} \int_{\{ \tau \leqslant r \}} \rho \laplacian_f \tau \dif \point{x} - \frac{1}{2n} \int_{\{ \tau \leqslant r \}} I_1 \rho \dif \point{x} \\
        {} & - \frac{1}{2n} \int_{\{ \tau \leqslant r \}} n \Inner{\nabla f}{\nabla \tau} \rho \dif \point{x} \\
        {} & + \frac{1}{2n} \int_{\{ \tau \leqslant r \}} \left[ 2n \left( 1 - e^{-K \sigma} \right) + e^{-K \sigma} \left( \epsilon(\point{x},\overline{\point{x}}) + 2K \sigma \right) \right] \rho \dif \point{x}.
    \end{aligned}
\end{equation}

Next, integration by parts on $ \mathbb{R}^n $ gives the following expression for the $ \laplacian_f $ term:
\[
    \begin{aligned}
        \frac{1}{2n} \int_{\{ \tau \leqslant r \}} \rho \laplacian_f \tau \dif \point{x} ={} & \frac{1}{2n} \int_{\{ \tau \leqslant r \}} \frac{1}{\sqrt{\det h_{ij}} e^{nf}} \partial_i \left[ h^{ij} \sqrt{\det h_{ij}} e^{nf} \partial_{x^j} \tau(\point{x}, \point{T}(\point{x})) \right] \rho \dif \point{x} \\
        ={} & \frac{1}{2n} \int_{\{ \tau \leqslant r \}} \partial_i \left[ h^{ij} \sqrt{\det h_{ij}} e^{nf} \partial_{x^j} \tau(\point{x}, \point{T}(\point{x})) \right] \dif \point{x} \\
        ={} & \frac{1}{2n} \int_{\{ \tau = r \} \cap \regs} \Inner{\nabla^{\Sigma} \tau}{\nabla^{\Sigma} \tau}_{h} \frac{\rho(\point{x})}{\norm{\nabla^{\mathbb{R}^n} \tau}} \dif \mathcal{H}^{n - 1}(\point{x}) \\
        {} & + \frac{1}{2n} \int_{\partial \regs \cap \{ \tau \leqslant r \}} h^{ij} \partial_{x^j} \tau(\point{x}, \point{T}(\point{x})) (\vec{\nu}_{\regs})_i \rho \dif \mathcal{H}^{n - 1}(\point{x}).
    \end{aligned}
\]
Note that here $ \tau = \tau(\point{x}, \point{T}(\point{x})) $ is regarded as a function on $ M $, so $ \partial_{x^j} \tau(\point{x}, \point{T}(\point{x})) $ does not differentiate only on the $ M $ side. Also, $ \vec{\nu}_{\regs} $ is the outward Euclidean unit normal.

\paragraph{The gradient term on the submanifold}
Next consider $ \Inner{\nabla^{\Sigma} \tau}{\nabla^{\Sigma} \tau}_{h} $. Below, $ \tau_i, \tau_{\overline{k}} $ denote partial derivatives in the ambient space. First recall
\[
    \nabla \tau = 2b^{i \overline{l}} \tau_{\overline{l}} \vec{e}_i + 2b^{\overline{k} j} \tau_j \vec{e}_{\overline{k}}.
\]
and
\[
    |\nabla \tau|^2 = 4 \tau_i b^{i \overline{k}} \tau_{\overline{k}} = 4 \sigma e^{-2K \sigma}.
\]
Here $ |\nabla \tau|^2 $ denotes the inner product under the ambient pseudo-Riemannian metric.
We next calculate $ \nabla^{\Sigma} \tau $. Using $ \vec{s}_i, \vec{s}_{\overline{k}} $, decompose $ \nabla \tau $ into two mutually ``orthogonal'' parts:
\[
    \nabla \tau = \Inner{\nabla \tau}{\vec{s}_i} h^{ij} \vec{s}_j + \Inner{\nabla \tau}{\vec{s}_{\overline{k}}} h^{\overline{k} \overline{l}} \vec{s}_{\overline{l}}.
\]
The first part is $ \nabla^{\Sigma} \tau $, and the second is the normal component. Recall that $ \vec{s}_i = \vec{e}_i - \vec{f}_i $, $ \vec{s}_{\overline{k}} = b^{\overline{k} i} \left( \vec{e}_i + \vec{f}_i \right) $, and $ \vec{f}_i = -h_{is} b^{s \overline{k}} \vec{e}_{\overline{k}} $. Thus the tangential component is
\[
    \nabla^{\Sigma} \tau = \left( \tau_i + h_{is} b^{s \overline{k}} \tau_{\overline{k}} \right) h^{ij} \vec{s}_j,
\]
Write $ \vec{v} := \tau_i h^{ik} \vec{s}_k $ and $ \vec{w} := \tau_{\overline{k}} b^{\overline{k} m} \vec{s}_m $. Then $ \nabla^{\Sigma} \tau = \vec{v} + \vec{w} $, and
\[
    \Inner{\vec{v}}{\vec{w}}_{h} = \tau_i b^{i \overline{k}} \tau_{\overline{k}} = \frac{1}{4} |\nabla \tau|^2,
\]
therefore
\[
    \begin{aligned}
        \Inner{\nabla^{\Sigma} \tau}{\nabla^{\Sigma} \tau}_{h} ={} & |\vec{v} + \vec{w}|_{h}^2 = |\vec{v} - \vec{w}|_{h}^2 + 4 \Inner{\vec{v}}{\vec{w}}_{h} \\
        ={} & |\tau_i h^{im} \vec{s}_m - \tau_{\overline{k}} b^{\overline{k} m} \vec{s}_m|_{h}^2 + |\nabla \tau|^2 \\
        ={} & |\tau_i h^{im} \vec{s}_m - \tau_{\overline{k}} b^{\overline{k} m} \vec{s}_m|_{h}^2 + 4 \sigma e^{-2K \sigma}.
    \end{aligned}
\]

\paragraph{The final identity}
Substituting the preceding calculations gives
\begin{equation}\label{equation::thesis::formula-de-derivative-of-monotonicity}
    \begin{aligned}
        \frac{\dif}{\dif r} \left[ r^{-\kappa} \int_{\{ \tau \leqslant r \}} e^{nf} \dif \vol_{h} \right] ={} &
        -\frac{\kappa}{2n} r^{-(\kappa + 1)} \int_{\{ \tau = r \} \cap \regs} |\tau_i h^{im} \vec{s}_m - \tau_{\overline{k}} b^{\overline{k} m} \vec{s}_m|_{h}^2 \frac{\rho(\point{x})}{\norm{\nabla^{\mathbb{R}^n} \tau}} \dif \mathcal{H}^{n - 1}(\point{x}) \\
        -{} & \frac{\kappa}{2n} r^{-(\kappa + 1)} \int_{\{ \tau = r \} \cap \regs} \left[ 4 \sigma e^{-2K \sigma} - \tau \frac{2n}{\kappa} \right] \frac{\rho(\point{x})}{\norm{\nabla^{\mathbb{R}^n} \tau}} \dif \mathcal{H}^{n - 1}(\point{x}) \\
        -{} & \frac{\kappa}{2n} r^{-(\kappa + 1)} \int_{\partial \regs \cap \{ \tau \leqslant r \}} h^{ij} \partial_{x^j} \tau(\point{x}, \point{T}(\point{x})) (\vec{\nu}_{\regs})_i \rho(\point{x}) \dif \mathcal{H}^{n - 1}(\point{x}) \\
        +{} & \frac{\kappa}{2n} r^{-(\kappa + 1)} \int_{\{ \tau \leqslant r \}} I_1 \rho(\point{x}) \dif \point{x} \\
        +{} & \frac{\kappa}{2n} r^{-(\kappa + 1)} \int_{\{ \tau \leqslant r \}} \left[ \Inner{\nabla nf}{\nabla \tau} - I_2(\tau) - 2n \right] \rho(\point{x}) \dif \point{x}.
    \end{aligned}
\end{equation}

\section{Step 3 for the monotonicity formula: estimates for the right hand side terms}\label{section::thesis::estimate-c}

In this section, we analyze the sign or an upper bound for each line of \ref{equation::thesis::formula-de-derivative-of-monotonicity}. 

\subsection{Choice of neighborhood}\label{subsection::thesis::choice-of-neighborhood}

This subsection explains how to choose the radius $ r $. The reader may skip it on a first reading and return when the choice of $ r $ is needed. We record the bounds for the cost Hessian, choose a geodesically convex neighborhood for the curves used below, and then estimate the diameter of $ D_r $ and the radial geodesics to its points.

\begin{lemma}[Bounds for the cost Hessian]\label{lemma::thesis::uniform-cost-bounds-on-U}
    There exist $ R>0 $ and $ C>0 $, independent of the base transport pair $ \point{o}=(\point{x}_0,\overline{\point{x}}_0) $, such that the normalized benefit is defined on a neighborhood of the closure of
    \[
        U:=B_R(\point{x}_0)\times B_R(\overline{\point{x}}_0),
    \]
    the cost assumptions hold there, and
    \[
        \|D^2b\|+
        \left\|\left(D_{\point{x}}D_{\overline{\point{x}}}b\right)^{-1}\right\|
        \leqslant C\quad\text{on }U.
    \]
    Moreover, for $ \point{z}\in U $,
    \[
        \|D^2b(\point{z})-D^2b(\point{o})\|
        \leqslant C|\point{z}-\point{o}|.
    \]
\end{lemma}
\begin{proof}
    The cost is smooth on an open neighborhood of $ \cl\regs\times\cl\regt $, and its mixed Hessian is invertible there. Choose a smaller neighborhood with compact closure in that open set. The derivatives of the cost and the inverse mixed Hessian are bounded on its closure. Compactness of $ \cl\regs\times\cl\regt $ gives a common $ R $ for which all the displayed product neighborhoods have closure in this smaller neighborhood. Subtracting the functions used in \cref{definition::thesis::normalization-at-a-contact} leaves the mixed Hessian unchanged and preserves uniform bounds for the other derivatives. The last estimate follows by integrating $ D^3b $ along the segment from $ \point{o} $ to $ \point{z} $.
\end{proof}

\begin{lemma}[Geodesic neighborhoods]\label{lemma::thesis::geodesic-neighborhood}\label{remark::thesis::small-height-and-geometric-neighborhood}
    Inside $ U $, choose a geodesically convex normal neighborhood $ V $ of $ \point{o} $. There is a smaller product neighborhood $ U_1 $ of $ \point{o} $ such that the radial geodesics and the source and target $ c $-segments used in the following local calculations, including their endpoint variations, lie in $ V $ when their initial endpoint data lie in $ U_1 $. This includes a radial geodesic followed by a fixed-coordinate $ c $-segment to an axis, and a fixed-coordinate $ c $-segment followed by radial geodesics to its points.

    For the radial geodesic $ \gamma_{\point{z}} $ from $ \point{o} $ to $ \point{z}\in U_1 $, parametrized on $ [0,1] $, we have
    \[
        \sup_{0\leqslant t\leqslant1}
        \left(|\gamma_{\point{z}}(t)-\point{o}|+
        |\dot\gamma_{\point{z}}(t)|\right)
        \leqslant C|\point{z}-\point{o}|.
    \]
    The neighborhood sizes and the constant can be chosen uniformly in the base transport pair. For the fixed constants $ K,k>0 $ used in the Hessian estimates below, the neighborhoods can be reduced so that
    \[
        K|\sigma|<1,\qquad k|\sigma|<\pi^2/4
    \]
    throughout the region traversed by these curves.
\end{lemma}
\begin{proof}
    The convex normal neighborhood theorem for semi-Riemannian metrics gives $ V $; see \cite{o1983semi}. On this neighborhood the exponential map at $ \point{o} $ has a smooth inverse. Write $ \vec{v}(\point{z})=\exp_{\point{o}}^{-1}(\point{z}) $. Since $ \vec{v}(\point{o})=0 $, smoothness gives $ |\vec{v}(\point{z})|\leqslant C|\point{z}-\point{o}| $. The formulas
    \[
        \gamma_{\point{z}}(t)=\exp_{\point{o}}(t\vec{v}(\point{z})),
        \qquad
        \dot\gamma_{\point{z}}(t)
        =D\exp_{\point{o}}(t\vec{v}(\point{z}))[\vec{v}(\point{z})]
    \]
    prove the stated displacement and velocity bounds after restricting to a smaller neighborhood.

    A source or target $ c $-segment is an ambient geodesic with one coordinate fixed. Its local construction depends smoothly on its endpoints and the fixed coordinate, by invertibility of the mixed Hessian. At the base pair these curves are constant. Choose a product neighborhood with closure in $ V $, and then a smaller product neighborhood $ U_1 $ so that each of the finitely many successive curve constructions just listed stays in the first product neighborhood. Such a choice follows from smooth dependence, uniformly for the curve parameter in $ [0,1] $. Replacing one coordinate by its base value keeps the point in the same product neighborhood. The resulting curves are the short geodesics in $ V $, and the strict containment also permits small endpoint variations.

    The cost bounds from \cref{lemma::thesis::uniform-cost-bounds-on-U} and compactness of the set of base pairs make these choices uniform. Finally, $ \sigma(\point{o})=0 $ and the smooth dependence of $ \sigma $ allow the neighborhoods to be reduced to satisfy the two displayed inequalities for any fixed $ K,k $.
\end{proof}

Although \cref{lemma::thesis::Dr-and-radial-geodesic-estimates} is stated before \cref{lemma::thesis::positiveness-of-sigma} and uses its conclusion, the proof of the latter uses only the geometric neighborhoods established above and does not depend on the former. Thus there is no circular argument.

\begin{lemma}[Diameter of $ D_r $ and radial geodesic estimates]\label{lemma::thesis::Dr-and-radial-geodesic-estimates}\label{lemma::thesis::localization-of-Dr}\label{corollary::thesis::diameter-of-Dr-by-r}
    Choose a product neighborhood $ U_2 $ on which the comparison in \cref{lemma::thesis::positiveness-of-sigma} and the lower bounds in \cref{lemma::thesis::4n-lowerbound-for-u-v} hold. Set $ D_r:=\{\point{x}\in\regs:(\point{x},\point{T}(\point{x}))\in U_2,\ \tau(\point{x},\point{T}(\point{x}))<r\} $. For a neighborhood $ U\Subset U_2 $ of the base pair chosen as in \cref{lemma::thesis::geodesic-neighborhood}, there exists $ r_0>0 $ such that
    \[
        \cl D_r\subset U\qquad\text{for every }0<r<r_0,
    \]
    where $ D_r $ is identified with its transport graph. Furthermore,
    \[
        \sup_{\point{x}\in\cl D_r}
        \left(|\point{x}-\point{x}_0|+
        |\vec{T}(\point{x})-\overline{\point{x}}_0|\right)
        \leqslant Cr^{\frac1{4n}}.
    \]
    In particular, the source diameter and the diameter of the corresponding graph are at most $ Cr^{1/(4n)} $. For each point $ \point{z} $ of this closed graph, the radial geodesic satisfies
    \[
        \sup_{0\leqslant t\leqslant1}
        \left(|\gamma_{\point{z}}(t)-\point{o}|+
        |\dot\gamma_{\point{z}}(t)|\right)
        \leqslant Cr^{\frac1{4n}},
        \qquad
        \int_0^1|\dot\gamma_{\point{z}}(t)|\dif t
        \leqslant Cr^{\frac1{4n}}.
    \]
    These distances and lengths are Euclidean in the source--target product. The constants can be chosen uniformly in the base transport pair.
\end{lemma}
\begin{proof}
    Let $ c_0>0 $ and $ C_0\geqslant0 $ be the constants in the lower bounds for the potentials and the comparison with $ \tau $, respectively. Choose $ d>0 $ so that $ \cl B_d(\point{o})\subset U $. Put $ C_1=2(2/c_0)^{1/(4n)} $, and choose $ r_0>0 $ so that
    \[
        C_0r_0\leqslant1,
        \qquad C_1r_0^{1/(4n)}<d.
    \]
    For $ \point{x}\in D_r $ and $ 0<r<r_0 $, the two estimates give
    \[
        \begin{aligned}
            c_0|\point{x}-\point{x}_0|^{4n}
            &\leqslant u(\point{x})\leqslant\tau+C_0\tau^2\leqslant2r,\\
            c_0|\vec{T}(\point{x})-\overline{\point{x}}_0|^{4n}
            &\leqslant v(\vec{T}(\point{x}))\leqslant\tau+C_0\tau^2\leqslant2r.
        \end{aligned}
    \]
    Hence
    \[
        |\point{x}-\point{x}_0|+
        |\vec{T}(\point{x})-\overline{\point{x}}_0|
        \leqslant C_1r^{1/(4n)}<d.
    \]
    By continuity of $ \vec{T} $, the same estimate holds on the closed graph over $ D_r $. This proves $ \cl D_r\subset U $ and the diameter bound. In particular, the local restriction to $ U_2 $ introduces no additional boundary for these sublevels.

    Applying \cref{lemma::thesis::geodesic-neighborhood} gives
    \[
        \sup_{0\leqslant t\leqslant1}
        \left(|\gamma_{\point{z}}(t)-\point{o}|+
        |\dot\gamma_{\point{z}}(t)|\right)
        \leqslant C|\point{z}-\point{o}|
        \leqslant Cr^{1/(4n)}.
    \]
    Integrating the velocity bound on $ [0,1] $ gives the length estimate. Choosing $ d $ uniformly in the base pair, as permitted by \cref{lemma::thesis::geodesic-neighborhood}, makes $ r_0 $ and all constants uniform.
\end{proof}

In what follows, $ D_r $ denotes this local sublevel. For $ r<r_0 $, its closed graph lies strictly inside the neighborhood on which $ \tau $ is defined.

\subsection{Positivity of \texorpdfstring{$ \sigma $}{σ}}

One major difference between choosing $ z = \sigma $ and choosing $ z = -c $ is that the latter is naturally positive on the optimal transport submanifold: it represents the increase in cost, compared with the optimal plan, after exchanging the destinations of $ (\point{x}, \overline{\point{x}}) $ and $ (\point{x}_0, \overline{\point{x}}_0) $. However, $ \sigma $ is not naturally positive, and nonpositivity would cause difficulties in defining $ D_r $. We therefore first prove this property.

\begin{lemma}[$ \sigma $ is nonnegative on the optimal transport submanifold]\label{lemma::thesis::positiveness-of-sigma}
    Fix a base transport pair $ \point{o}=(\point{x}_0,\overline{\point{x}}_0)\in\cl\Sigma $, and normalize $ u,v,b $ as in \cref{definition::thesis::normalization-at-a-contact}. There exist a neighborhood $ U $ of $ \point{o} $ and a constant $ C>0 $ such that, at every transport pair $ (\point{x},\overline{\point{x}})\in\cl\Sigma\cap U $,
    \[
        \sigma(\point{x},\overline{\point{x}})\geqslant0,
    \]
    and
    \[
        u(\point{x}),\ v(\overline{\point{x}})
        \leqslant \tau(\point{x},\overline{\point{x}})
        +C\tau(\point{x},\overline{\point{x}})^2.
    \]
    In the coordinates $ \point{q}(\point{x})=D_{\overline{\point{x}}}b(\point{x},\overline{\point{x}}_0) $, we also have
    \[
        u_{q^{\overline{k}}}q^{\overline{k}}
        =u_i(\point{x})b^{i\overline{k}}(\point{x},\overline{\point{x}}_0)
        b_{\overline{k}}(\point{x},\overline{\point{x}}_0)
        \geqslant\tau(\point{x},\overline{\point{x}}).
    \]
\end{lemma}

\begin{proof}

\noindent \textbf{Part 1: Positivity of $ \sigma $.} The first-variation calculation in the preceding section gives $ \nabla\sigma(\gamma(1)) = 2\dot{\gamma}(1) $, and hence $ \mathfrak{h}(\nabla\sigma,\nabla\sigma)=4\sigma $. Write $ \gamma(t)=(\point{x}(t),\overline{\point{x}}(t)) $ and set $ s:=\sigma(\gamma(1)) $. The geodesic energy is constant, so
\[
    s = \mathfrak{h}(\dot{\gamma}(t),\dot{\gamma}(t))
      = \dot{x}^i(t)b_{i\overline{k}}(\gamma(t))\dot{\overline{x}}^{\overline{k}}(t),
    \qquad \sigma(\gamma(t))=t^2s.
\]
To determine the sign of $ s $, recall from \cref{definition::thesis::normalization-at-a-contact} that
\[
    b(\point{x}_0,\overline{\point{x}}_0)=0, \quad
    \left.b\right|_{\Sigma}\geqslant0, \quad
    D_{\point{x}}b(\cdot,\overline{\point{x}}_0)
    =D_{\overline{\point{x}}}b(\point{x}_0,\cdot)=0.
\]

Consider
\[
    A(t) := (\dif b) (\vec{\dot{x}}(t)) = \dot{x}^i b_{i},
\]
Then
\[
    \begin{aligned}
        A^\prime(t) ={}& \ddot{x}^i(t) b_i(\gamma(t)) + \dot{x}^i(t) b_{ij}(\gamma(t)) \dot{x}^j(t) + \dot{x}^i(t) b_{i \overline{k}}(\gamma(t)) \dot{\overline{x}}^{\overline{k}}(t) \\
        ={}& \left( b_{jk} - b_i \Gamma^{i}_{jk} \right) \dot{x}^j \dot{x}^k + s.
    \end{aligned}
\]
Here \( b_{jk} - b_i \Gamma^{i}_{jk} \) can be regarded as \( \nabla^{2} b \). For this purpose, we prove the following lemma:

\begin{lemma}\label{lemma::thesis::hessian-of-b}
    In a sufficiently small neighborhood of the base transport pair, there exists $ k > 0 $ such that $ \nabla^2_{\point{x}\point{x}} b \leqslant k \dif_{\point{x}} b \otimes \dif_{\point{x}} b $, and $ \nabla^2_{\overline{\point{x}}\overline{\point{x}}} b \leqslant k \dif_{\overline{\point{x}}} b \otimes \dif_{\overline{\point{x}}} b $.
\end{lemma}
\begin{proof}
    Since the two conclusions are symmetric, we prove only the first. The other follows by exchanging $ \point{x}, \overline{\point{x}} $.

    Take $ (\point{x}, \overline{\point{x}}) $ sufficiently close to the base transport pair. Let $ (\point{x}, \overline{\point{x}}(t)) $, $ 0 \leqslant t \leqslant 1 $, be the target $ c $-segment with $ \overline{\point{x}}(0) = \overline{\point{x}}_0 $ and $ \overline{\point{x}}(1) = \overline{\point{x}} $. Take a source $ c $-geodesic $ (\point{x}(s), \overline{\point{x}}) $ with $ \point{x}(0) = \point{x} $ and $ \dot{\point{x}}(0) = \vec{w} $. By reducing the endpoint neighborhood, these curves remain in the ambient neighborhood where the cost assumptions hold. Define
    \[
        m(t) := \left. \frac{\dif^2}{\dif s^2} b(\point{x}(s), \overline{\point{x}}(t)) \right|_{s = 0}.
    \]
    In the following calculations, the source velocities and accelerations are evaluated at $ s = 0 $, and the derivatives of $ b $ at $ (\point{x},\overline{\point{x}}(t)) $ unless specified otherwise.

    Then
    \[
        \begin{aligned}
            m(t) ={}& b_{i}\ddot{x}^i + b_{ij}\dot{x}^{i} \dot{x}^{j}, \\
            m^{\prime}(t) ={}& b_{i \overline{k}} \ddot{x}^i \dot{\overline{x}}^{\overline{k}} + b_{ij \overline{k}} \dot{x}^{i} \dot{x}^{j} \dot{\overline{x}}^{\overline{k}} \\
            m^{\prime \prime}(t) ={}& b_{i \overline{k}} \ddot{x}^i \ddot{\overline{x}}^{\overline{k}} + b_{i \overline{k} \overline{l}} \ddot{x}^i \dot{\overline{x}}^{\overline{k}} \dot{\overline{x}}^{\overline{l}} + b_{ij \overline{k}} \dot{x}^{i} \dot{x}^{j} \ddot{\overline{x}}^{\overline{k}} + b_{ij \overline{k} \overline{l}} \dot{x}^{i} \dot{x}^{j} \dot{\overline{x}}^{\overline{k}} \dot{\overline{x}}^{\overline{l}}.
        \end{aligned}
    \]
    Along the geodesic $ (\point{x}(s), \overline{\point{x}}) $, the geodesic equation gives
    \[
        \left. \ddot{x}^{k} + b^{k \overline{m}} b_{\overline{m} p q} \dot{x}^{p} \dot{x}^{q} \right|_{(\point{x}(s), \overline{\point{x}})} = 0,
    \]
    Therefore, at $ t = 1 $,
    \[
        m(1) = \left. \left[ b_{ij} - b_k b^{k \overline{m}} b_{\overline{m} ij} \right] \dot{x}^{i} \dot{x}^{j} \right|_{(\point{x}, \overline{\point{x}})} = \nabla^2_{\point{x}\point{x}} b \left[ \vec{\dot{x}}, \vec{\dot{x}} \right]_{(\point{x}, \overline{\point{x}})} = \nabla^2_{\point{x}\point{x}} b \left[ \vec{w}, \vec{w} \right]_{(\point{x}, \overline{\point{x}})}.
    \]
    The normalization gives $ b(\point{x}(s),\overline{\point{x}}_0) = 0 $, hence $ m(0) = 0 $. The source geodesic equation also gives
    \[
        m'(1) = \left[ b_{i\overline{k}}\ddot{x}^i
        + b_{ij\overline{k}}w^i w^j \right]_{(\point{x},\overline{\point{x}})}
        \dot{\overline{x}}^{\overline{k}}(1) = 0.
    \]
    Thus the goal is to prove
    \[
        m(1) \leqslant k \left( \dif_{\point{x}} b(\vec{w}) \right)^2 = k \left( \dot{x}^i b_i \right)^2.
    \]

    Along the geodesic $ (\point{x}, \overline{\point{x}}(t)) $, the geodesic equation gives
    \[
        \left. \ddot{\overline{x}}^{\overline{k}} + b^{\overline{k} m} b_{m \overline{p} \overline{q}} \dot{\overline{x}}^{\overline{p}} \dot{\overline{x}}^{\overline{q}} \right|_{(\point{x}, \overline{\point{x}}(t))} = 0,
    \]
    Hence
    \[
        m^{\prime\prime}(t) = \left[ b_{ij \overline{k} \overline{l}} - b_{ij \overline{m}} b^{\overline{m} n} b_{n \overline{k}\overline{l}} \right]_{(\point{x}, \overline{\point{x}}(t))} \dot{x}^{i} \dot{x}^{j} \dot{\overline{x}}^{\overline{k}} \dot{\overline{x}}^{\overline{l}}
    \]
    The coefficient is twice the corresponding curvature component. The estimate in \cref{equation::thesis::lowerbound-of-Rm}, with the constants adjusted for this factor, gives a common $ k > 0 $ such that
    \begin{equation}\label{equation::thesis::second-derivative-of-m}
        m^{\prime\prime}(t) \geqslant -2k
        \left( w^i b_{i\overline{k}}(\point{x},\overline{\point{x}}(t))
        \dot{\overline{x}}^{\overline{k}}(t) \right)^2.
    \end{equation}
    By the definition of the target $ c $-segment and the normalization,
    \[
        b_i(\point{x},\overline{\point{x}}(t))
        = t b_i(\point{x},\overline{\point{x}}).
    \]
    Differentiating gives
    \[
        b_{i\overline{k}}(\point{x},\overline{\point{x}}(t))
        \dot{\overline{x}}^{\overline{k}}(t)
        = b_i(\point{x},\overline{\point{x}}).
    \]
    Thus the quantity in parentheses in \cref{equation::thesis::second-derivative-of-m} is independent of $ t $. Integrating by parts and using $ m(0)=0 $ and $ m'(1)=0 $, we obtain
    \[
        \begin{aligned}
            m(1)
            ={}& -\int_0^1 t m''(t)\dif t \\
            \leqslant{}& 2k\left(w^i b_i(\point{x},\overline{\point{x}})\right)^2
                \int_0^1 t\dif t \\
            ={}& k\left(\dif_{\point{x}}b(\vec{w})\right)^2.
        \end{aligned}
    \]
    Since $ \vec{w} $ was arbitrary, this proves the source inequality. Interchanging source and target and taking the larger of the two constants proves both assertions.
\end{proof}

The Hessian estimate gives
\[
    A'(t) \leqslant s+kA(t)^2, \qquad A(0)=0.
\]
Let $ y $ solve $ y'=s+ky^2 $ with $ y(0)=0 $, and set $ \phi:=A-y $. On any interval where $ y $ exists,
\[
    \phi'\leqslant k(A+y)\phi, \qquad \phi(0)=0.
\]
Multiplying by $ \exp(-k\int_0^t(A(\xi)+y(\xi))\dif\xi) $ and integrating gives $ \phi\leqslant0 $, so $ A\leqslant y $.

Suppose $ s<0 $. For $ k>0 $, the comparison solution on $ [0,1] $ is
\[
    y_s(t)=-\sqrt{\frac{-s}{k}}\tanh(\sqrt{-ks}\,t).
\]
For $ k=0 $, take $ y_s(t)=st $. In either case, $ A(t)\leqslant y_s(t)<0 $ for $ t>0 $. The target Hessian estimate gives the same inequality for
\[
    \overline{A}(t):=\dot{\overline{x}}^{\overline{k}}(t)b_{\overline{k}}(\gamma(t)).
\]
Since
\[
    \frac{\dif}{\dif t}b(\gamma(t))=A(t)+\overline{A}(t),
\]
and $ b(\gamma(0))=0 $, integration gives
\[
    b(\gamma(1))=\int_0^1\bigl(A(t)+\overline{A}(t)\bigr)\dif t<0.
\]
The endpoint $ \gamma(1) $ is a transport pair, where the normalization and optimality give $ b(\gamma(1))\geqslant0 $. This contradiction proves $ s\geqslant0 $.

\noindent \textbf{Part 2: Upper bounds for $ u $, $ v $.} Use the normalized potentials $ u, v $ from \cref{definition::thesis::normalization-at-a-contact}, with base transport pair $ (\point{x}_0, \overline{\point{x}}_0) $.
Along the ambient geodesic,
\[
    \frac{\dif}{\dif t}b(\point{x}(t),\overline{\point{x}}(t))
    =A(t)+\overline{A}(t).
\]
At the endpoint, $ u+v=b $ and $ u,v\geqslant0 $, so it suffices to estimate $ b(\gamma(1)) $.

Part 1 has already shown that $ s=\sigma(\gamma(1))\geqslant0 $. If $ s>0 $, the neighborhood choice in \cref{lemma::thesis::geodesic-neighborhood} gives $ \sqrt{ks}<\pi/2 $. We may therefore take the comparison function for $ A(t) $ to be
\[
    y_s(t) := \sqrt{\frac{s}{k}}\tan \left( \sqrt{k s} t \right),
\]
At the endpoint, $ u_i=b_i $ and the first-variation formula gives $ \sigma_{\overline{k}}b^{\overline{k}i}=\dot{x}^i(1) $. Therefore
\[
    \sigma_{\overline{k}}b^{\overline{k}i}u_i
    =A(1)\leqslant s+\frac{k}{3}s^2+O(s^3).
\]
Applying the same comparison to $ \overline{A} $ and integrating, using $ b(\point{o})=0 $, gives
\[
    \begin{aligned}
        b(\gamma(1))
        &=\int_0^1\left(A(t)+\overline{A}(t)\right)\dif t\\
        &\leqslant-\frac{2}{k}\log\cos\sqrt{ks}
        =s+\frac{k}{6}s^2+O(s^3).
    \end{aligned}
\]
Hence $ u(\point{x}),v(\overline{\point{x}})\leqslant s+Cs^2 $ at the endpoint.
Moreover,
\[
    \tau = \frac{1 - e^{-K \sigma}}{K} = s - \frac{K}{2} s^2 + O(s^3),
\]
Hence
\[
    \sigma_{\overline{k}} b^{\overline{k} i} u_i \leqslant \tau + C \tau^2, \quad u(\point{x}) \leqslant \tau + C \tau^2.
\]
Similarly, $ v(\overline{\point{x}}), \sigma_i b^{i \overline{k}} v_{\overline{k}} $ satisfy analogous estimates.

If $ s = 0 $, the comparison solution is $ y_s \equiv 0 $. Thus $ A(t) \leqslant 0 $, $ \overline{A}(t) \leqslant 0 $, and hence $ b \leqslant 0 $. By nonnegativity of $ b $, we obtain $ b = u = v = 0 $. Therefore the same upper bounds hold.

\noindent \textbf{Part 3: Lower bounds for the radial derivatives.}

The pure source and target blocks of $ \nabla^2\tau $ are negative semidefinite. Thus $ \tau $ is concave along $ c $-segments with either coordinate fixed, within the neighborhood chosen in \cref{lemma::thesis::geodesic-neighborhood}. Consider the segment $ (\point{x},\overline{\point{y}}(t)) $ from $ (\point{x},\overline{\point{x}}_0) $ to $ (\point{x},\overline{\point{x}}) $. Since $ \tau(\point{x},\overline{\point{x}}_0)=\sigma(\point{x},\overline{\point{x}}_0)=0 $, concavity gives
\[
    \tau(\point{x},\overline{\point{x}})
    \leqslant\left.\frac{\dif}{\dif t}\tau(\point{x},\overline{\point{y}}(t))\right|_{t=0}.
\]
The $ c $-segment satisfies $ D_{\point{x}}b(\point{x},\overline{\point{y}}(t))=tD_{\point{x}}b(\point{x},\overline{\point{x}}) $, so
\[
    \dot{\overline{y}}^{\overline{k}}(0)
    =b^{\overline{k}i}(\point{x},\overline{\point{x}}_0)b_i(\point{x},\overline{\point{x}}).
\]
At the initial point, the first-variation formula on the source $ c $-segment from $ \point{o} $ gives $ \sigma_{\overline{k}}(\point{x},\overline{\point{x}}_0)=b_{\overline{k}}(\point{x},\overline{\point{x}}_0) $. Since $ \sigma=0 $ there, the same identity holds for $ \tau_{\overline{k}} $. Hence
\[
    \begin{aligned}
        \tau(\point{x},\overline{\point{x}})
        &\leqslant\tau_{\overline{k}}(\point{x},\overline{\point{x}}_0)
            \dot{\overline{y}}^{\overline{k}}(0)\\
        &=b_{\overline{k}}(\point{x},\overline{\point{x}}_0)
            b^{\overline{k}i}(\point{x},\overline{\point{x}}_0)
            b_i(\point{x},\overline{\point{x}})\\
        &=u_i(\point{x})b^{\overline{k}i}(\point{x},\overline{\point{x}}_0)
            b_{\overline{k}}(\point{x},\overline{\point{x}}_0).
    \end{aligned}
\]
Here the last equality uses the transport identity at the terminal pair $ (\point{x},\overline{\point{x}}) $.

\end{proof}

\begin{remark}
    The normalization in \cref{definition::thesis::normalization-at-a-contact} makes the value and first derivatives of $ b $ vanish at the starting point $ \point{o} $. The vanishing value is used in the final application of the fundamental theorem of calculus. The vanishing first derivatives are used in the initial value of $ A $. The function itself also becomes the cost/benefit difference obtained by exchanging the destinations, serving several purposes at once.
\end{remark}

\subsection{The second line}
 
From now on, take $ \kappa = n/2 $. Write $ a(r) := -\frac{1}{K}(1 - Kr)^2 \log \left(1 - Kr\right) $ for $ K>0 $, and set $ a(r):=r $ when $ K=0 $. Since $ \tau = r $ on the integration region, the integrand is $ 4(a(r) - r) \rho(\point{x}) / \norm{\nabla^{\mathbb{R}^n} \tau} $. Taylor expansion gives
\[
    a(r) = r - \frac{3}{2} K r^2 + O(r^3).
\]

\subsection{The third line}

In this section we prove
\begin{lemma}[The third line is nonpositive]\label{lemma::thesis::line-3}
    There exists $ r>0 $ such that, for every $ \point{y}\in\partial\regs $ with $ \tau(\point{y},\vec{T}(\point{y}))<r $,
    \[
        \left.
        h^{ij}\partial_{x^j}\left[\tau(\point{x},\vec{T}(\point{x}))\right]
        (\vec{\nu}_{\regs})_i
        \right|_{\point{x}=\point{y}}\geqslant0.
    \]
    Thus the third line is nonpositive because of its preceding minus sign.
\end{lemma}

Expanding the total derivative along the transport graph, the term to be estimated equals
\[
    \tau_i h^{ij} (\vec{\nu}_{\regs})_j + \tau_{\overline{k}} b^{\overline{k} j} (\vec{\nu}_{\regs})_j.
\]

Since $ D\vec{T} $ maps boundary tangent spaces to boundary tangent spaces and $ \vec{T} $ maps the interior to the interior, in the smooth diffeomorphism setting we have $ (D\vec{T})^{\top}\vec{\nu}_{\regt}=\lambda\vec{\nu}_{\regs} $ for some $ \lambda>0 $. Using $ h_{ij}=b_{i\overline{k}}\partial_jT^{\overline{k}} $, this gives
\[
    b_{i\overline{k}}h^{ij}(\vec{\nu}_{\regs})_j
    =\lambda^{-1}(\vec{\nu}_{\regt})_{\overline{k}}.
\]
Thus the first term can also be written as
\[
    \tau_i h^{ij}(\vec{\nu}_{\regs})_j
    =\lambda^{-1}\tau_i b^{i\overline{k}}(\vec{\nu}_{\regt})_{\overline{k}}.
\]

By the first-variation formula, the two terms are therefore positive multiples of the inner products of the source and target components of the geodesic velocity with their respective outward normals.

In \citeauthor{collins2025boundary}'s paper \cite{collins2025boundary}, $ c = - \point{x} \cdot \overline{\point{x}} $. In this case, with the normalization at $ \point{o} = (\point{x}_0, \overline{\point{x}}_0) $, we have $ \tau = b = \Inner{\point{x} - \point{x}_0}{\overline{\point{x}} - \overline{\point{x}}_0} $. Thus $ D_{\overline{\point{x}}}\tau = \point{x} - \point{x}_0 $, while $ b_{i\overline{k}} $ is the identity matrix. At a source boundary point $ \point{y} $, convexity of $ \regs $ gives
\[
    \tau_{\overline{k}} b^{\overline{k} j} (\vec{\nu}_{\regs})_j = (\point{y} - \point{x}_0) \cdot \vec{\nu}_{\regs} \geqslant 0.
\]

If we took $ \tau = b $, then $ \tau_{\overline{k}} = b_{\overline{k}} $. By the normalization, $ D_{\overline{\point{x}}}b(\point{x}_0,\overline{\point{y}})=0 $. The $ c $-convexity of $ \regs $ therefore gives the nonnegativity of the second term, $ \tau_{\overline{k}}b^{\overline{k}j}(\vec{\nu}_{\regs})_j $, by the same argument as in the bilinear case.

However, here $ \tau $ is constructed from $ \sigma $. Its gradient with respect to the ambient metric is a positive multiple of the tangent vector to the radial geodesic. To determine the signs of the source and target components paired with the outward boundary normals, we show that the geodesic $ \gamma(t) $ from $ \point{o} $ to $ (\point{y},\overline{\point{y}}) $ stays in $ \cl\regs\times\cl\regt $. We therefore need to prove

\begin{lemma}[$ c $-convexity implies local geodesic convexity]\label{lemma::thesis::local-geodesic-convexity}
    There exists $ r > 0 $ such that, {for every boundary transport pair $ (\point{y}, \overline{\point{y}}) := (\point{y}, \point{T}(\point{y})) \in \cl\Sigma\cap U_2 $ with $ \point{y}\in\partial\regs $ and $ \tau(\point{y},\overline{\point{y}})<r $,} if the geodesic from $ \point{o} $ to $ (\point{y}, \overline{\point{y}}) $ is $ \gamma(t) = (\point{x}(t), \overline{\point{x}}(t)) $, then
    \[
        \point{x}(t) \in \cl \regs, \quad \overline{\point{x}}(t) \in \cl \regt.
    \]
\end{lemma}

\begin{proof}
    By \cref{lemma::thesis::Dr-and-radial-geodesic-estimates}, the radial geodesics to points of $ D_r $ satisfy $ \sup_{0\leqslant t\leqslant1}(|\gamma(t)-\point{o}|+|\dot\gamma(t)|)\leqslant Cr^\eta $ with $ \eta=1/(4n) $. Since $ \point{o}\in\cl\regs\times\cl\regt $, their source and target components have distance at most $ Cr^\eta $ from the corresponding domains. We use this estimate below to choose $ r $ so that the geodesics remain in the required ambient neighborhood.

    \noindent \textbf{An auxiliary function}

    Suppose there exists $ t_0 $ such that $ \point{y} := \point{x}(t_0) \not \in \cl \regs $. Choose a point $ \point{x}^{\ast} \in \partial \regs $ with the shortest Euclidean distance to $ \point{y} $. Choose $ \overline{\point{z}} \in \cl \regt $, to be specified, and define its corresponding $ \point{q} $ coordinates by $ q_{\overline{k}}(\point{x}) := b_{\overline{k}}(\point{x}, \overline{\point{z}}) $. Then $ Q := \point{q}(\regs) $ is convex. Take the supporting plane at $ \point{q}(\point{x}^{\ast}) $,
    \[
        \ell(\point{q}) := \left( \point{q} - \point{q}(\point{x}^{\ast}) \right) \cdot \vec{\nu}_{Q}(\point{q}(\point{x}^{\ast})),
    \]
    where $ \vec{\nu}_{Q} $ is an outward normal of $ Q $, not necessarily of unit length, chosen so that $ (\vec{\nu}_Q)^{\overline{k}} = (\vec{\nu}_{\regs})_i b^{i \overline{k}}(\point{x}^{\ast}, \overline{\point{z}}) $. Then $ Q \subset \{ \ell \leqslant 0\} $. In the original coordinates this reads
    \[
        \ell(\point{x}) = \left[ b_{\overline{k}}(\point{x}, \overline{\point{z}}) - b_{\overline{k}}(\point{x}^{\ast}, \overline{\point{z}}) \right] b^{i \overline{k}}(\point{x}^{\ast},\overline{\point{z}}) (\vec{\nu}_{\regs})_i.
    \]
    
    \noindent \textbf{Smoothness needed for $ \ell $}

    Write $ \vec{n}=\vec{\nu}_{\regs}(\point{x}^{\ast}) $. The defining function above satisfies $ \ell(\point{x}^{\ast})=0 $ and $ D\ell(\point{x}^{\ast})=\vec{n} $. By \ref{lemma::thesis::uniform-cost-bounds-on-U}, the inverse mixed Hessian is uniformly bounded on the common neighborhood $ U $. We may therefore choose common constants $ R_0>0 $ and $ M_0\geqslant1 $ such that $ \cl B_{R_0}(\point{x}^{\ast})\times\{\overline{\point{z}}\}\subset U $ and $ \|D^2\ell\|\leqslant M_0 $ on this ball. Choose
    \[
        0<R\leqslant\min\{R_0/2,1/(4M_0)\},
        \qquad
        \varepsilon_0:=\min\{R/8,1/(4M_0)\}.
    \]
    Then $ |D\ell-\vec{n}|\leqslant1/4 $ on $ B_R(\point{x}^{\ast}) $, so $ H:=\{\ell=0\}\cap B_R(\point{x}^{\ast}) $ is smooth. Set
    \[
        K_{\ell}:=\{\point{z}\in\cl B_R(\point{x}^{\ast}):\ell(\point{z})\leqslant0\}.
    \]
    For $ 0<s<\varepsilon_0 $, put $ \point{x}_s=\point{x}^{\ast}+s\vec{n} $. For every $ \point{z}\in K_{\ell} $, Taylor's formula gives
    \[
        \ell(\point{z})
        =\vec{n}\cdot(\point{z}-\point{x}^{\ast})
        +\int_0^1(1-t)D^2\ell\bigl(\point{x}^{\ast}+t(\point{z}-\point{x}^{\ast})\bigr)
        [\point{z}-\point{x}^{\ast},\point{z}-\point{x}^{\ast}]\dif t.
    \]
    Since $ \ell(\point{z})\leqslant0 $ and $ \|D^2\ell\|\leqslant M_0 $, it follows that
    \[
        \vec{n}\cdot(\point{z}-\point{x}^{\ast})
        \leqslant\frac{M_0}{2}|\point{z}-\point{x}^{\ast}|^2,
        \qquad
        |\point{x}_s-\point{z}|^2
        \geqslant s^2+(1-M_0s)|\point{z}-\point{x}^{\ast}|^2.
    \]
    Thus $ \point{x}^{\ast} $ is the unique nearest point of $ K_{\ell} $ to $ \point{x}_s $. Also $ \ell(\point{x}_s)\geqslant s-M_0s^2/2>0 $.

    For each $ \point{z}\in H $, define the unit normal $ \vec{n}(\point{z})=D\ell(\point{z})/|D\ell(\point{z})| $. Consider the map $ F(\point{z},s)=\point{z}+s\vec{n}(\point{z}) $. For every tangent vector $ \vec{v}\in T_{\point{x}^{\ast}}H $,
    \[
        D_{\point{z}}F(\point{x}^{\ast},s)[\vec{v}]
        =\vec{v}+sD_{\vec{v}}\vec{n}(\point{x}^{\ast}).
    \]
    Here $ D_{\vec{v}}\vec{n} $ denotes the directional derivative of the unit normal along $ H $. Since $ |D_{\vec{v}}\vec{n}(\point{x}^{\ast})|\leqslant M_0|\vec{v}| $ and $ sM_0<1/4 $, this map is invertible on the tangent space. The derivative of $ F $ with respect to $ s $ is $ \vec{n}(\point{x}^{\ast}) $, perpendicular to that tangent space. The inverse function theorem therefore gives smooth normal coordinates near $ \point{x}_s $. By compactness of $ K_{\ell} $ and uniqueness of the nearest point at $ \point{x}_s $, nearest points for nearby points lie near $ \point{x}^{\ast} $. They lie on $ H $. If $ \point{z}\in H $ is a nearest point to a nearby exterior point $ \point{x} $, then $ \point{x}-\point{z} $ is perpendicular to the tangent space of $ H $ at $ \point{z} $, so $ \point{x}=\point{z}+t\vec{n}(\point{z})=F(\point{z},t) $, where $ t=|\point{x}-\point{z}|>0 $. The local inverse of $ F $ therefore makes this nearest point unique. Consequently the signed distance $ d $ to $ H $, positive on the side $ \ell>0 $, is smooth near $ \point{x}_s $, with
    \[
        d(\point{x}_s)=s,
        \qquad Dd(\point{x}_s)=\vec{n}.
    \]

    At the selected exterior point $ \point{y} $, write $ a=\dist_{\mathbb{R}^n}(\point{y},\cl\regs) $. Then $ \point{y}=\point{x}^{\ast}+a\vec{n} $. Once the preliminary geodesic estimate gives $ a<\varepsilon_0 $, the preceding construction applies. Since $ \ell\leqslant0 $ on $ \cl\regs $, nearest points of $ \cl\regs $ for points near $ \point{y} $ lie in $ K_{\ell} $. Indeed, the displayed distance inequality makes $ \point{x}^{\ast} $ the unique nearest point of $ \cl\regs\cap\cl B_R(\point{x}^{\ast}) $ to $ \point{y} $, while points outside this ball are farther away. Compactness then gives the same localization for nearby points. Hence
    \[
        d(\point{y})=\dist_{\mathbb{R}^n}(\point{y},\cl\regs),
        \qquad
        d(\point{x})\leqslant\dist_{\mathbb{R}^n}(\point{x},\cl\regs)
        \quad\text{for $ \point{x} $ near $ \point{y} $.}
    \]

    We next compute the second derivatives of $ d $. Since $ |Dd|=1 $, we have $ D^2d[\vec{n},\cdot]=0 $ at $ \point{x}^{\ast} $. Given a tangent vector $ \vec{w}\in T_{\point{x}^{\ast}}H $, choose a curve $ \point{w}(s)\subset H $ with $ \point{w}(0)=\point{x}^{\ast} $ and $ \point{w}'(0)=\vec{w} $. Differentiating $ d(\point{w}(s))=\ell(\point{w}(s))=0 $ twice, and using $ Dd(\point{x}^{\ast})=D\ell(\point{x}^{\ast})=\vec{n} $, gives
    \[
        \begin{aligned}
            D^2d(\point{x}^{\ast})[\vec{w},\vec{w}]
            & =D^2\ell(\point{x}^{\ast})[\vec{w},\vec{w}] \\
            & =b_{\overline{k}ij}(\point{x}^{\ast},\overline{\point{z}})
            b^{m\overline{k}}(\point{x}^{\ast},\overline{\point{z}})n_m w^i w^j \\
            & =n_k\Gamma^k_{ij}(\point{x}^{\ast},\overline{\point{z}})w^i w^j.
        \end{aligned}
    \]
    Define the symmetric linear map $ \matr{G}:T_{\point{x}^{\ast}}H\to T_{\point{x}^{\ast}}H $ by $ \matr{G}\vec{w}=D_{\vec{w}}\vec{n}(\point{x}^{\ast}) $. Thus, for tangent vectors $ \vec{v},\vec{w} $,
    \[
        \Inner{\vec{v}}{\matr{G}\vec{w}}
        =D^2d(\point{x}^{\ast})[\vec{v},\vec{w}]
        =n_k\Gamma^k_{ij}(\point{x}^{\ast},\overline{\point{z}})v^i w^j.
    \]
    These identities identify $ \matr{G} $ with the restriction of the bilinear form $ n_k\Gamma^k_{ij} $ to tangent vectors.

    At $ \point{y}=\point{x}^{\ast}+a\vec{n} $, we likewise have $ D^2d(\point{y})[\vec{n},\cdot]=0 $. For tangential directions, consider $ \tilde{\point{w}}(s)=\point{w}(s)+a\vec{n}(\point{w}(s)) $. Differentiating $ Dd(\tilde{\point{w}}(s))=\vec{n}(\point{w}(s)) $ at $ s=0 $ gives
    \[
        D^2d(\point{y})(\vec{w}+a\matr{G}\vec{w})=\matr{G}\vec{w}.
    \]
    Therefore, for every tangent vector $ \vec{w} $,
    \[
        D^2d(\point{y})[\vec{w},\vec{w}]
        =\Inner{\vec{w}}{\matr{G}(\matr{I}+a\matr{G})^{-1}\vec{w}},
    \]
    where $ \matr{I} $ is the identity on $ T_{\point{x}^{\ast}}H $. For an arbitrary vector $ \vec{v}=\vec{w}+f\vec{n} $ with $ \vec{w} $ tangent, the same expression equals $ D^2d(\point{y})[\vec{v},\vec{v}] $, since the normal component contributes zero.

    \noindent \textbf{Constructing the inequality}

    Define $ A(t):=\dist_{\mathbb{R}^n}(\point{x}(t),\cl\regs) $ and $ \overline{A}(t):=\dist_{\mathbb{R}^n}(\overline{\point{x}}(t),\cl\regt) $. Write $ M_A:=\max_{[0,1]}A $, $ M_{\overline{A}}:=\max_{[0,1]}\overline{A} $, and $ M:=M_A+M_{\overline{A}} $. We must show that $ M=0 $. Suppose instead that $ M>0 $.

    Take any $ t_0\in(0,1) $ with $ A(t_0)>0 $. Set $ \point{y}:=\point{x}(t_0) $, $ \overline{\point{y}}:=\overline{\point{x}}(t_0) $, and $ a:=A(t_0) $. Choose a nearest point $ \point{x}^{\ast}\in\cl\regs $ to $ \point{y} $ and a nearest point $ \overline{\point{z}}\in\cl\regt $ to $ \overline{\point{y}} $. Construct $ d $ as above and set $ s(t):=d(\point{x}(t)) $ near $ t_0 $. Then $ s(t_0)=A(t_0) $ and $ s(t)\leqslant A(t) $ near $ t_0 $. All velocities in the following calculations are evaluated at $ t_0 $. The chain rule and the geodesic equation give
    \[
        \begin{aligned}
            s'(t_0)&=\dot{\point{x}}\cdot\vec{n},\\
            s''(t_0)&=\ddot{\point{x}}\cdot\vec{n}
                +D^2d(\point{y})[\dot{\point{x}},\dot{\point{x}}]\\
                &=-n_k\Gamma^k_{ij}(\point{y},\overline{\point{y}})\dot{x}^i\dot{x}^j
                +\Inner{\vec{w}}{\matr{G}(\matr{I}+a\matr{G})^{-1}\vec{w}},
        \end{aligned}
    \]
    where $ \dot{\point{x}}=\vec{w}+f\vec{n} $, $ \vec{w} $ is tangent to $ H $ at $ \point{x}^{\ast} $, and $ f=s'(t_0) $. Using the tangential identity for $ \matr{G} $, we obtain
    \[
        \begin{aligned}
            s''(t_0)={}&n_k\left[\Gamma^k_{ij}(\point{x}^{\ast},\overline{\point{z}})
                    -\Gamma^k_{ij}(\point{y},\overline{\point{y}})\right]\dot{x}^i\dot{x}^j\\
                &+\Inner{\vec{w}}{\left[\matr{G}(\matr{I}+a\matr{G})^{-1}-\matr{G}\right]\vec{w}}\\
                &-2f n_k\Gamma^k_{ij}(\point{x}^{\ast},\overline{\point{z}})w^i n^j
                 -f^2 n_k\Gamma^k_{ij}(\point{x}^{\ast},\overline{\point{z}})n^i n^j.
        \end{aligned}
    \]
    The bounds for $ b $ and its inverse mixed Hessian give a common bound $ C_1 $ for $ \|\matr{G}\| $ and for the bilinear form with coefficients $ n_k\Gamma^k_{ij} $. By \cref{lemma::thesis::Dr-and-radial-geodesic-estimates}, we may choose $ r $ so that $ a\leqslant M_A\leqslant Cr^\eta\leqslant1/(2C_1) $. Since
    \[
        \matr{G}(\matr{I}+a\matr{G})^{-1}-\matr{G}
        =-a\matr{G}^2(\matr{I}+a\matr{G})^{-1},
    \]
    we have
    \[
        \left\|\matr{G}(\matr{I}+a\matr{G})^{-1}-\matr{G}\right\|
        \leqslant2M_A C_1^2.
    \]
    Consequently the last three terms in the expression for $ s''(t_0) $ are bounded below by
    \[
        -2M_A C_1^2|\vec{w}|^2-2C_1|\vec{w}||f|-C_1|f|^2
        \geqslant-2M C_1^2|\dot{\point{x}}|^2-C_2|\dot{\point{x}}||s'(t_0)|.
    \]
    The derivative bounds for the Christoffel symbols therefore give
    \[
        \begin{aligned}
            s''(t_0)\geqslant{}&-C_3\left(|\point{y}-\point{x}^{\ast}|+
                    |\overline{\point{y}}-\overline{\point{z}}|\right)|\dot{\point{x}}|^2
                    -2M C_1^2|\dot{\point{x}}|^2
                    -C_2|\dot{\point{x}}||s'(t_0)|\\
                \geqslant{}&-C_4 M r^{2\eta}-C_4r^\eta|s'(t_0)|.
        \end{aligned}
    \]
    Here we used $ |\point{y}-\point{x}^{\ast}|\leqslant M_A $, $ |\overline{\point{y}}-\overline{\point{z}}|\leqslant M_{\overline{A}} $, and the velocity bound from \cref{lemma::thesis::Dr-and-radial-geodesic-estimates}. The same argument applies with source and target interchanged, with a common constant $ C_4 $.

    \noindent \textbf{Obtaining a contradiction}

    Set $ a_1:=C_4r^{2\eta} $ and $ a_2:=C_4r^\eta $. Choose $ r $ so that $ 0<a_1<2 $ and $ a_2<1 $.
    \begin{claim}
        $ A(t)\leqslant a_1Mt(1-t) $ for $ 0\leqslant t\leqslant1 $.
    \end{claim}
    \begin{proof}
        Otherwise, since $ A(0)=A(1)=0 $, the function $ A(t)-a_1Mt(1-t) $ attains a positive maximum at some $ t_1\in(0,1) $. In particular $ A(t_1)>0 $. Apply the preceding construction at $ t_1 $ to obtain a smooth function $ s $ with $ s(t_1)=A(t_1) $ and $ s\leqslant A $ near $ t_1 $. Then $ s(t)-a_1Mt(1-t) $ has a local maximum at $ t_1 $, so
        \[
            s'(t_1)=a_1M(1-2t_1),
            \qquad s''(t_1)\leqslant-2a_1M.
        \]
        The differential inequality gives
        \[
            -2a_1M\geqslant s''(t_1)
            \geqslant-a_1M-a_2|s'(t_1)|
            =-\left(1+a_2|1-2t_1|\right)a_1M
            >-2a_1M,
        \]
        a contradiction.
    \end{proof}
    Thus $ M_A\leqslant a_1M/4 $. By symmetry, $ M_{\overline{A}}\leqslant a_1M/4 $. Hence
    \[
        M=M_A+M_{\overline{A}}\leqslant\frac{a_1}{2}M<M,
    \]
    again a contradiction. Therefore $ M=0 $, proving the containment assertion.
\end{proof}

\begin{remark}
    This proof is similar to \cite[Lemma 3.5]{caponio2013infinitesimal}: it uses the distance function and obtains a similar differential inequality.
\end{remark}

\begin{proof}[Proof of \ref{lemma::thesis::line-3}]
    With the identification of the boundary normals above, the two terms to be estimated are positive multiples of the inner products of the source and target geodesic directions with their respective outward normals. These inner products are nonnegative by \cref{lemma::thesis::local-geodesic-convexity}. Hence the third line is nonpositive.
\end{proof}

\subsection{The fifth line}

In \cite{collins2025boundary}, only the term $ \Inner{\nabla nf}{\nabla \tau} $ needs to be considered. In their case this equals $ \rho_i x^i / \rho + \overline{\rho}_{\overline{k}} \overline{x}^{\overline{k}} / \overline{\rho} $, so it can be written as derivatives of $ \rho $ and $ \overline{\rho} $ along the straight line $ t \point{x} $ from $ \point{o} $. Writing the integral over $ D_r $ by first integrating along rays from $ o $ within $ D_r $ and then integrating over directions gives their ``subhomogeneity inequality''. The fact that each ray meets $ D_r $ in an interval allows the integral to be written in this way. The choice $ z = b $ may fail because proving $ D_r $ is ``convex'' with respect to the required curves becomes difficult. But with $ \tau $ we can do it. Actually we use a divergence theorem directly to prove a similar result.

More precisely, the conclusion to be proved in this section is
\begin{lemma}
    There exists $ R > 0 $ such that, when $ r \leqslant R $,
    \[
        \int_{D_r} \left[ \Inner{\nabla nf}{\nabla \tau} - I_2(\tau) - 2n \right] \rho(\point{x}) \dif \point{x} \leqslant \delta(r)\int_{D_r}\rho(\point{x})\dif\point{x}.
    \]
\end{lemma}

\begin{proof}
    By definition, these terms expand as
    \[
        \begin{aligned}
            \Inner{\nabla nf}{\nabla \tau} ={}& \tau_i b^{i \overline{k}} \overline{\rho}_{\overline{k}} / \overline{\rho} + \tau_{\overline{k}} b^{\overline{k} i} \rho_i / \rho - \tau_i b^{i \overline{k}} D_{\overline{k}} \log\left|\det D^2_{\point{x}\overline{\point{x}}} b\right| - \tau_{\overline{k}} b^{\overline{k} i} D_{i} \log\left|\det D^2_{\point{x}\overline{\point{x}}} b\right| \\
            I_2(\tau) ={}& -2 b^{i \overline{k}} \tau_{i \overline{k}}.
        \end{aligned}
    \]

    We can change the integration variable from $ \point{x} $ to $ \overline{\point{x}} $, with the integration weight $ \rho(\point{x}) \dif \point{x} $ changing to $ \overline{\rho}(\overline{\point{x}}) \dif \overline{\point{x}} $. Thus the integrand in the lemma can be split into two parts, namely,
    \[
        \tau_{\overline{k}} b^{\overline{k} i} \rho_i / \rho - \tau_{\overline{k}} b^{\overline{k} i} D_{i} \log\left|\det D^2_{\point{x}\overline{\point{x}}} b\right| + b^{i \overline{k}} \tau_{i \overline{k}} - n
    \] 
    integrated against $ \rho(\point{x}) \dif \point{x} $, and
    \[
        \tau_i b^{i \overline{k}} \overline{\rho}_{\overline{k}} / \overline{\rho} - \tau_i b^{i \overline{k}} D_{\overline{k}} \log\left|\det D^2_{\point{x}\overline{\point{x}}} b\right| + b^{i \overline{k}} \tau_{i \overline{k}} - n
    \]
    integrated against $ \overline{\rho}(\overline{\point{x}}) \dif \overline{\point{x}} $. By symmetry, it suffices to analyze the first part below.

    \noindent \textbf{Rewriting the expression}
    Denote the vector field on $ \regs $ by
    \[
        \vec{X}^{i}(\point{x}, \point{T}(\point{x})) := \tau_{\overline{k}}(\point{x}, \point{T}(\point{x})) b^{\overline{k} i}(\point{x}, \point{T}(\point{x})),
    \]
    namely, the restriction of $ \tau_{\overline{k}} b^{\overline{k} i} $ to $ \Sigma $. Taking its Euclidean divergence in the source coordinates gives
    \[
        \begin{aligned}
            \divr_{\mathbb{R}^n} \vec{X} ={}& \partial_i \left[ \tau_{\overline{k}} b^{\overline{k} i} \right] \\
            ={}& \tau_{i \overline{k}} b^{\overline{k} i} + \tau_{\overline{k} \overline{l}} b^{\overline{k} i} b^{\overline{l} m} h_{mi} - \tau_{\overline{k}} b^{\overline{k} m} b_{m \overline{n} i} b^{\overline{n} i} - \tau_{\overline{k}} b^{\overline{k} m} b_{m \overline{n} \overline{l}} b^{\overline{n} i} b^{\overline{l} a} h_{ai} \\
            ={}& \tau_{i \overline{k}} b^{\overline{k} i} - \tau_{\overline{k}} b^{\overline{k} m} D_m \log\left|\det D^2_{\point{x}\overline{\point{x}}} b\right| + \left( \nabla^2 \tau \right)_{\overline{k} \overline{l}} b^{\overline{k} i} h_{mi} b^{\overline{l} m},
        \end{aligned}
    \]
    The contraction of the target Hessian can be written as the trace of four matrix factors:
    \[
        \left(\nabla^2\tau\right)_{\overline{k}\overline{l}}
        b^{\overline{k}i}h_{mi}b^{\overline{l}m}
        =\tr\left[
            \left((\nabla^2\tau)_{\overline{k}\overline{l}}\right)
            \left(D^2_{\point{x}\overline{\point{x}}}b\right)^{-1}
            (h_{im})
            \left(D^2_{\point{x}\overline{\point{x}}}b\right)^{-\top}
        \right]\leqslant0.
    \]
    Indeed, the first factor is negative semidefinite by \cref{subsection::thesis::analyse-of-I_1}, while the product of the last three factors is positive definite. Therefore
    \[
        \divr_{\mathbb{R}^n} \vec{X} \leqslant \tau_{i \overline{k}} b^{\overline{k} i} - \tau_{\overline{k}} b^{\overline{k} m} D_m \log\left|\det D^2_{\point{x}\overline{\point{x}}} b\right|.
    \]

    Applying the divergence theorem to the term $ \vec{X}^i \rho_i / \rho $ to be estimated gives
    \[
        \begin{aligned}
            \int_{D_r} \vec{X}^i \rho_i \rho^{-1} \rho \dif \point{x} ={}& \int_{D_r} \vec{X}^i \rho_i \dif \point{x} \\
            ={}& \int_{\partial D_r} \vec{X} \cdot \vec{\nu} \rho \dif \mathcal{H}^{n - 1}(\point{x}) - \int_{D_r} \rho \divr_{\mathbb{R}^n} \vec{X} \dif \point{x} \\
        \end{aligned}
    \]
    The first boundary integral splits into two parts. On the part {$ \{\point{x}\in\partial\regs:(\point{x},\point{T}(\point{x}))\in U_2,\ \tau(\point{x},\point{T}(\point{x}))<r\} $} of the domain boundary,
    \[
        \vec{X} \cdot \vec{\nu} = \tau_{\overline{k}} b^{\overline{k} i} \nu_i,
    \]
    This is nonnegative by Lemma \ref{lemma::thesis::line-3} in the preceding section. On the boundary part $ \{ \tau = r \} \cap \regs $,
    \[
        \norm{D \tau}_{\mathbb{R}^n} \vec{X} \cdot \vec{\nu} = \tau_{\overline{k}} b^{\overline{k} i} \tau_i + \tau_{\overline{k}} b^{\overline{k} i} h_{ij} b^{j \overline{l}} \tau_{\overline{l}},
    \]
    The first part expands as $ e^{-2K\sigma} \sigma_{\overline{k}} b^{\overline{k} i} \sigma_{i} = \sigma e^{-2K\sigma} $, which is nonnegative by Lemma \ref{lemma::thesis::positiveness-of-sigma}. The second part is a positive definite matrix multiplied on both sides by the vector $ \tau_{\overline{k}} $, and is also nonnegative. Hence the entire integrand $ \vec{X} \cdot \vec{\nu} $ is nonnegative on the boundary. Taking $ \rho_{\ast} := \sup_{D_r} \rho $, we obtain
    \[
        \begin{aligned}
            \int_{D_r} \vec{X}^i \rho_i \dif \point{x} \leqslant{}& \rho_{\ast} \int_{\partial D_r} \vec{X} \cdot \vec{\nu} \dif \mathcal{H}^{n - 1}(\point{x}) - \int_{D_r} \rho \divr_{\mathbb{R}^n} \vec{X} \dif \point{x} \\
            ={}& \int_{D_r} (\rho_{\ast} - \rho) \divr_{\mathbb{R}^n} \vec{X} \dif \point{x}.
        \end{aligned}
    \]

    \noindent \textbf{Estimating the remaining terms}
    We next need an upper bound for $ \int_{D_r} \divr_{\mathbb{R}^n} \vec{X} \dif \point{x} $ and an estimate for the remaining part $ - \tau_{\overline{k}} b^{\overline{k} i} D_{i} \log\left|\det D^2_{\point{x}\overline{\point{x}}} b\right| + b^{i \overline{k}} \tau_{i \overline{k}} - n $. Define
    \[
        E_{\point{x}}^{\tau} := b^{i \overline{k}} \tau_{i \overline{k}} - \tau_{\overline{k}} b^{\overline{k} m} D_m \log\left|\det D^2_{\point{x}\overline{\point{x}}} b\right| - n
    \]
    Then the original integrand is precisely $ \tau_{\overline{k}} b^{\overline{k} i} \rho_i / \rho + E_{\point{x}}^{\tau} $, and the preceding analysis shows $ \divr_{\mathbb{R}^n} \vec{X} \leqslant E^{\tau}_{\point{x}} + n $.

    Define analogous notation for $ \sigma $,
    \[
        E_{\point{x}}^{\sigma} := b^{i \overline{k}} \sigma_{i \overline{k}} - \sigma_{\overline{k}} b^{\overline{k} m} D_m \log\left|\det D^2_{\point{x}\overline{\point{x}}} b\right| - n,
    \]
    Then
    \[
        E_{\point{x}}^{\tau} = e^{-K \sigma} E_{\point{x}}^{\sigma}
                + n \left( e^{-K \sigma} - 1 \right)
                - K \sigma e^{-K \sigma}.
    \]
    Thus, once we have a bound for $ E^{\sigma}_{\point{x}} $, the bound for $ E^{\tau}_{\point{x}} $ merely combines it with some elementary functions of $ r $.

    For the geodesic $ \gamma(t) $ defining $ \sigma $, if one coordinate of the endpoint agrees with that of the starting point $ o $, the geodesic is a $ c $-segment. For example, when the endpoint is $ (\point{x}_0, \overline{\point{x}}) $, $ D_{\point{x}} b (\point{x}_0, \overline{\point{x}}(t)) $ is linear in $ t $, so $ \dot{\overline{x}}^{\overline{k}}(t) = b^{i \overline{k}} b_i(\point{x}_0, \overline{\point{x}}) $. Therefore $ \sigma_i(\point{x}_0, \overline{\point{x}}) = b_{i \overline{k}} \dot{\overline{x}}^{\overline{k}} = b_i(\point{x}_0, \overline{\point{x}}) $, and hence
    \[
        \sigma_{i \overline{k}}(\point{x}_0, \overline{\point{x}}) = b_{i \overline{k}}(\point{x}_0, \overline{\point{x}}).
    \]
    Similarly,
    \[
        \sigma_{i \overline{k}}(\point{x}, \overline{\point{x}}_0) = b_{i \overline{k}}(\point{x}, \overline{\point{x}}_0).
    \]
    Therefore
    \[
        \left. \sigma_{i \overline{k}} b^{i \overline{k}} \right|_{\point{x}_0, \overline{\point{x}}_0} = n, \quad |\sigma_{i \overline{k}} b^{i \overline{k}} - n| \leqslant C(n, \norm{b}_{C^4}, \norm{(D^2_{\point{x}\overline{\point{x}}}b)^{-1}}) \norm{\point{x} - \point{x}_0} \norm{\overline{\point{x}} - \overline{\point{x}}_0}.
    \]
    
    For the remaining term in $ E^{\sigma}_{\point{x}} $, we have
    \[
        |\sigma_{\overline{k}} b^{\overline{k} m} D_m \log\left|\det D^2_{\point{x}\overline{\point{x}}} b\right|| \leqslant C(n, \norm{b}_{C^4}, \norm{(D^2_{\point{x}\overline{\point{x}}}b)^{-1}}) \norm{\point{x} - \point{x}_0}.
    \]
    Combining these gives
    \[
        |E^{\sigma}_{\point{x}}| \leqslant C \left( \norm{\point{x} - \point{x}_0} + \norm{\point{x} - \point{x}_0} \norm{\overline{\point{x}} - \overline{\point{x}}_0} \right)
    \]
    Returning to $ E^{\tau}_{\point{x}} $, we obtain
    \[
        |E^{\tau}_{\point{x}}| \leqslant C \left( \norm{\point{x} - \point{x}_0} + \norm{\point{x} - \point{x}_0} \norm{\overline{\point{x}} - \overline{\point{x}}_0} + |\sigma| \right).
    \]
    
    With the upper bound for $ E^{\tau}_{\point{x}} $,
    \[
        \begin{aligned}
            {}& \int_{D_r} \tau_{\overline{k}} b^{\overline{k} i} \rho_i + \left[ -\tau_{\overline{k}} b^{\overline{k} i} D_{i} \log\left|\det D^2_{\point{x}\overline{\point{x}}} b\right| + b^{i \overline{k}} \tau_{i \overline{k}} - n \right] \rho(\point{x}) \dif \point{x} \\
            ={}& \int_{D_r} \vec{X}^i \rho_i + E^{\tau}_{\point{x}} \rho(\point{x}) \dif \point{x} \\
            \leqslant{}& \int_{D_r} \left( \rho_{\ast} - \rho \right) \divr_{\mathbb{R}^n} \vec{X} + E^{\tau}_{\point{x}} \rho(\point{x}) \dif \point{x} \\
            \leqslant{}& \int_{D_r} \left( \rho_{\ast} - \rho \right) (E^{\tau}_{\point{x}} + n) + E^{\tau}_{\point{x}} \rho(\point{x}) \dif \point{x} \\
            ={}& n \int_{D_r} \rho_{\ast} - \rho \dif \point{x} + \rho_{\ast} \int_{D_r} E^{\tau}_{\point{x}} \dif \point{x} \\
            \leqslant{}& \delta(r) \int_{D_r} \rho \dif \point{x}.
        \end{aligned}
    \]
    In the last line, denote this more complicated coefficient, composed of powers of $ r $, by $ \delta(r) $. By \cref{lemma::thesis::Dr-and-radial-geodesic-estimates}, the source and target displacements on $ D_r $ are bounded by $ Cr^\eta $, where we denote the shrinking exponent by $ \eta>0 $ and retain $ \alpha $ for the H\"older exponent of the densities. We can therefore take
    \[
        \delta(r) = C\left(r^{\alpha\eta}+r^\eta+r^{2\eta}+r\right),
    \]
    where $ C>0 $ depends only on the fixed problem data.

    On the target side, define
    \[
        Y^{\overline{k}}:=\tau_i b^{i\overline{k}},
        \qquad
        E^\tau_{\overline{\point{x}}}
        :=b^{i\overline{k}}\tau_{i\overline{k}}
        -Y^{\overline{k}}D_{\overline{k}}\log\left|\det D^2_{\point{x}\overline{\point{x}}}b\right|-n,
    \]
    evaluated at $ (\vec{T}^{-1}(\overline{\point{x}}),\overline{\point{x}}) $. Interchanging source and target in the preceding calculation gives
    \[
        \int_{\vec{T}(D_r)}
        \left(Y^{\overline{k}}\overline{\rho}_{\overline{k}}
        +E^\tau_{\overline{\point{x}}}\overline{\rho}\right)\dif\overline{\point{x}}
        \leqslant\delta(r)\int_{\vec{T}(D_r)}\overline{\rho}\dif\overline{\point{x}},
    \]
    with the target density's H\"older constant in the estimate. Adding the source and target estimates and increasing the constant in $ \delta $ gives the stated inequality.
\end{proof}

\subsection{The final form}

Write $ F(r):=\int_{D_r}\rho(\point{x})\dif\point{x} $. Combining the preceding estimates and increasing the constant in $ \delta $ if necessary gives
\[
    \frac{\dif}{\dif r}\left[r^{-\frac n2}F(r)\right]
    \leqslant-(a(r)-r)r^{-\frac n2-1}F'(r)
    +r^{-\frac n2-1}\delta(r)F(r).
\]
Define $ \chi(r):=r^{-n/2}F(r) $. Then
\[
    \frac{a(r)}r\chi'(r)
    \leqslant\left[\frac{\delta(r)}r
    +\frac n2\frac{r-a(r)}{r^2}\right]\chi(r).
\]
Set
\[
    \beta(r):=\frac r{a(r)}\left(\delta(r)
    +\frac n2\frac{r-a(r)}r\right).
\]
Since $ a(r)=r-\frac32Kr^2+O(r^3)>0 $ for sufficiently small $ r $, we obtain
\[
    \chi'(r)\leqslant\frac{\beta(r)}r\chi(r).
\]
Consequently,
\[
    \exp\left(-\int_0^r\frac{\beta(t)}t\dif t\right)\chi(r)
\]
is nonincreasing.

For sufficiently small $ r<1 $, we have $ r/a(r)<2 $ and $ (r-a(r))/r=O(r) $. Therefore
\[
    0\leqslant\beta(r)
    \leqslant C\left(r^{\alpha\eta}+r^\eta+r^{2\eta}+r\right)
    \leqslant Cr^{\alpha\eta},
\]
and
\[
    \int_0^r\frac{\beta(t)}t\dif t
    \leqslant C\int_0^r\frac{t^{\alpha\eta}+t^\eta+t^{2\eta}+t}{t}\dif t
    \leqslant Cr^{\alpha\eta}.
\]
Put $ \epsilon:=\alpha\eta $. Since $ \beta(r)\leqslant Cr^\epsilon $, increasing $ C $ if necessary gives
\[
    \frac{\dif}{\dif r}\left(e^{-Cr^\epsilon}\chi(r)\right)\leqslant0.
\]
Here $ 0<\alpha\leqslant1 $ is the H\"older exponent of the densities; if both densities are constant, we may take $ \alpha=1 $. The diameter estimate gives $ \eta=1/(4n) $, so $ 0<\epsilon<1 $.

\subsection{The monotonicity formula for H\"older densities}\label{subsection::thesis::weak-monotonicity}

We now establish the monotonicity formula under the regularity assumptions of \cref{theorem::thesis::maintheorem}.

\begin{lemma}[Interior regularity and integrability of derivatives]\label{lemma::thesis::integrable-derivatives-for-monotonicity}
    Under the assumptions above and the preliminary uniqueness of the transport pairs and $ C^{1,\theta} $ estimates, the transport map takes $ \regs $ onto $ \regt $, and
    \[
        u \in C^{2,\alpha}_{\mathrm{loc}}(\regs) \cap W^{2,1}(\regs), \qquad v \in C^{2,\alpha}_{\mathrm{loc}}(\regt) \cap W^{2,1}(\regt).
    \]
    Moreover, $ \vec{T} $ and $ \vec{T}^{-1} $ are $ C^{1,\alpha} $ local diffeomorphisms in the interiors, and
    \[
        \vec{T} \in W^{1,1}(\regs), \qquad \vec{T}^{-1} \in W^{1,1}(\regt).
    \]
\end{lemma}
\begin{proof}
    By \cref{lemma::thesis::unique-closure-contacts}, the transport map is a homeomorphism of the closures. Invariance of domain makes the image of $ \regs $ open in $ \mathbb{R}^n $. Since this image lies in $ \cl \regt $, it lies in $ \regt $. Applying the same argument to the inverse gives $ \vec{T}(\regs)=\regt $.

    The potentials are $ C^1 $ by \cref{theorem::thesis::C1a-estimate-under-A3s}. A supporting mountain cannot touch at a second source point, by uniqueness of the source point paired with a fixed target. More quantitatively, \cref{lemma::thesis::uniform-positive-contact-gap} supplies a positive lower bound for the difference between the potential and its supporting mountain outside every fixed neighborhood of the point of tangency. On each compact subset of $ \regs $, its transport image is compactly contained in $ \regt $. These facts give strict $ c $-convexity and its local modulus. The interior regularity theorem of Liu, Trudinger, and Wang, \cite[Theorem 1.2 and Remark 4.1]{liu2010interior}, applies under A3w, hence under A3s, and gives interior $ C^2 $ regularity for the positive Dini-continuous densities. 
    
    Differentiating the identities $ u_i=b_i(\point{x},\vec{T}(\point{x})) $ and $ v_{\overline{k}}=b_{\overline{k}}(\vec{T}^{-1}(\overline{\point{x}}),\overline{\point{x}}) $ gives
    \[
        \partial_i T^{\overline{k}} = b^{\overline{k}m}h_{mi}, \qquad h_{\overline{k}\overline{l}} = b_{i\overline{k}}h^{im}b_{m\overline{l}}.
    \]
    Thus the map and its inverse are interior $ C^1 $ local diffeomorphisms. On compact interior subsets they are Lipschitz, so the equation
    \[
        \det h = |\det b_{i\overline{k}}|\frac{\rho(\point{x})}{\overline{\rho}(\vec{T}(\point{x}))}
    \]
    has a $ C^{0,\alpha} $ right hand side. The H\"older conclusion of \cite[Theorem 1.1 and Remark 4.1]{liu2010interior} then gives the stated interior $ C^{2,\alpha} $ regularity, and these identities give the corresponding map regularity.

    We next prove integrability up to the boundary. Boundedness of $ b_{ij} $ and $ h\geqslant0 $ imply that, for a common $ C_0 $,
    \[
        D^2u+C_0\matr{I}\geqslant0 \quad \text{in } \regs.
    \]
    Take nested $ C^1 $ inner domains $ \regs_\varepsilon\Subset\regs $ exhausting $ \regs $, with uniformly bounded boundary areas. They can be constructed by flowing the boundary along a smooth vector field that points strictly inward. The divergence theorem and the bounded closure gradient give
    \[
        0 \leqslant \int_{\regs_\varepsilon}\tr(D^2u+C_0\matr{I})\dif\point{x} = \int_{\partial\regs_\varepsilon}Du\cdot\vec{\nu}\dif\mathcal{H}^{n-1}+nC_0|\regs_\varepsilon| \leqslant C.
    \]
    For a positive semidefinite matrix its norm is bounded by its trace. Exhaustion therefore gives $ D^2u\in L^1(\regs) $. Interior $ C^2 $ regularity identifies these classical derivatives with the distributional derivatives everywhere in the open domain. Hence $ u\in W^{2,1}(\regs) $. The same proof gives $ v\in W^{2,1}(\regt) $.

    Finally, the displayed identity for $ D\vec{T} $, boundedness of the inverse mixed matrix, and $ h\in L^1 $ give $ D\vec{T}\in L^1 $. Interior regularity again identifies the distributional derivative. Since the map is bounded, $ \vec{T}\in W^{1,1}(\regs) $. Interchanging source and target proves the assertion for $ \vec{T}^{-1} $.
\end{proof}

Fix a base transport pair $ (\point{x}_0,\overline{\point{x}}_0) $. To distinguish ambient partial derivatives from derivatives along the transport graph, write
\[
    s(\point{x}) := \tau(\point{x},\vec{T}(\point{x})), \qquad \overline{s}(\overline{\point{x}}) := \tau(\vec{T}^{-1}(\overline{\point{x}}),\overline{\point{x}}).
\]
These functions are defined where the corresponding transport pairs lie in $ U_2 $. With this neighborhood restriction understood, $ D_r=\{s<r\}\cap\regs $ and $ \vec{T}(D_r)=\{\overline{s}<r\}\cap\regt $. By \cref{lemma::thesis::localization-of-Dr}, the small sublevels stay away from the artificial boundary of $ U_2 $. Write
\[
    F(r) := \int_{D_r}\rho\dif\point{x} = \int_{\vec{T}(D_r)}\overline{\rho}\dif\overline{\point{x}},
\]
and recall that $ \chi(r)=r^{-n/2}F(r) $.
Recall the ambient fields
\[
    X^i := b^{i\overline{k}}\tau_{\overline{k}}, \qquad Y^{\overline{k}} := b^{\overline{k}i}\tau_i,
\]
restricted to the source and inverse target graphs respectively. A smooth ambient cutoff extends these fields outside the geometric neighborhood, without changing them near the closed small sublevels. In the following formulas all occurrences of $ \tau $ and its ambient derivatives are evaluated in that neighborhood. Recall the scalar identity
\[
    \tau_i X^i = \tau_{\overline{k}}Y^{\overline{k}} = a(\tau), \qquad a(r) = -\frac{(1-Kr)^2\log(1-Kr)}{K} = r-\frac{3}{2}Kr^2+O(r^3),
\]
with $ a(r)=r $ when $ K=0 $.

Denote the two nonpositive Hessian terms by
\[
    \begin{aligned}
        I_{\point{x}} :={}& (\nabla^2\tau)_{\overline{k}\overline{l}}b^{\overline{k}i}b^{\overline{l}m}h_{im}, \\
        I_{\overline{\point{x}}} :={}& (\nabla^2\tau)_{im}b^{i\overline{k}}b^{m\overline{l}}h_{\overline{k}\overline{l}} = (\nabla^2\tau)_{im}h^{im}.
    \end{aligned}
\]
The second is written as a function on the target graph. Thus
\[
    I_1(\point{x}) = I_{\point{x}}(\point{x}) + I_{\overline{\point{x}}}(\vec{T}(\point{x})).
\]
Both terms are nonpositive because the pure covariant Hessian blocks of $ \tau $ are nonpositive and the source and target matrices of $ h $ are positive definite.

Define
\[
    \begin{aligned}
        N :={}& (\tau_i-h_{im}X^m)h^{ij}(\tau_j-h_{jl}X^l) = \tau_i h^{ij}\tau_j+X^ih_{ij}X^j-2a(s), \\
        \mathcal{P}(r) :={}& \int_{\{s=r\}\cap\regs}N\frac{\rho}{|Ds|}\dif\mathcal{H}^{n-1}, \\
        \mathcal{Q}(r) :={}& \int_{\partial\regs\cap\{s<r\}}\rho\vec{X}\cdot\vec{\nu}_{\regs}\dif\mathcal{H}^{n-1}+\int_{\partial\regt\cap\{\overline{s}<r\}}\overline{\rho}\vec{Y}\cdot\vec{\nu}_{\regt}\dif\mathcal{H}^{n-1}, \\
        \mathcal{R}(r) :={}& -\int_{D_r}I_1\rho\dif\point{x}.
    \end{aligned}
\]
Here $ Ds $ is the total source derivative, and the boundary normals are outward. The level integrals are defined for almost every $ r $.

\begin{theorem}[Weak monotonicity formula]\label{theorem::thesis::weak-monotonicity}
    Under the conditions of \cref{theorem::thesis::maintheorem}, there is a common $ r_0>0 $ such that, for almost every $ 0<r<r_0 $,
    \begin{equation}\label{equation::thesis::weak-mass-inequality}
        4a(r)F'(r)+\mathcal{P}(r)+\mathcal{Q}(r)+\mathcal{R}(r) \leqslant (2n+4\delta(r))F(r),
    \end{equation}
    where, with $ \eta=1/(4n) $, one may take
    \[
        \delta(r) := C\left(r^{\alpha\eta}+r^\eta+r^{2\eta}+r\right).
    \]
    Consequently,
    \begin{equation}\label{equation::thesis::weak-monotonicity-with-dropped-terms}
        \chi'(r)+\frac{\mathcal{P}(r)+\mathcal{Q}(r)+\mathcal{R}(r)}{4a(r)r^{\frac{n}{2}}} \leqslant \frac{\beta(r)}{r}\chi(r),
    \end{equation}
    where
    \begin{equation}\label{equation::thesis::weak-monotonicity-error}
        \beta(r) := \frac{r}{a(r)}\delta(r)+\frac{n}{2}\left(\frac{r}{a(r)}-1\right), \qquad 0\leqslant\beta(r)\leqslant Cr^{\frac{\alpha}{4n}}.
    \end{equation}
    In particular, $ \int_0^{r_0}\beta(r)\dif r/r<\infty $, and
    \[
        \exp\left(-\int_0^r\frac{\beta(t)}{t}\dif t\right)\chi(r)
    \]
    is nonincreasing.
\end{theorem}
\begin{proof}
    \noindent \textbf{Step 1: Integrable vector fields and divergence identities.}
    By \cref{lemma::thesis::integrable-derivatives-for-monotonicity} and the Sobolev chain rule, the cutoff fields satisfy $ \vec{X}\in W^{1,1}(\regs) $ and $ \vec{Y}\in W^{1,1}(\regt) $. They are continuous on the closures, so their Sobolev traces agree with their continuous boundary values. On the source and target portions corresponding to $ U_2 $, the functions $ s,\overline{s} $ have integrable weak derivatives and continuous boundary values. The ambient coefficients and their derivatives are bounded on these neighborhoods.

    Here and below, $ \divr_{\Sigma}\vec{X} $ denotes the coordinate divergence $ \sum_i\partial_{x^i}[X^i(\point{x},\vec{T}(\point{x}))] $, and $ \divr_{\Sigma}\vec{Y} $ denotes $ \sum_{\overline{k}}\partial_{\overline{x}^{\overline{k}}}[Y^{\overline{k}}(\vec{T}^{-1}(\overline{\point{x}}),\overline{\point{x}})] $. The same convention applies in the cost coordinates and after blow-up.

    The interior calculation from the fifth-line estimate therefore gives, both almost everywhere and distributionally,
    \[
        \divr_{\Sigma}\vec{X}=n+E^\tau_{\point{x}}+I_{\point{x}}, \qquad \divr_{\Sigma}\vec{Y}=n+E^\tau_{\overline{\point{x}}}+I_{\overline{\point{x}}}
    \]
    on $ D_r $ and $ \vec{T}(D_r) $, respectively, for $ r<r_0 $. Recall the two errors
    \[
        \begin{aligned}
            E^\tau_{\point{x}} :={}& b^{i\overline{k}}\tau_{i\overline{k}}-n-X^i\partial_i\log|\det b_{m\overline{l}}|, \\
            E^\tau_{\overline{\point{x}}} :={}& b^{i\overline{k}}\tau_{i\overline{k}}-n-Y^{\overline{k}}\partial_{\overline{k}}\log|\det b_{m\overline{l}}|.
        \end{aligned}
    \]
    The derivatives of the determinant here are ambient derivatives. The two Hessian terms are integrable on their respective domains, since the source and target components of $ h $ are integrable and all their other coefficients are bounded. Their sum on the source graph is integrable against $ \rho $ by transport of mass.

    We use the scalar error estimate already obtained in the fifth-line calculation:
    \[
        |E^\tau_{\point{x}}|+|E^\tau_{\overline{\point{x}}}| \leqslant C\left(|\point{x}-\point{x}_0|+|\overline{\point{x}}-\overline{\point{x}}_0|+|\point{x}-\point{x}_0|\,|\overline{\point{x}}-\overline{\point{x}}_0|+r\right)
    \]
    at transport pairs with $ s\leqslant r $. By \cref{lemma::thesis::localization-of-Dr}, this is bounded by $ C(r^\eta+r^{2\eta}+r) $.

    \noindent \textbf{Step 2: The boundary terms.}
    The chain rule gives
    \[
        \partial_i s=\tau_i+h_{im}X^m, \qquad \partial_{\overline{k}}\overline{s}=\tau_{\overline{k}}+h_{\overline{k}\overline{l}}Y^{\overline{l}},
    \]
    and hence
    \[
        Ds\cdot\vec{X}=a(s)+X^ih_{ij}X^j\geqslant0, \qquad D\overline{s}\cdot\vec{Y}=a(\overline{s})+Y^{\overline{k}}h_{\overline{k}\overline{l}}Y^{\overline{l}}\geqslant0.
    \]
    For positive sufficiently small $ s $, $ a(s)>0 $, so $ Ds\ne0 $ in the interior. The analogous statement holds on the target side.

    On the domain boundaries,
    \[
        \vec{X}\cdot\vec{\nu}_{\regs}\geqslant0, \qquad \vec{Y}\cdot\vec{\nu}_{\regt}\geqslant0.
    \]
    Indeed, \cref{lemma::thesis::local-geodesic-convexity} gives nonnegative scalar products of the source and target endpoint velocities with their respective outward normals. The first-variation identities $ \vec{X}=e^{-K\sigma}\dot{\point{x}}(1) $ and $ \vec{Y}=e^{-K\sigma}\dot{\overline{\point{x}}}(1) $ therefore give these signs.

    The Sobolev divergence theorem and coarea now give, for almost every $ r $,
    \[
        \int_{D_r}\divr_{\Sigma}\vec{X}\dif\point{x} = \int_{\{s=r\}\cap\regs}\frac{Ds\cdot\vec{X}}{|Ds|}\dif\mathcal{H}^{n-1}+\int_{\partial\regs\cap\{s<r\}}\vec{X}\cdot\vec{\nu}_{\regs}\dif\mathcal{H}^{n-1},
    \]
    and the analogous target identity. More explicitly, take $ r+\varepsilon<r_0 $ and let $ \zeta_\varepsilon $ equal one below $ r $ and zero above $ r+\varepsilon $, decreasing smoothly between them. Form $ \zeta_\varepsilon(s)\vec{X} $ where the transport pair lies in $ U_2 $, and extend this product by zero to the rest of $ \regs $. By \cref{lemma::thesis::localization-of-Dr}, its support stays away from the artificial boundary of $ U_2 $. Boundedness of $ \vec{X} $, integrability of $ Ds $ on this neighborhood, and $ \vec{X}\in W^{1,1} $ show that the extended product belongs to $ W^{1,1}(\regs) $. Apply the divergence theorem to this product. The same construction applies on the target side. Letting $ \varepsilon\to0 $, coarea gives the level-set integral for almost every $ r $, while continuity of the traces gives the integral over the domain boundary. Levels for which $ \{s=r\}\cap\partial\regs $ has positive boundary measure form at most a countable set and may be excluded. 

    \noindent \textbf{Step 3: Estimate the density terms.}
    Fix a radius $ r $ for which the preceding integral identities hold. Define
    \[
        \begin{aligned}
            J(r):={}&\int_{\{s=r\}\cap\regs}
                \frac{Ds\cdot\vec{X}}{|Ds|}\dif\mathcal{H}^{n-1}
                +\int_{\partial\regs\cap\{s<r\}}
                \vec{X}\cdot\vec{\nu}_{\regs}\dif\mathcal{H}^{n-1},\\
            J_\rho(r):={}&\int_{\{s=r\}\cap\regs}
                \rho\frac{Ds\cdot\vec{X}}{|Ds|}\dif\mathcal{H}^{n-1}
                +\int_{\partial\regs\cap\{s<r\}}
                \rho\vec{X}\cdot\vec{\nu}_{\regs}\dif\mathcal{H}^{n-1}.
        \end{aligned}
    \]
    All integrands in these boundary integrals are nonnegative. Let
    \[
        \omega(r):=\sup_{\cl D_r}|\rho-\rho(\point{x}_0)|.
    \]
    Since $ J(r)=\int_{D_r}\divr_{\Sigma}\vec{X}\dif\point{x} $, subtracting the constant density $ \rho(\point{x}_0) $ gives
    \[
        \begin{aligned}
            \left|J_\rho(r)-\int_{D_r}\rho\divr_{\Sigma}\vec{X}\dif\point{x}\right| &\leqslant \omega(r)\left(J(r)+\int_{D_r}|\divr_{\Sigma}\vec{X}|\dif\point{x}\right) \\
            &=2\omega(r)\int_{D_r}(\divr_{\Sigma}\vec{X})_+\dif\point{x} \\
            &\leqslant 2\omega(r)\left(n+\sup_{D_r}|E^\tau_{\point{x}}|\right)|D_r|.
        \end{aligned}
    \]
    The last inequality uses $ I_{\point{x}}\leqslant0 $. Since $ |D_r|\leqslant\lambda F(r) $, this error is bounded by $ C\omega(r)F(r) $. The target calculation is identical, using $ \overline{\rho}(\overline{\point{x}}_0) $ and $ \vec{T}(D_r) $. The diameter estimate and H\"older continuity bound the variation of each density from its value at the base point by $ Cr^{\alpha\eta} $.

    Adding the source and target estimates, substituting the divergence formulas, and changing variables in the integral over $ \vec{T}(D_r) $, we obtain
    \[
        J_\rho(r)+J_{\overline{\rho}}(r)+\mathcal{R}(r) \leqslant (2n+4\delta(r))F(r),
    \]
    after increasing the constant in $ \delta $. Here $ J_{\overline{\rho}} $ denotes the sum of the two weighted target boundary integrals.

    \noindent \textbf{Step 4: Combine the integrals over the level sets.}
    The two integrals over the domain boundaries give $ \mathcal{Q}(r) $. For the integrals over the level sets, the change of variables and coarea formula give
    \[
        \begin{aligned}
            &\int_{\{s=r\}\cap\regs}\rho\frac{Ds\cdot\vec{X}}{|Ds|}\dif\mathcal{H}^{n-1}+\int_{\{\overline{s}=r\}\cap\regt}\overline{\rho}\frac{D\overline{s}\cdot\vec{Y}}{|D\overline{s}|}\dif\mathcal{H}^{n-1} \\
            &\qquad=\int_{\{s=r\}\cap\regs}\frac{\rho}{|Ds|}\left(2a(r)+X^ih_{ij}X^j+\tau_i h^{ij}\tau_j\right)\dif\mathcal{H}^{n-1} \\
            &\qquad=4a(r)F'(r)+\mathcal{P}(r).
        \end{aligned}
    \]
    This equality follows from the mass-preserving change of variables and the coarea formula, for almost every $ r $.

    Since $ s $ is $ C^1 $ and $ Ds\ne0 $ at interior points where $ 0<s<r_0 $, the coarea formula gives
    \[
        F'(r)=\int_{\{s=r\}\cap\regs}\frac{\rho}{|Ds|}\dif\mathcal{H}^{n-1}
    \]
    for almost every $ r $, and $ F $ is locally absolutely continuous on $ (0,r_0) $. The quantity $ \mathcal{P}(r) $ is integrable on compact radius intervals: its $ X^ih_{ij}X^j $ part is controlled by $ h\in L^1 $, and its $ \tau_i h^{ij}\tau_j $ part is controlled by the corresponding target expression with $ h\in L^1 $ and mass transport. Thus all integrals over the level sets used above are finite for almost every $ r $. Substitution in Step 3 proves \cref{equation::thesis::weak-mass-inequality}.

    \noindent \textbf{Step 5: The integrating factor.}
    Substitute $ F(r)=r^{n/2}\chi(r) $ into \cref{equation::thesis::weak-mass-inequality} and divide by $ 4a(r)r^{n/2} $. This gives \cref{equation::thesis::weak-monotonicity-with-dropped-terms,equation::thesis::weak-monotonicity-error}. The expansion of $ a $ gives $ r/a(r)-1=O(r) $. Taking the correction constant $ K\geqslant0 $ and reducing $ r_0 $ ensures $ 0<a(r)\leqslant r $. Consequently $ \beta\geqslant0 $ and $ \beta(r)\leqslant Cr^{\alpha/(4n)} $. Since $ \chi $ is locally absolutely continuous on positive radius intervals, multiplication by the integrating factor proves the monotonicity assertion.
\end{proof}

\section{Blow-up analysis}\label{section::thesis::blow-up}

With the monotonicity formula proved, we can proceed to analyse sections of $ u $ and $ v $. Like in \cite{figalli2013holder}, from now on we will consider the problem in the cost coordinates determined by $ (\point{x}_0, \overline{\point{x}}_0) $ since under this coordinate system the potential can obtain ``more convexity''. These coordinates are defined as follows:
\begin{definition}\label{definition::thesis::setting-q-coordinate}
    Fix $ (\point{x}_0, \overline{\point{x}}_0) \in \cl \Sigma $ and take $ u, v, b $ from \cref{definition::thesis::normalization-at-a-contact}.

    Write $ \matr{B}_0 := D_{\point{x}}D_{\overline{\point{x}}}b(\point{x}_0,\overline{\point{x}}_0) $, with rows indexed by $ \point{x} $ and columns by $ \overline{\point{x}} $. By \ref{condition::thesis::A2}, this matrix is invertible. Define
    \[
        \point{q}(\point{x}) := D_{\overline{\point{x}}}b(\point{x},\overline{\point{x}}_0), \qquad
        \point{p}(\overline{\point{x}}) := \matr{B}_0^{-1}D_{\point{x}}b(\point{x}_0,\overline{\point{x}}).
    \]
    Then $ \point{q}(\point{x}_0) = \point{p}(\overline{\point{x}}_0) = \point{0} $. Denote $ Q := \point{q}(\regs) $ and $ P := \point{p}(\regt) $. In these coordinates, the benefit satisfies
    \[
        D_{\point{q}}b(\point{0},\point{p})=\point{p}, \qquad
        D_{\point{p}}b(\point{q},\point{0})=\point{q},
    \]
    and hence
    \[
        D_{\point{q}}D_{\point{p}}b(\point{0},\point{p})
        =D_{\point{q}}D_{\point{p}}b(\point{q},\point{0})=\matr{I}.
    \]
    The $ c $-segments on $ \{\point{0}\}\times P $ and $ Q\times\{\point{0}\} $ are straight lines.
\end{definition}

\begin{remark}
Changing coordinates to $ (\point{q}, \point{p}) $ does not violate the conditions \ref{condition::thesis::A1} - \ref{condition::thesis::A3s}, the only distinction is that the differentiability of $ b $ may decrease by one order, for example if $ c(\point{x}, \overline{\point{x}}) \in C^4 $ then $ b(\point{q}, \point{p}) \in C^3 $, see \cite[Remark 4.2]{figalli2013holder} for more details. But we have assumed $ C^{\infty} $ for the cost, so this will not be a problem. We still do \textbf{not} claim that $ C^\infty $ is necessary for this regularity, we do not attempt to optimize cost regularity here.

After the coordinate change, the constants used in \ref{condition::thesis::A0} - \ref{condition::thesis::A3s} may change but we still use the same notation. The transformed source density is $ \rho_{\mathrm{new}}(\point{q}) = \rho_{\mathrm{old}}(\point{x}(\point{q})) |\det D_{\point{q}}\point{x}(\point{q})| $, and the target density transforms in the same way. We continue to denote them by $ \rho $ and $ \overline{\rho} $. With this understood, one can repeat the proof in \cref{subsection::thesis::weak-monotonicity} to obtain the same formula but in $ (\point{q}, \point{p}) $ coordinates. And the constants in the formula may also change but can be chosen uniformly as the base point pair ranges over a compact set. So from now on we will do the calculations entirely in the new coordinates. We continue to write $ u, v, b $ for these functions in the new coordinates.

We write $ u_i $ for $ D_{q^i}u $ and use unbarred indices for source coordinates and barred indices for target coordinates.
\end{remark}

\subsection{Convexity of the sections}

We first prove the convexity of sections of $ u $ and $ v $ in the new coordinates. For this we first prove a lemma
\begin{lemma}\label{lemma::thesis::positiveness-of-D2b-DbDb}
    Under the $ \point{q}, \point{p} $ coordinate,
    \[
        D^2_{\point{q}\point{q}} b \geqslant -k D_{\point{q}} b \otimes D_{\point{q}} b, \quad D^2_{\point{p}\point{p}} b \geqslant -k D_{\point{p}} b \otimes D_{\point{p}} b, 
    \]
\end{lemma}
\begin{proof}
    Pick any vector $ \vec{v} $. Take a geodesic $ (\point{q}(t), \point{p}) $, with $ \point{q}(0) = \point{0} $, $ \point{q}(1) = \point{q} $ and $ D_{\point{p}} b(\point{q}(t), \point{p}) = t D_{\point{p}} b(\point{q}, \point{p}) $. 
    Differentiating the defining equation of the geodesic segment gives
    \[
        b_{i\overline{k}}\dot q^i
        =b_{\overline{k}}(\point{q},\point{p}),
        \qquad
        b_{i\overline{k}}\ddot q^i
        +b_{ij\overline{k}}\dot q^i\dot q^j=0.
    \]
    Construct
    \[
        m(t) := \vec{v}^{\top} D^{2}_{\point{p} \point{p}} b(\point{q}(t), \point{p}) \vec{v},
    \]
    Differentiating with respect to $ t $,
    \[
        \begin{aligned}
            m^{\prime}(t)
            &=
            b_{i\overline{k}\overline{l}}
            \dot q^i v^{\overline{k}}v^{\overline{l}},
            \\
            m^{\prime\prime}(t)
            &=
            \left(
                b_{ij\overline{k}\overline{l}}
                -
                b_{ij\overline{a}}
                b^{\overline{a}m}
                b_{m\overline{k}\overline{l}}
            \right)
            \dot q^i\dot q^j
            v^{\overline{k}}v^{\overline{l}}
            \\
            &\geqslant
            -K\left(
                b_{i\overline{k}}\dot q^i
                v^{\overline{k}}
            \right)^2
            \\
            &=
            -K\left(
                \vec{v}\cdot
                D_{\point{p}}b(\point{q},\point{p})
            \right)^2,
        \end{aligned}
    \]
    
    The normalization in \cref{definition::thesis::normalization-at-a-contact} gives $ b(\point{0},\point{p}) = 0 $, hence $ m(0) = 0 $. Moreover, since $ D_{\point{q}}b (\point{0}, \point{p})$ is a linear function of $ \point{p} $, $ b_{i\overline{k}\overline{\ell}}(\point{0},\point{p}) = 0 $ and therefore $ m^{\prime}(0) = 0 $.
    Integrating twice, we obtain
    \[
        \begin{aligned}
            \vec{v}^{\top}
            D^2_{\point{p}\point{p}}b(\point{q},\point{p})
            \vec{v}
            &=m(1)
            =\int_0^1(1-t)m^{\prime\prime}(t) \dif t
            \\
            &\geqslant
            -\frac K2
            \left(
                \vec{v}\cdot
                D_{\point{p}}b(\point{q},\point{p})
            \right)^2.
        \end{aligned}
    \]
    Taking $ k = K / 2 $ proves the second inequality. The first inequality follows by interchanging $ \point{q} $ with $ \point{p} $.
\end{proof}

Using this estimate, we can convexify the potential functions.
\begin{proposition}\label{proposition::thesis::convexification}
    Define 
    \[
        w(\point{q}) := \frac{e^{k u(\point{q})} - 1}{k}, \quad z(\point{p}) := \frac{e^{k v(\point{p})} - 1}{k},
    \]
    if $ k = 0 $ then just define $ w := u $ and $ z := v $. Then $ w $ and $ z $ are convex. Moreover, their Monge--Amp\`ere measures satisfy $ \det D^2w\geqslant c>0 $ on $ Q $ and $ \det D^2z\geqslant c>0 $ on $ P $, with $ c $ independent of the base transport pair.
\end{proposition}
\begin{proof}
    Calculating the second derivative of $ w $ gives
    \[
        D^2 w = e^{k u} \left[ D^2 u + k D u \otimes D u \right],
    \]
    while $ D^2 u = D^2 b + H $ where $ H = (h_{ij}) $ is positive definite, thus
    \[
        D^2 u + k D u \otimes D u = H + D^2 b + k D u \otimes D u \geqslant H \geqslant 0.
    \]
    Thus $ w $ is convex.

    The lower bound for the Monge--Amp\`ere measure follows from the fact that the Monge-Amp\'ere equation for $ u $ shows that
    \[
        \lambda^{-1} \leqslant \det \left[ u_{ij} - b_{ij} \right] = \det h_{ij} \leqslant \lambda,
    \]
    and $ D^2 w \geqslant H $ so $ \det D^2 w \geqslant \lambda^{-1} $.

    The same results for $ z $ follow by a dual argument.
\end{proof}

Thus the sections of $ w $ are convex. Since
\[
    S_h[u] = \{ \point{q}\in Q \colon u(\point{q}) < h \} = S_{(e^{kh} - 1) / k}[w],
\]
the sections of $ u $ are the same as those of $ w $ (with a change of height). 
\begin{remark}\label{remark::thesis::equivalence-of-sections-of-uw}
Boundedness of the original potentials and benefit gives a common bound $ 0\leqslant u,v\leqslant M $ for all normalized base pairs. Hence
\[
    u \leqslant \frac{e^{ku} - 1}{k} \leqslant e^{kM} u,
\]
Thus the sections of $ u $ and $ w $ are comparable after multiplying their heights by fixed constants. The same holds for $ v $ and $ z $. After blow-up, $ s_j\hat{v}\leqslant M $ gives $ \hat{v}\leqslant\hat{z}\leqslant e^{kM}\hat{v} $, and similarly for $ \hat{u},\hat{w} $. We use these comparisons to apply convexity estimates to the sections.
\end{remark}

\begin{lemma}[Comparison of sections and radius sublevels]\label{lemma::thesis::section-radius-comparison}
    For all sufficiently small $ r > 0 $,
    \begin{equation}\label{equation::thesis::section-radius-inclusions}
        \frac{1}{4}S_r[u] \subset D_r \subset S_{2r}[u],
    \end{equation}
    where $ D_r $ is defined in \cref{lemma::thesis::localization-of-Dr}. The same holds for $ v $.
\end{lemma}
In \cite[Proposition 2.6]{collins2025boundary}, they have $ \frac{1}{2}S_r[u]\subset D_r\subset S_r[u] $ by using the convexity of $ u $. We obtain a similar containment using the convexification of $ u $.

\begin{proof}
Use the convexification from \cref{proposition::thesis::convexification},
\[
    w(\point{q})=\frac{e^{ku(\point{q})}-1}{k},
\]
with $ w=u $ when $ k=0 $. By \cref{lemma::thesis::uniform-positive-contact-gap} and continuity of the transport map, taking $ r $ small places the graph over $ S_r[u] $ in the neighborhood where the estimates in \cref{lemma::thesis::positiveness-of-sigma} hold and the local sublevels $ D_r $ are defined as in \cref{lemma::thesis::localization-of-Dr}.

Take $ \point{q}\in S_r[u] $ and put $ \point{z}:=\point{q}/4 $. Convexity of the coordinate domain puts $ \point{z} $ in that domain, and convexity of $ w $ gives $ w(\point{z})\leqslant w(\point{q})/4 $, so $ u(\point{z})<r $. The supporting inequality at $ \point{z} $ gives
\[
    3\point{z}\cdot Dw(\point{z})
    \leqslant w(\point{q})-w(\point{z})
    \leqslant w(\point{q}).
\]
Using the third estimate in \cref{lemma::thesis::positiveness-of-sigma} and $ Du=e^{-ku}Dw $, we obtain
\[
    \tau(\point{z},\vec{T}(\point{z}))
    \leqslant \point{z}\cdot Du(\point{z})
    \leqslant \frac{w(\point{q})}{3}
    < \frac{e^{kr}-1}{3k}
    < r
\]
for sufficiently small $ r $. For $ k=0 $, the second to last term is $ r/3 $. This proves the first inclusion. For the second, \cref{lemma::thesis::positiveness-of-sigma} gives on $ D_r $
\[
    u\leqslant\tau+C\tau^2<r+Cr^2\leqslant2r
\]
after reducing $ r $ so that $ Cr\leqslant1 $. The target argument is identical.
\end{proof}

\subsection{Preparation for blow-up}

We first prove an $ h^{\frac{n}{2}} $ volume bound for $ S_h[u] $.
\begin{lemma}\label{lemma::thesis::estimate-of-sections-and-limit-of-chi}
    There exist $ h_0,C_1,C_2>0 $ such that, for $ 0<h<h_0 $,
    \[
        C_1 h^{\frac n2}\leqslant |S_h[u]|\leqslant C_2 h^{\frac n2}.
    \]
    These constants can be chosen independently of the base transport pair. Moreover,
    \[
        0<\lim_{r\to0}\chi(r)=\chi(0)<\infty.
    \]
\end{lemma}
\begin{proof}
    We first prove the upper bound for $ S_h[w] $. Write its John ellipsoid as $ \point{z}_h+\matr{A}_hB_1 $, choosing $ \matr{A}_h $ symmetric positive definite, and put $ E_h:=\point{z}_h+\frac14\matr{A}_hB_1 $. For $ \point{q}\in E_h $, $ \vec{\xi}\in\partial w(\point{q}) $, and any unit vector $ \vec{e} $, the points $ \point{q}\pm\frac12\matr{A}_h\vec{e} $ belong to $ S_h[w] $. Since $ w\geqslant0 $, the supporting inequalities give
    \[
        |\vec{\xi}\cdot\matr{A}_h\vec{e}|\leqslant2h,
        \qquad \matr{A}_h\partial w(E_h)\subset B_{2h}.
    \]
    By the lower bound for the Monge--Amp\`ere measure from \cref{proposition::thesis::convexification},
    \[
        c\det\matr{A}_h
        \leqslant\lambda^{-1}|E_h|
        \leqslant|\partial w(E_h)|
        \leqslant|\matr{A}_h^{-1}B_{2h}|
        =\frac{|B_1|(2h)^n}{\det\matr{A}_h}.
    \]
    Hence $ \det\matr{A}_h\leqslant Ch^{n/2} $. John's theorem gives
    \[
        |S_h[w]|\leqslant n^n|B_1|\det\matr{A}_h\leqslant Ch^{n/2}.
    \]
    Since $ S_h[u]=S_{(e^{kh}-1)/k}[w] $ and $ (e^{kh}-1)/k\leqslant Ch $ for bounded $ h $, the same upper bound holds for $ S_h[u] $, with the usual interpretation when $ k=0 $.

    Choose a common radius $ r_\ast>0 $ within the range of the monotonicity formula and \cref{equation::thesis::section-radius-inclusions}. We first give a uniform positive lower bound for $ \chi(r_\ast) $. Boundedness of the original potential and benefit gives a common bound $ 0\leqslant w\leqslant M $ for the normalized convexifications. Choose $ 0<\lambda_0<1 $ such that $ \lambda_0M<(e^{kr_\ast}-1)/k $. Convexity and $ w(\point{0})=0 $ give
    \[
        \lambda_0Q\subset S_{r_\ast}[u], \qquad
        \frac{\lambda_0}{4}Q\subset D_{r_\ast}.
    \]
    Here, when $ k=0 $, the threshold is $ r_\ast $. The volume $ |Q| $ has a common positive lower bound under these coordinate changes, since the Jacobian of the $ \point{q} $-coordinate map and its inverse are bounded in terms of the cost norms and the nondegeneracy bound. The transformed density also has a common positive lower bound, so
    \[
        \chi(r_\ast)=r_\ast^{-n/2}\int_{D_{r_\ast}}\rho\dif\point{q}\geqslant c_\ast>0.
    \]
    By the monotonicity formula, for some common $ C,\gamma>0 $ and $ 0<r\leqslant r_\ast $,
    \[
        \chi(r)\geqslant e^{C r^\gamma-C r_\ast^\gamma}\chi(r_\ast)
        \geqslant e^{-C r_\ast^\gamma}c_\ast.
    \]
    The upper bound $ \rho\leqslant\lambda $ therefore implies
    \[
        |D_r|\geqslant\lambda^{-1}r^{n/2}\chi(r)\geqslant c r^{n/2}.
    \]
    Using $ D_{h/2}\subset S_h[u] $ proves the lower bound for $ |S_h[u]| $ after choosing $ h_0 $ sufficiently small.

    Finally, $ D_r\subset S_{2r}[u] $ and the section upper bound give
    \[
        \chi(r)=r^{-n/2}\int_{D_r}\rho\dif\point{q}
        \leqslant\lambda r^{-n/2}|S_{2r}[u]|\leqslant C.
    \]
    The corrected function $ e^{-Cr^\gamma}\chi(r) $ is nonincreasing and bounded above and below by positive constants. It therefore has a positive finite limit at zero. Since $ e^{-Cr^\gamma}\to1 $, the same holds for $ \chi(r) $.
\end{proof}

We now define the blow-up used here. Volume estimates from \cref{lemma::thesis::estimate-of-sections-and-limit-of-chi} help us determine the coefficient used in blow-up.
\begin{definition}\label{definition::thesis::definition-of-blow-up}
    Choose transport pairs $ (\point{x}_j, \overline{\point{x}}_j) \in \cl\Sigma $ and a sequence $ s_j \to 0 $. For every pair apply the $ (\point{q}, \point{p}) $-coordinate transform defined in \cref{definition::thesis::setting-q-coordinate} with respect to $ (\point{x}_j, \overline{\point{x}}_j) $. Denote the resulting functions by $ u_j, v_j, b_j, \tau_j$, and the coordinate domains and transport map by $ Q_j^0, P_j^0, \point{T}_j^0$.
    
    Fix a constant $ L > 2 $, let $ S_{Ls_j}[u_j] = \{ u_j < Ls_j \} $, denote its John ellipsoid by $ \point{z}_j + \matr{A}_j B_1 $, where $ \point{z}_j $ is its center and $ \matr{A}_j $ is symmetric positive definite. Denote the ellipsoid $ \mathcal{E}_j := \matr{A}_j B_1 $, which is the John ellipsoid translated to the origin. Define the normalized potentials and sections by
    \[
        \begin{gathered}
            u^{(j)}(\point{q}) := \frac{u_j(\matr{A}_j \point{q})}{s_j}, \quad v^{(j)}(\point{p}) := \frac{v_j(s_j \matr{A}_j^{-\top} \point{p})}{s_j}\\
            K_j := \matr{A}_j^{-1} S_{Ls_j}[u_j]
        \end{gathered}
    \]
    The normalization uses the matrix of the John ellipsoid and keeps the base point at $ \point{0} $.

    After deciding the transformation, we have
    \[
        \begin{gathered}
            b^{(j)}(\point{q}, \point{p}) := s_j^{-1} b_j(\matr{A}_j \point{q}, s_j \matr{A}_j^{-\top} \point{p}), \\
            \tau^{(j)}(\point{q}, \point{p}) := s_j^{-1} \tau_j(\matr{A}_j \point{q}, s_j \matr{A}_j^{-\top} \point{p}), \\
            a^{(j)}(r) := \frac{a(s_j r)}{s_j}.
        \end{gathered}
    \]
    And the transport mapping
    \[
        \point{T}^{(j)}(\point{q}) := s_j^{-1}\matr{A}_j^\top \point{T}_j^{0}(\matr{A}_j \point{q}).
    \]

    Let $ Q_j := \matr{A}_j^{-1} Q_j^0 $ and $ P_j := s_j^{-1} \matr{A}_j^{\top} P_j^0$ be the normalized regions. For the monotonicity formula, let $ D_{j, r}^0 $ be the unscaled set defined by $ \tau_j $ and $ \rho_j^0,\overline{\rho}_j^0 $ be the source and target densities in the corresponding $ (\point{q}, \point{p}) $-coordinates. Let
    \[
        F_j^{0}(r) = \int_{D_{j,r}^0}\rho_j^0 \dif \point{q}, \quad \chi_j^0(r) = r^{-\frac{n}{2}} F_j^0(r).
    \]
    
    After blow-up denote $ D^{(j)}_r = \matr{A}_j^{-1} D_{j, s_j r}^0 $, and
    \[
        F^{(j)}(r) = \frac{F_j^0(s_j r)}{\det \matr{A}_j}, \quad \chi^{(j)}(r) = r^{-\frac{n}{2}} F^{(j)}(r).
    \]
    The transformed densities are
    \[
        \rho^{(j)}(\point{q}) := \rho_j^0(\matr{A}_j\point{q}), \qquad
        \overline{\rho}^{(j)}(\point{p})
        :=\frac{s_j^n}{(\det\matr{A}_j)^2}
        \overline{\rho}_j^0(s_j\matr{A}_j^{-\top}\point{p}).
    \]
    With these definitions, $ \point{T}^{(j)} $ transports $ \rho^{(j)}\dif\point{q} $ to $ \overline{\rho}^{(j)}\dif\point{p} $. The section volume bounds give $ c s_j^{n/2}\leqslant\det\matr{A}_j\leqslant C s_j^{n/2} $, so both transformed densities have positive lower and upper bounds independent of $ j $.
\end{definition}

\begin{remark}
    We chose to use superscript $ ^{(j)} $ because subscripts may be confused with derivatives, which we often use. But in the following proof, writing superscripts will make the formulae messy. Therefore we introduce a new notation: adding a $ \hat{} $ means the term varies with $ j $, for example
    \[
        \hat{\tau}_i \hat{b}^{i \overline{k}} \hat{\tau}_{\overline{k}}
    \]
    is equivalent to
    \[
        (\tau^{(j)})_i (b^{(j)})^{i \overline{k}}(\tau^{(j)})_{\overline{k}}.
    \]
    For constants and matrices that use subscripts, since there is no risk of confusion with derivatives, we shall use subscripts together with the hat notation.

    In sections concerning the blow-up analysis, we shall use this hat notation.
\end{remark}

\subsection{Homogeneity of blow-up limit}

In this subsection we prove that the blow-up limit is similar to the case in \cite{collins2025boundary}, a homogeneous function. Homogeneity in \cite{collins2025boundary} is proved by noticing that in the blow-up picture, the monotonic function $ \chi $ becomes a constant and thus the dropped nonpositive term (i.e. the first line and third line in \cref{equation::thesis::formula-de-derivative-of-monotonicity}, the fourth line in our formula vanishes in their case) becomes zero, the square in the first line gives homogeneity and the inner product with boundary normal gives conic shape of the blow-up region. The idea of this subsection is similar, however since we do not have $ \tau_i(\point{q}, \point{p}) \to \point{p} $ directly, we have to prove a similar result. To be precise we show that $ D_{\point{q}} u \cdot \point{q} - \tau \longrightarrow 0 $ in distribution, with this the vanishing dropped terms give the homogeneity we want.

The subsection is organized as follows. We first prove the existence of blow-up limits and make definitions. Then we show why and how the dropped terms tend to zero in our case. Finally we use those results to prove the properties of the limit function $ u_\infty $.

\paragraph{Notation}

We first define the convergence limits.

\begin{proposition}[Blow-up limits]\label{proposition::thesis::limit-existence-of-blow-up}
    Take $ s_j, \hat{u}, Q_j $ from a blow-up sequence defined in \cref{definition::thesis::definition-of-blow-up}, and denote by 
    \[
        \hat{w} := \frac{e^{k s_j \hat{u}} - 1}{k s_j}
    \]
    the convexification of $ \hat{u} $ as in \cref{proposition::thesis::convexification}. Then after passing to a subsequence, there is an open convex set $ K_\infty $  such that
    \[
        \cl K_j \to \cl K_{\infty}
    \]
    in Hausdorff sense, and a convex function $ 0 \leqslant u_\infty < \infty $ defined on $ K_\infty $ satisfying
    \[
        \hat{u} \longrightarrow u_\infty, \quad \hat{w} \longrightarrow u_\infty
    \]
    locally uniformly. Moreover,
    \[
        D \hat{u} \longrightarrow D u_\infty, \quad D \hat{w} \longrightarrow D u_\infty
    \]
    pointwise at differentiability points of $ u_\infty $ and in $ L^1_{\mathrm{loc}}(K_\infty) $.
\end{proposition}
\begin{proof}
    The existence of $ K_\infty $ comes from Blaschke's selection theorem and the fact that $ K_j $ are bounded by $ B_{2n}(\point{0}) $ by normalization. The theorem gives $ \cl K_j\to K $ for a compact convex body $ K $. The normalized inner John ellipsoid of $ K_j $ ensures that $ K_\infty:=\interior K $ is nonempty, and $ K=\cl K_\infty $.

    Take any compact set $ G \subset K_\infty $, since $ G $ is compact and $ K_\infty $ is open, there exists a $ \delta $ such that $ \{ \point{q} \colon \dist(\point{q}, G) < 2\delta \} \subset K_{\infty} $. By convexity of $ K_j $ and the Hausdorff convergence $ \cl K_j \to \cl K_\infty $, there exists $ \delta > 0 $ such that, for all sufficiently large $ j $, $ G \subset \{ \point{q} \colon \dist(\point{q}, G) < \delta \} \Subset K_j $. Therefore we can talk about convergence over any compact subset of $ K_\infty $.
    
    Since $ K_j $ is a normalization of $ S_{L s_j}[u] $, we have $ \hat{u} \leqslant L $ on $ K_j $, and thus $ \hat{w} \leqslant 2L $ on $ K_j $. Take a compact subset $ G \subset K_\infty $, for $ j $ sufficiently large $ \hat{w} $ is defined on $ G $. Since they are uniformly bounded, we can extract a subsequence that converges to a convex function. Now apply Arzel\'a--Ascoli on a compact exhaustion of $ K_\infty $ and taking a diagonal subsequence, we have for a subsequence $ \hat{w} \to u_\infty $ for some $ u_\infty $ convex locally uniformly on $ K_\infty $.

    For the gradient of $ \hat{w} $, by its convexity and the $ \delta $ margin we have left around $ G $,
    \[
        \frac{\hat{w}(\point{q}) - \hat{w}(\point{q} - \delta \vec{e})}{\delta} \leqslant D \hat{w}(\point{q}) \cdot \vec{e} \leqslant \frac{\hat{w}(\point{q} + \delta \vec{e}) - \hat{w}(\point{q})}{\delta},
    \]
    let $ j \to \infty $ and then $ \delta \to 0 $ at a differentiability point of $ u_\infty $; both bounds converge to $ D u_\infty \cdot \vec{e} $, so $ D\hat{w} \to D u_\infty $. Since $ u_\infty $ is convex, it is differentiable almost everywhere. On each compact subset with a $ \delta $ margin, we have $ |D \hat{w}| \leqslant 2L / \delta $ uniformly bounded, so
    \[
        \int_{G} |D \hat{w} - D u_\infty| \to 0
    \]
    follows from the dominated convergence theorem. This proves $ L^1_{\mathrm{loc}}(K_\infty) $ convergence of $ D \hat{w} $.

    For the convergence of $ \hat{u} $, notice that $ \hat{w} - \hat{u} \to 0 $ and $ D \hat{w} = e^{k s_j \hat{u}} D\hat{u} $ with $ e^{k s_j \hat{u}} \to 1 $, the result follows.
\end{proof}

A direct corollary from the proof of the proposition is 
\begin{lemma}\label{lemma::thesis::bound-for-Du-on-arbitrary-G}
    Take a compact subset $ G \subset \interior K_{\infty} $, then there are constants $ C_G > 0 $, $ j_{G} > 0 $. As long as $ j > j_G $, we have
    \[
        \sup_{G} |D \hat{u}|, \sup_{G} |D \hat{w}| \leqslant C_{G}.
    \]
    And there is $ r_G > 0 $ such that
    \[
        G \subset \hat{D}_{r_G}.
    \]
\end{lemma}
\begin{proof}
    The bound for $ D \hat{u} $ follows from the proof of \cref{proposition::thesis::limit-existence-of-blow-up}. By \cref{lemma::thesis::positiveness-of-sigma}, $ 0 \leqslant \hat{\tau} \leqslant \point{q} \cdot D \hat{u} $ on $ G $, so a sufficiently large $ r_G $ gives $ G \subset \hat{D}_{r_G} $.
\end{proof}

\paragraph{The dropped terms}

Now we start with the dropped term in the monotonicity formula. We first write the change of those terms after blow-up explicitly,
\begin{proposition}\label{proposition::thesis::remaining-terms}
    Under the blow-up from \cref{definition::thesis::definition-of-blow-up}, recall that $ \hat{h}_{ij} = \hat{u}_{ij} - \hat{b}_{ij} $, and
    \[
        \hat{X}^i = \hat{b}^{i \overline{k}} \hat{\tau}_{\overline{k}}, \hat{Y}^{\overline{k}} = \hat{b}^{\overline{k} i} \hat{\tau}_{i}.
    \]

    Recall $ N,\mathcal{P},\mathcal{Q},\mathcal{R} $ from \cref{subsection::thesis::weak-monotonicity}. We write $ \mathcal{P}_j,\mathcal{Q}_j,\mathcal{R}_j $ for their values after blow-up. Their transformations are as follows:
    \begin{enumerate}
        \item \textbf{The first line, integral of a square on the level set boundary.} Recall that
        \[
            N = (\tau_i h^{im} - \tau_{\overline{k}} b^{\overline{k} m}) h_{mn} (\tau_j h^{jn} - \tau_{\overline{l}} b^{\overline{l} n}) = \tau_i h^{ij} \tau_j + X^i h_{ij} X^j - 2 a(\tau)
        \]
        as part of the integrand in the first line of \cref{equation::thesis::formula-de-derivative-of-monotonicity}. After blow-up,
        \[
            \hat{N} = \hat{\tau}_i \hat{h}^{ij} \hat{\tau}_j + \hat{X}^i \hat{h}_{ij} \hat{X}^j - 2 \hat{a}(\hat{\tau}),
        \]
        which is $ s_j^{-1} $ times the value of $ N $ evaluated at the corresponding transformed points. The corresponding integral is
        \[
            \mathcal{P}_j(r) = \int_{\{\hat{\tau}(\point{q}, \hat{\point{T}}(\point{q})) = r \}} \hat{N} \frac{\hat{\rho}(\point{q})}{|\partial_{\point{q}} \hat{\tau}(\point{q}, \hat{\point{T}}(\point{q}))|} \dif \mathcal{H}^{n - 1}.
        \]
        {Here $ \partial_{\point{q}} $ denotes differentiation of $ \hat{\tau}(\point{q},\hat{\point{T}}(\point{q})) $. By the coarea formula and the change of variables in the volume integral, for every $ 0<a<b $, $ \int_a^b\mathcal{P}_j(r)\dif r=(s_j\det\matr{A}_j)^{-1}\int_{s_ja}^{s_jb}\mathcal{P}(t)\dif t $. Hence $ \mathcal{P}_j(r)=\mathcal{P}(s_jr)/\det\matr{A}_j $ for almost every $ r $.}
        \item \textbf{The third line, boundary normal term.} Recall that
        \[
            \mathcal{Q}(r) = \int_{\partial Q \cap \{ \tau < r \}} \rho \vec{X} \cdot \vec{\nu}_{Q} \dif \mathcal{H}^{n - 1} + \int_{\partial P \cap \{ \tau < r \}} \overline{\rho} \vec{Y} \cdot \vec{\nu}_{P} \dif \mathcal{H}^{n - 1}.
        \] 
        After blow-up, with $ \vec{\nu}_{Q_j},\vec{\nu}_{P_j} $ denoting the outward unit normals, 
        \[
            \mathcal{Q}_j(r) = \int_{\partial Q_j \cap \{ \hat{\tau} < r \}} \hat{\rho} \hat{\vec{X}} \cdot \vec{\nu}_{Q_j} \dif \mathcal{H}^{n - 1} + \int_{\partial P_j \cap \{ \hat{\tau} < r \}} \hat{\overline{\rho}} \hat{\vec{Y}} \cdot \vec{\nu}_{P_j} \dif \mathcal{H}^{n - 1},
        \]
        similar to $ \mathcal{P} $, $ \mathcal{Q}_j(r) = \mathcal{Q}(s_j r) / \det \matr{A}_j $ after the change of coordinates.
        \item \textbf{The fourth line, trace of Hessian term from $ I_1 $.} After blow-up transform $ I_1 $ is
        \[
            \hat{I}_1 = (\nabla^2_{\point{q}\point{q}} \hat{\tau})_{ij} \hat{h}^{ij} + (\nabla^2_{\point{p}\point{p}} \hat{\tau})_{\overline{k} \overline{l}} \hat{b}^{\overline{k} i} \hat{b}^{\overline{l} j} \hat{h}_{ij}.
        \]
        which has the same value after the change of coordinates. And the integral with respect to the density is
        \[
            \mathcal{R}_j(r) = -\int_{\hat{D}_r} \hat{I}_1 \hat{\rho} \dif \point{q},
        \]
        and similarly $ \mathcal{R}_j(r) = \mathcal{R}(s_j r) / \det \matr{A}_j $ after the change of coordinates. 
    \end{enumerate}

    Finally, $ \hat{\chi}(r) $ satisfies
    \begin{equation}\label{equation::thesis::formula-j-with-dropped-terms}
        \hat{\chi}^{\prime}(r) + \frac{1}{4 \hat{a}(r) r^{\frac{n}{2}}} \left[ \mathcal{P}_j + \mathcal{Q}_j + \mathcal{R}_j \right] \leqslant \frac{\beta(s_j r)}{r} \hat{\chi}(r).
    \end{equation}
    It holds on each fixed compact interval for all sufficiently large $ j $.
\end{proposition}

With this notation, we prove that these nonnegative terms tend to zero under blow-up.
\begin{lemma}\label{lemma::thesis::limit-chi-picture-and-zero-negative-terms}
Assume that the monotonicity inequality applies to $ \hat{D}_r $ for all sufficiently large $ j $. The constants in the upper bound for $ \beta(r) $ in \cref{equation::thesis::weak-monotonicity-error} depend only on the norms of $ c,\rho,\overline{\rho} $, which are uniformly bounded in the cost coordinates as the base pair varies. We may therefore use a common upper bound, still denoted by $ \beta(r) $, satisfying
\[
    \int_0^r \frac{\beta(t)}{t} \dif t < \infty
\]
and 
\[
    \hat{\chi}^{\prime}(r) \leqslant \frac{\beta(s_jr)}{r} \hat{\chi}(r)
\]
on each fixed compact interval for all sufficiently large $ j $. 
Assume for every fixed $ 0 < a < b $,
\[
    \frac{\hat{\chi}(a)}{\hat{\chi}(b)} \longrightarrow 1.
\]
\textbf{Conclusion: } After passing to a further subsequence, there exists a constant $ M > 0 $ such that for any $ r $
\[
    \lim_{j \to \infty} \hat{\chi}(r) = M
\]
locally uniformly for $ r > 0 $. 

Besides, using the notation from \cref{proposition::thesis::remaining-terms}, we have
\begin{enumerate}
    \item For any $ 0 < a < b $, 
    \[
        \int_{a}^{b} \mathcal{P}_j(t) \dif t = \int_{\{a < \hat{\tau} < b\}} \hat{N} \hat{\rho} \dif \point{q} \longrightarrow 0;
    \]
    \item For every fixed $ r > 0 $, 
    \[
        \mathcal{Q}_j(r), \mathcal{R}_j(r) \longrightarrow 0.
    \]
    Note that in $ \mathcal{R}_j $, the two parts in $ \hat{I}_1 $ are both nonpositive, thus their separate integrals also tend to zero. We will use this property in the proof of \cref{lemma::thesis::integral-of-x-Du-tau-tends-to-zero}.
\end{enumerate}
\end{lemma}
\begin{proof}
    The proof is done in three steps. First we prove the limit of $ \hat{\chi}(r) $ with respect to $ j $ exists for any $ r $. Then we prove their limits are the same. Finally we prove the three terms tend to zero in different senses.

    \noindent \textbf{Step 1.} By definition of the blow-up and volume bounds from \cref{lemma::thesis::estimate-of-sections-and-limit-of-chi}, we have
    \[
        c r^{\frac{n}{2}} \leqslant |S_r[\hat{u}]| \leqslant C r^{\frac{n}{2}},
    \]
    and by upper and lower bounds of the density,
    \[
        \lambda^{-1} |\hat{D}_r| \leqslant \hat{F}(r) \leqslant \lambda |\hat{D}_r|.
    \]
    By \cref{equation::thesis::section-radius-inclusions}, $ \frac{1}{4} S_r[\hat{u}] \subset \hat{D}_r \subset S_{2r}[\hat{u}] $, so we have $ \hat{\chi}(r) = r^{-\frac{n}{2}} \hat{F}(r) $ is then bounded from above and below by constants independent of $ j $. Fix $ r_\ast > 0 $. After passing to a subsequence, the uniform upper and positive lower bounds give
    \[
        \hat{\chi}(r_\ast) \longrightarrow M
    \]
    for some $ 0 < M < \infty $.

    \noindent \textbf{Step 2.} By the ratio assumption, for every fixed $ r > 0 $,
    \[
        \frac{\hat{\chi}(r)}{\hat{\chi}(r_\ast)} \longrightarrow 1.
    \]
    Therefore, along the same subsequence,
    \[
        \hat{\chi}(r) \longrightarrow M.
    \]

    For uniformity, set $ I_j(r) := \int_0^{s_j r} \frac{\beta(t)}{t} \dif t $, and $ \mathcal{M}_j(r) := e^{-I_j(r)} \hat{\chi}(r) $. We have by the monotonicity formula that $ \mathcal{M}_j $ is nonincreasing. For any fixed $ 0 < a < b $, 
    \[
        \sup_{r \in [a, b]} I_j(r) \leqslant I_j(b) \longrightarrow 0
    \]
    uniformly as $ j \to \infty $, so since $ \hat{\chi}(a) $ and $ \hat{\chi}(b) $ tends to the same limit, so does $ \mathcal{M}_j(a) $ and $ \mathcal{M}_j(b) $. And by monotonicity, for $ a \leqslant r \leqslant b $ we have $ \mathcal{M}_j(a) \geqslant \mathcal{M}_j(r) \geqslant \mathcal{M}_j(b) $, which proves uniform convergence of $ \mathcal{M}_j $. Thus $ \hat{\chi} $ also converges to $ M $ uniformly.

    \noindent \textbf{Step 3.} Multiply \cref{equation::thesis::formula-j-with-dropped-terms} by $ e^{-I_j(r)} $ and integrate we have 
    \[
        0 \leqslant \int_{a}^{b} \frac{e^{-I_j(r)}}{4\hat{a}(r) r^{n / 2}} \left[ \mathcal{P}_j + \mathcal{Q}_j + \mathcal{R}_j \right] \dif r \leqslant \mathcal{M}_j(a) - \mathcal{M}_j(b) \longrightarrow 0 
    \]
    as $ j \to \infty $. Since the coefficient $ \frac{e^{-I_j(r)}}{4\hat{a}(r) r^{n / 2}} $ is bounded from above and below by positive constants independent of $ j $, we have
    \[
        \int_{a}^{b} \mathcal{P}_j + \mathcal{Q}_j + \mathcal{R}_j \dif r \longrightarrow 0.
    \]

    For $ \mathcal{P}_j $, by coarea
    \[
        \int_{a}^{b} \mathcal{P}_j \dif r = \int_{\{ a < \hat{\tau} < b \}}  \hat{N} \hat{\rho} \dif \point{q} \longrightarrow 0.
    \]

    For $ \mathcal{Q}_j $ and $ \mathcal{R}_j $, since their integration regions are monotonically increasing with radius, and their integrands are nonnegative, both functions are nondecreasing with $ r $. So for any $ r $ fix any $ b > r $
    \[
        b\left[\mathcal{Q}_j(r) + \mathcal{R}_j(r) \right] \leqslant \int_{b}^{2b} (\mathcal{Q}_j + \mathcal{R}_j) \dif r \longrightarrow 0.
    \]
    therefore $ \mathcal{Q}_j(r), \mathcal{R}_j(r) \longrightarrow 0 $.
\end{proof}

We now know the integral of $ \hat{N} $ over any strip converges to zero, but we cannot judge the convergence if the radius goes to zero. The next lemma addresses this problem.

The assumptions of \cref{lemma::thesis::limit-chi-picture-and-zero-negative-terms}, including the ratio hypothesis, remain in force through \cref{lemma::thesis::2-homogeneity-of-limit}.

\begin{lemma}\label{lemma::thesis::int-Nj->0}
    We have for any fixed $ r > 0 $, 
    \[
        \lim_{j \to \infty} \int_{\hat{D}_r} \hat{N} \dif \point{q} = 0.
    \]
\end{lemma}
\begin{proof}
    By \cref{lemma::thesis::limit-chi-picture-and-zero-negative-terms} and the uniform density bounds, the integral of $ \hat{N} $ over every fixed strip $ \{a<\hat{\tau}<b\} $ tends to zero. It therefore suffices to prove
    \[
        \lim_{r\downarrow0}\limsup_{j\to\infty}
        \int_{\hat{D}_r}\hat{N}\dif\point{q}=0.
    \]
    Recall that
    \[
        \hat{N}=\hat{\tau}_i\hat{h}^{ij}\hat{\tau}_j
        +\hat{X}^i\hat{h}_{ij}\hat{X}^j-2\hat{a}(\hat{\tau}).
    \]
    Fix $ 0<r<L/2 $. On $ \hat{D}_r $, the product rule gives
    \[
        \divr_{\Sigma}((r-\hat{\tau})^+\hat{\vec{X}})
        =(r-\hat{\tau})\divr_{\Sigma}\hat{\vec{X}}
        -\hat{X}^i\hat{h}_{ij}\hat{X}^j-\hat{a}(\hat{\tau}).
    \]
    Integrating, we obtain
    \[
        \begin{aligned}
            \int_{\hat{D}_r}\hat{X}^i\hat{h}_{ij}\hat{X}^j\dif\point{q}
            ={}&-\int_{\partial Q_j\cap\{\hat{\tau}<r\}}
            (r-\hat{\tau})\hat{\vec{X}}\cdot\vec{\nu}_{Q_j}\dif\mathcal{H}^{n-1}\\
            &+\int_{\hat{D}_r}(r-\hat{\tau})\divr_{\Sigma}\hat{\vec{X}}\dif\point{q}
            -\int_{\hat{D}_r}\hat{a}(\hat{\tau})\dif\point{q}.
        \end{aligned}
    \]
    The factor $ (r-\hat{\tau})^+ $ vanishes on the level set $ \{\hat{\tau}=r\} $. Since $ \hat{\vec{X}}\cdot\vec{\nu}_{Q_j}\geqslant0 $ and $ \hat{a}(\hat{\tau})\geqslant0 $, it follows that
    \[
        \int_{\hat{D}_r}\hat{X}^i\hat{h}_{ij}\hat{X}^j\dif\point{q}
        \leqslant\int_{\hat{D}_r}(r-\hat{\tau})\divr_{\Sigma}\hat{\vec{X}}\dif\point{q}.
    \]
    Denote the unscaled vector field at the $ j $-th base point by $ \vec{X} $. The divergence calculation in \cref{subsection::thesis::weak-monotonicity} and the nonpositive pure Hessian term give
    \[
        \divr_{\Sigma}\vec{X}\leqslant n+E^\tau_{\point{q}}\leqslant C.
    \]
    Under the blow-up,
    \[
        D\hat{\vec{X}}(\point{q})
        =\matr{A}_j^{-1}D\vec{X}(\matr{A}_j\point{q})\matr{A}_j,
    \]
    so taking the trace yields $ \divr_{\Sigma}\hat{\vec{X}}\leqslant C $, with $ C $ independent of $ j $. Consequently,
    \[
        \int_{\hat{D}_r}\hat{X}^i\hat{h}_{ij}\hat{X}^j\dif\point{q}
        \leqslant Cr|\hat{D}_r|\leqslant Cr,
    \]
    where we used $ \hat{D}_r\subset K_j $ and the uniform bound for $ |K_j| $, for all sufficiently large $ j $.

    The same estimate holds for $ \hat{\tau}_i\hat{h}^{ij}\hat{\tau}_j $ by interchanging source and target and using the transport change of variables, together with the uniform density bounds in \cref{definition::thesis::definition-of-blow-up}. These bounds also control the volume of $ \hat{\vec{T}}(\hat{D}_r) $ by a constant times $ |\hat{D}_r| $. Dropping the nonpositive term $ -2\hat{a}(\hat{\tau}) $, we obtain
    \[
        0\leqslant\limsup_{j\to\infty}\int_{\hat{D}_r}\hat{N}\dif\point{q}\leqslant Cr.
    \]
    Letting $ r\downarrow0 $ proves the required estimate. For any fixed $ R>0 $, split $ \hat{D}_R $ into $ \hat{D}_r $ and $ \{r\leqslant\hat{\tau}<R\} $. The integral over the latter set tends to zero by the fixed-strip convergence, applying it on $ \{r/2<\hat{\tau}<R\} $. First letting $ j\to\infty $ and then $ r\downarrow0 $ proves the lemma.
\end{proof}

\paragraph{Equation that gives homogeneity}

We first show 
\begin{lemma}[$ \tau $ is similar to $ \point{q} \cdot D u $ after blow-up]\label{lemma::thesis::integral-of-x-Du-tau-tends-to-zero}
    On $ K_\infty $, as $ j \to \infty $,
    \[
        D_{\point{q}} \hat{u} \cdot \point{q} - \hat{\tau} \longrightarrow 0
    \]
    converges in distribution and $ L^1_{\mathrm{loc}}(K_\infty) $.
\end{lemma}
\begin{proof}
    The idea is to use
    \[
        q^i \hat{u}_i - \hat{\tau} = (q^i - \hat{X}^i) \hat{u}_i + (\hat{X}^i \hat{u}_i - \hat{\tau}),
    \]
    where \cref{lemma::thesis::positiveness-of-sigma} gives $ X^i u_i \leqslant \tau + C \tau^2 $. After blow-up, both the value of $ \hat{X}^i \hat{u}_i $ and $ \hat{\tau} $ become $ s_j^{-1} $ of their corresponding original value, thus
    \[
        \hat{X}^i \hat{u}_i - \hat{\tau} \leqslant s_j^{-1} \left( \tau + C\tau^2 \right) - s_j^{-1} \tau = C s_j \hat{\tau}^2.
    \]
    This gives
    \[
        q^i \hat{u}_i - \hat{\tau} \leqslant (q^i - \hat{X}^i) \hat{u}_i + C s_j \hat{\tau}^2.
    \]
    Apply the divergence theorem to the first term, and the second term tends to zero whenever $ \hat{\tau} $ is uniformly bounded. We are considering the convergence over $ K_\infty $, i.e. the limit of $ K_j $, and we want to control the value of $ \hat{\tau} $, but $ K_j $ itself is not contained in $ \hat{D}_L $, so we take a cutoff: take a smooth nonincreasing $ f_M(\tau) $ such that $ f_M = 1 $ where $ \tau \leqslant M $, $ f_M = 0 $ where $ \tau > 2M $, and the derivative $ |f_M^{\prime}| \leqslant C / M $. Apply the divergence theorem to $ f_M(\hat{\tau}) (\point{q} - \hat{\vec{X}}) (L - \hat{u})^{+} $
    \[
    \begin{aligned}
        \int_{K_j} f_M(\hat{\tau}) (q^i - \hat{X}^i) \hat{u}_i \dif \point{q}
        ={}&    \int_{K_j}(L - \hat{u}) f_M(\hat{\tau}) (n - \divr_{\Sigma} \hat{\vec{X}}) \dif \point{q} \\
        &+ \int_{K_j}(L - \hat{u}) f_M^{\prime}(\hat{\tau})(q^i - \hat{X}^i)\partial_i[\hat{\tau}(\point{q},\hat{T}(\point{q}))] \dif \point{q} \\
        &-\int_{\partial Q_j \cap \{\hat{u} < L \}}(L - \hat{u}) f_M(\hat{\tau}) (q^i - \hat{X}^i) \nu_i \dif \mathcal{H}^{n-1}.
    \end{aligned}
    \]
    \begin{enumerate}
        \item For the first line we use proved bounds for $ E^{\tau}_{\point{q}} $. For convenience we use terms without a hat, like $ u, \vec{X} $ to refer to the $ j $-th unscaled term(which should be $ u_j $ as defined). Recall from the proof of \cref{lemma::thesis::int-Nj->0} that $ \divr_{\Sigma} \hat{\vec{X}}(\point{q}) = \divr_{\Sigma} \vec{X}(\matr{A}_j \point{q}) $, and from the proof of the monotonicity formula there is
        \[
            \divr_{\Sigma} \vec{X} = n + E^{\tau}_{\point{q}} + \left( \nabla^2 \tau \right)_{\overline{k} \overline{l}} b^{\overline{k} i} b^{\overline{l} j} h_{ij},
        \]
        so 
        \[
            n - \divr_{\Sigma} \vec{X} = -E^{\tau}_{\point{q}} - \left( \nabla^2 \tau \right)_{\overline{k} \overline{l}} b^{\overline{k} i} b^{\overline{l} j} h_{ij}.
        \]
        the second term is part of $ I_1 $, so by \cref{lemma::thesis::limit-chi-picture-and-zero-negative-terms} we have its integral tends to zero as $ j \to \infty $; For the term $ E^{\tau}_{\point{q}} $, the support of $ f_M(\hat{\tau}) $ lies in $ \hat{D}_{2M} $. Thus the corresponding original source points lie in $ D_{2Ms_j} $, and their transport images lie in $ T(D_{2Ms_j}) $. For fixed $ M $, the diameter estimates from \cref{corollary::thesis::diameter-of-Dr-by-r} imply that both original coordinates approach $ (\point{0}, \point{0}) $ uniformly. Since $ E^{\tau}_{\point{q}} $ is uniformly continuous and $ E^{\tau}_{\point{q}}(\point{0}, \point{0}) = 0 $, the scaling gives
        \[
            \sup_{\point{q} \in \hat{D}_{2M}} \left| \hat{E}^{\hat{\tau}}_{\point{q}}(\point{q}, \hat{T}(\point{q})) \right| = \sup_{\point{q} \in \hat{D}_{2M}} \left| E^{\tau}_{\point{q}}(\matr{A}_j \point{q}, T(\matr{A}_j \point{q})) \right| \longrightarrow 0.
        \]
        Since $ 0 \leqslant (L - \hat{u}) f_M(\hat{\tau}) \leqslant L $ and $ |K_j| $ is uniformly bounded, the integral containing $ \hat{E}^{\hat{\tau}}_{\point{q}} $ tends to zero as $ j \to \infty $.
        
        Combining both estimates we have the first line tends to zero.
        
        \item In the third line, the convexity of $ Q_j $ implies $ q^i \nu_i $ nonnegative(because $ q_i $ is pointing from $ \point{0} $ to a boundary point, which is pointing ``outward'' when considering the outer normal at $ \point{q} $), we have
        \[
            -\int_{\partial Q_j \cap \{\hat{u} < L \}}(L - \hat{u}) f_M(\hat{\tau}) (q^i - \hat{X}^i) \nu_i \dif \mathcal{H}^{n-1} \leqslant L \int_{\partial Q_j \cap \{\hat{u} < L \}} f_M(\hat{\tau}) \hat{X}^i \nu_i \dif \mathcal{H}^{n-1}
        \]
        which tends to zero as $ \mathcal{Q}_j(2M) \to 0 $ proved in \cref{lemma::thesis::limit-chi-picture-and-zero-negative-terms}.
        \item Finally for the second line, we split it into two parts, containing $ q^i \partial_i \hat{\tau} $ and $ \hat{X}^i \partial_i \hat{\tau} $ separately. Applying divergence theorem on $ (L - \hat{u})^+ (1 - f_M(\hat{\tau})) \point{q} $ for the first part gives
        \[
        \begin{aligned}
            0 \leqslant{}& \int_{\partial Q_j \cap \{\hat{u} < L\}}(L - \hat{u})(1 - f_M(\hat{\tau})) q^i \nu_i \dif \mathcal{H}^{n-1} \\
            ={}& n\int_{K_j}(L - \hat{u})(1 - f_M(\hat{\tau})) \dif \point{q} \\
            &-\int_{K_j}(1 - f_M(\hat{\tau})) q^i \hat{u}_i \dif \point{q} \\
            &-\int_{K_j}(L - \hat{u})f_M^{\prime}(\hat{\tau}) q^i \partial_i[\hat{\tau}(\point{q},\hat{T}(\point{q}))] \dif \point{q}.
        \end{aligned}
        \]
        The first line is bounded by $ nL |K_j \cap \{ \hat{\tau} > M \}| $; The second line is nonpositive so we just drop it; The third by inserting $ q^i \hat{u}_i \geqslant \hat{\tau} \geqslant 0 $ from \cref{lemma::thesis::positiveness-of-sigma}, we obtain
        \[
            \begin{aligned}
            &\int_{K_j}(L - \hat{u}) f_M^{\prime}(\hat{\tau}) q^i \partial_i[\hat{\tau}(\point  {q},\hat{T}(\point{q}))] \dif \point{q} \\
            &\quad \leqslant nL |K_j \cap \{ \hat{\tau} > M \}|.
            \end{aligned}
        \]
        For the volume of the intersection, we have
        \[
            \begin{aligned}
                M |K_j \cap \{ \hat{\tau} > M \}| \leqslant{}& \int_{K_j \cap \{ \hat{\tau} > M \}} \hat{\tau} \dif \point{q} \\
                \leqslant{}& \int_{K_j} q^i \hat{u}_i \dif \point{q} \\
            \end{aligned}
        \]
        Applying the divergence theorem to $ (L - \hat{u})^+ \point{q} $ and using $ q^i \nu_i \geqslant 0 $, we obtain $ \int_{K_j} q^i \hat{u}_i \dif \point{q} \leqslant nL |K_j| \leqslant C_L $. Thus $ |K_j \cap \{ \hat{\tau} > M \}| \leqslant C_L / M $. So this part of the second line is bounded by $ C / M $ with $ C $ independent of $ j $. By taking $ j \to \infty $ we have this part of the second line bounded by $ C / M $.

        For the part containing $ \hat{X}^i \partial_i \hat{\tau} $, expanding it gives $ \hat{X}^i \partial_i[\hat{\tau}(\point{q},\hat{T}(\point{q}))] = \hat{a}(\hat{\tau}) + \hat{X}^i \hat{h}_{ij} \hat{X}^j $, and by Cauchy-Schwarz
        \[
            \begin{aligned}
                \left|\hat{X}^i(\hat{\tau}_i-\hat{h}_{im}\hat{X}^m) \right|
                &\leqslant \left(\hat{h}_{im}\hat{X}^i\hat{X}^m\right)^{1/2}\left[(\hat{\tau}_i-\hat{h}_{im}\hat{X}^m)\hat{h}^{i l}(\hat{\tau}_l - \hat{h}_{l s} \hat{X}^s)\right]^{1/2} \\
                &= \left(\hat{h}_{im}\hat{X}^i\hat{X}^m\right)^{1/2}\hat{N}^{1/2} \leqslant \frac{1}{2} \hat{h}_{im}\hat{X}^i\hat{X}^m+\frac{1}{2} \hat{N}.
            \end{aligned}
        \]
        Therefore
        \[
            \hat{X}^i \hat{h}_{ij} \hat{X}^j = \hat{a}(\hat{\tau}) - \hat{X}^i (\hat{\tau}_i-\hat{h}_{im}\hat{X}^m) \leqslant \hat{a}(\hat{\tau}) + \frac{1}{2} \hat{h}_{im}\hat{X} ^i\hat{X}^m+\frac{1}{2} \hat{N},
        \] 
        inserting back, we have
        \[
            \hat{X}^i \partial_i[\hat{\tau}(\point{q},\hat{T}(\point{q}))] \leqslant 3 \hat{a}(\hat{\tau}) + \hat{N},
        \]
        Since $f_M'$ is supported in $\{M<\hat{\tau}<2M\}$ and $0\leqslant-f_M'\leqslant C/M$, we obtain
        \[
            \begin{aligned}
                &-\int_{K_j}(L-\hat{u})f_M'(\hat{\tau})
                \hat{X}^i\partial_i[\hat{\tau}(\point{q},\hat{T}(\point{q}))]
                \dif\point{q} \\
                &\qquad\leqslant
                \frac{CL}{M}
                \int_{K_j\cap\{M<\hat{\tau}<2M\}}
                \bigl(3\hat{a}(\hat{\tau})+\hat{N}\bigr)
                \dif\point{q}.
            \end{aligned}
        \]
        For each fixed $M$, the estimate $\hat{a}(t)=t+O(s_jt^2)$ implies $\hat{a}(\hat{\tau})\leqslant CM$ on this transition region for all sufficiently large $j$. Therefore,
        \[
            \begin{aligned}
                &-\int_{K_j}(L-\hat{u})f_M'(\hat{\tau})
                \hat{X}^i\partial_i[\hat{\tau}(\point{q},\hat{T}(\point{q}))]
                \dif\point{q} \\
                &\qquad\leqslant
                C_L|K_j\cap\{\hat{\tau}>M\}|
                +\frac{CL}{M}\int_{\hat{D}_{2M}}\hat{N}\dif\point{q}\\
                &\qquad\leqslant \frac{C_L}{M}+o_j(1),
            \end{aligned}
        \]
        where $o_j(1)\to0$ as $j\to\infty$ with $M$ fixed, by \cref{lemma::thesis::int-Nj->0}
    \end{enumerate}

    \noindent On the support of $f_M(\hat\tau)$ we have $\hat\tau\leqslant2M$. Then
    \[
        0\leqslant q^i\hat u_i-\hat\tau
        \leqslant(q^i-\hat X^i)\hat u_i+Cs_j\hat\tau^2,
    \]
    multiplying by $f_M(\hat\tau)$ and integrating, we obtain
    \[
        0 \leqslant
        \int_{K_j}f_M(\hat\tau)(q^i\hat u_i-\hat\tau)
            \dif\point q \leqslant
        \int_{K_j}f_M(\hat\tau)(q^i-\hat X^i)\hat u_i
            \dif\point q
        +4Cs_jM^2|K_j|.
    \]
    Taking $\limsup$ as $j\to\infty$ with $M$ fixed gives
    \[
        0\leqslant
        \limsup_{j\to\infty}
        \int_{K_j}f_M(\hat\tau)(q^i\hat u_i-\hat\tau)
            \dif\point q
        \leqslant\frac{C_L}{M}.
    \]

    Consider any test function $ \phi $ with $ G := \supp \phi \Subset K_\infty $, from \cref{lemma::thesis::bound-for-Du-on-arbitrary-G} we have $ |D \hat{u}| \leqslant C_G $ over $ G $, so $ \hat{\tau} \leqslant \point{q} \cdot D \hat{u} \leqslant C^{\prime} $ independent of $ j $. Whenever $ M > C^{\prime} $ we have $ f_M = 1 $ over $ G $. Now we have
    \[
        |\int_{K_\infty} \phi \left( D_{\point{q}} \hat{u} \cdot \point{q} - \hat{\tau} \right) \dif \point{q}| \leqslant \norm{\phi}_{L^{\infty}} \int_{G} \left( D_{\point{q}} \hat{u} \cdot \point{q} - \hat{\tau} \right) \dif \point{q} \leqslant \norm{\phi}_{L^{\infty}} \int_{K_j} f_M(\hat{\tau}) \left( D_{\point{q}} \hat{u} \cdot \point{q} - \hat{\tau} \right) \dif \point{q}.
    \]
    Now apply the convergence bound to the RHS, its $ \limsup $ is bounded by $ \norm{\phi}_{L^{\infty}} C / M $, now we take $ M \to \infty $ and obtain
    \[
        \int_{K_\infty} \phi \left( D_{\point{q}} \hat{u} \cdot \point{q} - \hat{\tau} \right) \dif \point{q} \longrightarrow 0.
    \]

    That is
    \[
        D_{\point{q}} u^{(j)} \cdot \point{q} - \tau^{(j)} \longrightarrow 0
    \]
    in distribution on  $ K_\infty $. Since $ D_{\point{q}} \hat{u} \cdot \point{q} - \hat{\tau} \geqslant 0 $, its distributional convergence to zero also implies convergence to zero in $ L^1_{\mathrm{loc}}(K_\infty) $.
\end{proof}

\begin{lemma}[Equation for homogeneity]\label{lemma::thesis::distributional-zero-2-homogeneity}
    \[
        D^2_{\point{q}\point{q}} u_\infty \cdot \point{q} - D_{\point{q}} u_\infty = 0
    \]
    in distribution.
\end{lemma}
\begin{remark}\label{remark::thesis::idea-of-2-homogeneity}
    In \cite{collins2025boundary} this equality is obtained directly by $ \hat{N} \to 0 $. But in our case we could only have something like
    \[
        \tau_i = \tau_{\overline{k}} b^{\overline{k} m} h_{mi}.
    \]
    In the case $ b = \Inner{\point{x}}{\overline{\point{x}}} $, we have $ Y^{\overline{k}} = \tau_i b^{i \overline{k}} = \overline{x}^{\overline{k}} $, $ X^{i} = \tau_{\overline{k}} b^{\overline{k} i} = x^i $ and $ h_{im} = u_{im} $. This motivates estimating the contracted differences $ \hat{h}_{im}(q^m-\hat{X}^m) $ and $ (\hat{u}_{im}-\hat{h}_{im})q^m $, which tend to zero in $ L^1_{\mathrm{loc}} $ in the proof below.
\end{remark}
\begin{proof} 
    By \cref{lemma::thesis::integral-of-x-Du-tau-tends-to-zero}, the derivative $ \partial_i \left[ D_{\point{q}} \hat{u} \cdot \point{q} - \hat{\tau} \right] $ tends to zero in distribution. 

    Consider 
    \[
        q^i \partial_i \left[ D_{\point{q}} \hat{u} \cdot \point{q} - \hat{\tau} \right],
    \]
    which still tends to zero in distribution, expanding it gives
    \[
        \begin{aligned}
            q^i \partial_i \left[ D_{\point{q}} \hat{u} \cdot \point{q} - \hat{\tau} \right] = q^i \hat{u}_{im} q^{m} + q^i \hat{u}_i - q^i \hat{\tau}_i - q^i \hat{h}_{im} \hat{X}^m
        \end{aligned}
    \]
    Following the idea from \cref{remark::thesis::idea-of-2-homogeneity}, we group terms to estimate $ q^i - \hat{X}^i $. Recall that $ \hat{\tau}_i \hat{X}^i = \hat{\tau}_i \hat{b}^{i \overline{k}} \hat{\tau}_{\overline{k}} = \hat{a}(\hat{\tau}) $, 
    \begin{equation}\label{equation::thesis::radial-derivative-of-qDu-tau}
        \begin{aligned}
            {}& q^i \partial_i \left[ D_{\point{q}} \hat{u} \cdot \point{q} - \hat{\tau} \right] \\
            ={}& (q^i - \hat{X}^i) \hat{h}_{im} (q^m - \hat{X}^m) + q^i \hat{b}_{im} q^m + (q^i - \hat{X}^i) \hat{h}_{im} \hat{X}^m + q^i \hat{u}_i - q^i \hat{\tau}_i \\
            ={}& |\point{q} - \hat{\vec{X}}|^2_{\hat{h}} + q^i \hat{b}_{im} q^m + (q^i - \hat{X}^i)(\hat{h}_{im} \hat{X}^m - \hat{\tau}_i) - \hat{X}^i \hat{\tau}_i + q^i \hat{u}_i \\
            ={}& |\point{q} - \hat{\vec{X}}|^2_{\hat{h}} + q^i \hat{b}_{im} q^m + (q^i - \hat{X}^i)(\hat{\tau}_{\overline{k}} \hat{b}^{\overline{k} m} \hat{h}_{mi} - \hat{\tau}_i) + (\hat{\tau} - \hat{a}(\hat{\tau})) + (q^i \hat{u}_i - \hat{\tau}),
        \end{aligned}
    \end{equation}
    here we denote $ |\point{q} - \hat{\vec{X}}|^2_{\hat{h}} := (q^i - \hat{X}^i) \hat{h}_{im} (q^m - \hat{X}^m) $.
    
    For the last two terms in \cref{equation::thesis::radial-derivative-of-qDu-tau}, $ \hat{\tau} - \hat{a}(\hat{\tau}) $ can be controlled by its Taylor expansion; while the third inequality from \cref{lemma::thesis::positiveness-of-sigma} gives $ q^i \hat{u}_i - \hat{\tau} \geqslant 0 $. 

    We then consider the first two terms in \cref{equation::thesis::radial-derivative-of-qDu-tau}. Although $ \hat{b}_{im} $ itself may not be semi-positive definite, as proved in \cref{lemma::thesis::positiveness-of-D2b-DbDb}, by adding $ k s_j \hat{u}_i \hat{u}_m $ it is semi-positive definite. The Cauchy--Schwarz inequality then gives
    \[
        \begin{aligned}
            {}& |\point{q} - \hat{\vec{X}}|^2_{\hat{h}} + q^i \left[ \hat{b}_{im} + k s_j \hat{u}_i \hat{u}_m \right] q^m \\
            \leqslant{}& q^i \partial_i \left[ D_{\point{q}} \hat{u} \cdot \point{q} - \hat{\tau} \right] - (q^i - \hat{X}^i)(\hat{\tau}_{\overline{k}} \hat{b}^{\overline{k} m} \hat{h}_{mi} - \hat{\tau}_i) + k s_j (q^i \hat{u}_i)^2 + |\hat{\tau} - \hat{a}(\hat{\tau})| \\
            \leqslant{}& q^i \partial_i \left[ D_{\point{q}} \hat{u} \cdot \point{q} - \hat{\tau} \right] + \frac{1}{2} |\point{q} - \hat{\vec{X}}|^2_{\hat{h}} + \frac{1}{2} \hat{N} + k s_j (q^i \hat{u}_i)^2 + |\hat{\tau} - \hat{a}(\hat{\tau})|.
        \end{aligned}
    \]
    Rearranging gives
    \begin{equation}\label{equation::thesis::label-for-four-terms}
        \frac{1}{2}|\point{q} - \hat{\vec{X}}|^2_{\hat{h}} + q^i \left[ \hat{b}_{im} + k s_j \hat{u}_i \hat{u}_m \right] q^m \leqslant q^i \partial_i \left[ D_{\point{q}} \hat{u} \cdot \point{q} - \hat{\tau} \right] + \frac{1}{2} \hat{N} + ks_j (q^i \hat{u}_i)^2 + |\hat{\tau} - \hat{a}(\hat{\tau})|. 
    \end{equation}

    Take a smooth test function $ 0 \leqslant \phi \in C^{\infty}_c (K_{\infty}) $, and $ G := \supp \phi \Subset K_{\infty} $. Now we consider the four terms in the RHS of \cref{equation::thesis::label-for-four-terms}, after multiplying by $ \phi $ and integrating over $ K_\infty $
    \begin{itemize}
        \item The first tends to zero since it tends to zero in the sense of distribution; 
        \item The second tends to zero by \cref{lemma::thesis::int-Nj->0}, using \cref{lemma::thesis::bound-for-Du-on-arbitrary-G} to contain $ G $ in a fixed $ \hat{D}_{r_G} $;
        \item The third by using the bound for $ D \hat{u} $ from \cref{lemma::thesis::bound-for-Du-on-arbitrary-G}, we have $ (q^i \hat{u}_i)^2 $ bounded over $ G $, so since $ s_j \to 0 $ its integral tends to zero;
        \item The fourth by Taylor expansion we have $ |\hat{\tau} - \hat{a}(\hat{\tau})| = s_j^{-1} |s_j\hat{\tau} - a(s_j\hat{\tau})| \leqslant C s_j \hat{\tau}^2 \leqslant C s_j r_G^2 $, so its integral also tends to zero.
    \end{itemize}
    
    Since both terms on the left-hand side of \cref{equation::thesis::label-for-four-terms} are nonnegative and the integral of the right-hand side multiplied by $ \phi $ tends to zero, we have as $ j \to \infty $
    \[
        (q^i-\hat{X}^i) \hat{h}_{im} (q^m - \hat{X}^m) \longrightarrow 0,
    \]
    as well as
    \[
        q^i\left[ \hat{b}_{im} + ks_j \hat{u}_i \hat{u}_m \right]q^m \longrightarrow 0,
    \]
    in the sense $ L^1_{\mathrm{loc}}(K_\infty) $(which is equivalent to in distribution as the two terms are nonnegative).

    To estimate the remaining contracted differences, we use the following matrix inequality: Let $ \matr{H} $ be a semi-positive definite matrix, $ \vec{v} $ be a vector, then
    \[
        |\matr{H} \vec{v}|^2 \leqslant \tr \matr{H} \vec{v}^{\top}\matr{H} \vec{v}.
    \]
    Applying it to matrix $ \hat{h}_{i m} $ and vector $ q^i - \hat{X}^i $ gives
    \[
        \int_{G} |\hat{h}_{im}(q^m - \hat{X}^m)| \dif \point{q} \leqslant \left( \int_{G} \tr \hat{h} \right)^{\frac{1}{2}} \left( \int_{G} (q^i - \hat{X}^i) \hat{h}_{im} (q^m - \hat{X}^m) \right)^{\frac{1}{2}}.
    \]
    Applying it to matrix $ \hat{b}_{im} + ks_j \hat{u}_i \hat{u}_m $ and $ q^i $ gives
    \[
        \int_{G} |(\hat{b}_{im} + ks_j \hat{u}_i \hat{u}_m)q^i| \dif \point{q} \leqslant \left( \int_{G} \tr(\hat{b}_{im} + ks_j \hat{u}_i \hat{u}_m) \dif \point{q} \right)^{\frac{1}{2}} \left( \int_{G} q^i\left[ \hat{b}_{im} + k s_j \hat{u}_i \hat{u}_m \right] q^m \right)^{\frac{1}{2}}
    \]
    If the integrals over $ G $ of $ \tr\hat{h} $ and $ \tr(\hat{b}_{im}+ks_j\hat{u}_i\hat{u}_m) $ are bounded uniformly in $ j $, then the two left-hand sides tend to zero as $ j\to\infty $. 
    
    Consider the convexification $ \hat{w} $ of $ \hat{u} $, we have $ \hat{h}_{im} + \hat{b}_{im} + ks_j \hat{u}_i \hat{u}_m = e^{-ks_j \hat{u}} \hat{w}_{im} $, so we can bound the two traces by $ \tr D^2 \hat{w} = \laplacian \hat{w} $. Take a new test smooth function $ 0 \leqslant \psi \in C_c^\infty(K_\infty) $ with $ \supp \psi \Subset K_\infty $ and $ \psi = 1 $ over $ G $, using convexity of $ \hat{w} $ and uniform boundedness from \cref{lemma::thesis::bound-for-Du-on-arbitrary-G}, 
    \[
        \int_{G} \laplacian \hat{w} \dif \point{q} \leqslant \int_{K_\infty} \psi \laplacian \hat{w} \dif \point{q} = \int_{K_\infty} \hat{w} \laplacian \psi \dif \point{q} \leqslant C_G. 
    \]
    Since both matrices are semi-positive definite, the integrals of their traces over $ G $ are bounded uniformly in $ j $. Thus
    \[
        |\hat{h}_{i m}(q^m-(\hat{X})^m)| \longrightarrow 0, \quad |(\hat{b}_{im} + k s_j \hat{u}_i \hat{u}_m) q^m| \longrightarrow 0.
    \]
    in distribution. With $ k s_j \hat{u}_i \hat{u}_m \longrightarrow 0 $, so additionally
    \[
        |\hat{b}_{im} q^m| \longrightarrow 0.
    \]
    
    Apply the same matrix inequality to $ \hat{N} = |\hat{\tau}_{\overline{k}} \hat{b}^{\overline{k} m} \hat{h}_{mi} - \hat{\tau}_i|^2_{\hat{h}} $, taking matrix $ \hat{h}_{im} $ and vector $ \hat{\tau}_m \hat{h}^{m i} - \hat{\tau}_{\overline{k}} \hat{b}^{\overline{k} i} $. By the same argument as above
    \[
        \int_{G} |\hat{\tau}_i - \hat{\tau}_{\overline{k}} \hat{b}^{\overline{k} m} \hat{h}_{mi}| \dif \point{q} \leqslant \left( \int_{G} \tr \hat{h} \dif \point{q} \right)^{\frac{1}{2}} \left( \int_{G} \hat{N} \dif \point{q} \right)^{\frac{1}{2}} \longrightarrow 0,
    \]
    Therefore we have
    \[
        |\hat{\tau}_i - \hat{\tau}_{\overline{k}} \hat{b}^{\overline{k} m} \hat{h}_{mi}| \longrightarrow 0
    \]
    in distribution.

    Now we rearrange terms in $ \partial_i \left[  D_{\point{q}} \hat{u} \cdot \point{q} - \hat{\tau} \right] $ to isolate an expression of the form $ \hat{u}_{mi} q^m - \hat{u}_i $:
     \[
        \hat{u}_{mi} q^m - \hat{u}_i = -\partial_i \left[ D_{\point{q}} \hat{u} \cdot \point{q} - \hat{\tau} \right] + 2 \hat{b}_{mi} q^m - (\hat{\tau}_i - \hat{\tau}_{\overline{k}} \hat{b}^{\overline{k} m} \hat{h}_{mi}) + 2 \hat{h}_{mi}(q^m - \hat{X}^m).
    \]
    Each term on the RHS tends to zero in distribution on $ K_\infty $. Since $ \hat{u} \to u_\infty $ locally uniformly, passing to the limit gives
    \[
        D^2_{\point{q}\point{q}} u_\infty \cdot \point{q} - D_{\point{q}} u_\infty = 0
    \]
    in distribution on $ K_\infty $.
\end{proof}

\paragraph{The homogeneity}

Now we prove
\begin{lemma}[Homogeneity of the limit]\label{lemma::thesis::2-homogeneity-of-limit}
    The limit $ u_\infty $ is homogeneous of degree two: for every $ \point{q} \in K_\infty $ and $ 0 < t \leqslant 1 $,
    \[
        u_\infty(t\point{q}) = t^2 u_\infty(\point{q}).
    \]
    Moreover, as $ j \to \infty $,
    \[
        \hat{\tau} \longrightarrow 2 u_\infty
    \]
    in distribution on $ K_\infty $ and also in $ L^1_{\mathrm{loc}}(K_\infty) $.
\end{lemma}
\begin{proof}
By \cref{lemma::thesis::distributional-zero-2-homogeneity},
\[
    D^2_{\point{q}\point{q}} u_\infty \cdot \point{q} - D_{\point{q}} u_\infty = 0
\]
in distribution. Rewriting it gives
\[
    D_{\point{q}} \left( D_{\point{q}} u_\infty \cdot \point{q} - 2 u_\infty \right) = D^2_{\point{q}\point{q}} u_\infty \cdot \point{q} - D_{\point{q}} u_\infty,
\]
Since $ K_\infty $ is connected, there is a constant $ a_\ast $ such that
\[
    D_{\point{q}} u_\infty \cdot \point{q} - 2 u_\infty = 2 a_\ast
\]
a.e. in $ K_\infty $, i.e. $ \point{q} \cdot D_{\point{q}} u_\infty = 2 (u_\infty + a_\ast) $. Using Fubini's theorem in polar coordinates, over almost every closed radial segment $ t\point{q} $ inside $ K_\infty $, $ 0 < a \leqslant t \leqslant b $, there is 
\[
    t \frac{\dif}{\dif t} (u_{\infty}(t\point{q}) + a_\ast) = 2 (u_\infty(t\point{q}) + a_\ast),
\]
Since $ u_\infty $ is convex and thus absolutely continuous on each compact radial segment in the domain, we can apply the fundamental theorem of calculus to the formula above, giving
\[
    u_{\infty}(t\point{q}) + a_\ast = t^2 (u_\infty(\point{q}) + a_\ast)
\]
on almost every closed radial segment in $ K_{\infty} $. 

By the convexity of $ \hat{w} $, we have
\[
    0 \leqslant \hat{w}(t\point{q}) \leqslant (1 - t) \hat{w}(0) + t \hat{w}(\point{q}) = t \hat{w}(\point{q}),
\]
Taking the limit as $ j \to \infty $ gives
\[
    0 \leqslant u_\infty(t\point{q}) \leqslant t u_\infty(\point{q}).
\]
Therefore
\[
    a_\ast = \lim_{t \to 0} (u_\infty(t\point{q}) + a_\ast) = \lim_{t \to 0} t^2 (u_\infty(\point{q}) + a_\ast) = 0.
\]
so $ u_\infty(t \point{q}) = t^2 u_{\infty}(\point{q}) $ a.e. in $ K_\infty $. And then the 2-homogeneity holds over the whole $ K_\infty $ since $ u_\infty $ is continuous.

Now recall from \cref{lemma::thesis::integral-of-x-Du-tau-tends-to-zero}
\[
    D_{\point{q}} \hat{u} \cdot \point{q} - \hat{\tau} \longrightarrow 0
\]
in distribution. Since $ \hat{u} \to u_\infty $ locally uniformly, $ D_{\point{q}} \hat{u} \cdot \point{q} \to D_{\point{q}} u_\infty \cdot \point{q} = 2 u_\infty $ in distribution. Thus $ \hat{\tau} \to 2 u_\infty $ in distribution. Combining this with $ D_{\point{q}} \hat{u} \to D_{\point{q}} u_\infty $ in $ L^1_{\mathrm{loc}}(K_\infty) $, we also obtain $ \hat{\tau} \to 2 u_\infty $ in $ L^1_{\mathrm{loc}}(K_\infty) $.
\end{proof}

And we prove 
\begin{lemma}[Convergence of sections]\label{lemma::thesis::convergence-of-sections-at-any-t}
    Under the assumptions of \cref{lemma::thesis::limit-chi-picture-and-zero-negative-terms}, for every $ 0 < t < L / 2 $,
    \[
        \cl S_t[\hat{u}] \longrightarrow \cl S_t[u_\infty]
    \]
    in Hausdorff sense. Moreover,
    \[
        S_t[u_\infty] = \sqrt{t} S_1[u_\infty].
    \]
\end{lemma}
\begin{proof}
    We first prove convergence of the volumes of the radius sublevels. The estimate $ u\leqslant\tau+C\tau^2 $ from \cref{lemma::thesis::positiveness-of-sigma} becomes
    \[
        \hat{u}\leqslant\hat{\tau}+Cs_j\hat{\tau}^2.
    \]
    Thus, for every fixed $ 0<r<L $, we have $ \hat{D}_r\subset K_j $ for all sufficiently large $ j $, since $ r+Cs_jr^2<L $.

    By \cref{lemma::thesis::2-homogeneity-of-limit}, $ \hat{\tau}\to2u_\infty $ in $ L^1_{\mathrm{loc}}(K_\infty) $. For every $ r>0 $, the level set $ \{2u_\infty=r\} $ has zero volume: two-homogeneity implies that each ray from the origin meets this positive level at most once, so the assertion follows from polar Fubini. Consequently, on every compact set $ G\subset K_\infty $,
    \[
        \mathbf{1}_{\hat{D}_r}\longrightarrow
        \mathbf{1}_{S_{r/2}[u_\infty]}
        \quad\text{in }L^1(G).
    \]
    Hausdorff convergence of the full-dimensional convex bodies $ \cl K_j $ to $ \cl K_\infty $ gives $ |K_j|\to|K_\infty| $. Given $ \epsilon>0 $, choose a compact set $ G\subset K_\infty $ with $ |K_\infty\setminus G|<\epsilon $. Then $ G\subset K_j $ for all sufficiently large $ j $, and $ |K_j\setminus G|<2\epsilon $. Since $ \hat{D}_r\subset K_j $, the parts outside $ G $ contribute at most $ 3\epsilon $ to the difference of the two volumes. Letting $ j\to\infty $ and then $ \epsilon\downarrow0 $, we obtain
    \[
        |\hat{D}_r|\longrightarrow|S_{r/2}[u_\infty]|,
        \qquad 0<r<L.
    \]

    Recall the unscaled source densities $ \rho_j^0 $ from \cref{definition::thesis::definition-of-blow-up}. Their uniform positive lower and upper bounds allow us to pass to a further subsequence such that $ \rho_j^0(\point{0})\to\rho_\ast>0 $. The diameter estimate in \cref{corollary::thesis::diameter-of-Dr-by-r} shows that the original points corresponding to $ \hat{D}_r $ approach the base point uniformly. The uniform continuity of the densities in the cost coordinates therefore gives
    \[
        \sup_{\hat{D}_r}|\hat{\rho}-\rho_\ast|\longrightarrow0.
    \]
    By \cref{lemma::thesis::limit-chi-picture-and-zero-negative-terms}, $ \hat{\chi}(r)\to M $ for some $ M>0 $. Hence
    \[
        M r^{n/2}
        =\lim_{j\to\infty}\int_{\hat{D}_r}\hat{\rho}\dif\point{q}
        =\rho_\ast|S_{r/2}[u_\infty]|,
        \qquad 0<r<L.
    \]
    Equivalently,
    \[
        |S_t[u_\infty]|=\frac{M}{\rho_\ast}(2t)^{n/2},
        \qquad 0<t<L/2.
    \]

    For $ 0<a<b<L/2 $, convexity of $ K_\infty $, with $ \point{0}\in\cl K_\infty $, and two-homogeneity give
    \[
        \sqrt{\frac ab}\,S_b[u_\infty]\subset S_a[u_\infty].
    \]
    The volume formula shows that these open convex sets have equal volume, so they are equal. Since $ L>2 $, comparison with height $ 1 $ gives
    \[
        S_t[u_\infty]=\sqrt{t}\,S_1[u_\infty],
        \qquad 0<t<L/2.
    \]

    We now show Hausdorff convergence of $ \cl S_t[\hat{u}] $ for $ 0<t<L/2 $. For an interior point $ \point{q} $ of $ S_t[u_\infty] $, we have $ u_\infty(\point{q}) < t $, so $ \hat{u}(\point{q}) < t $ for all $ j $ sufficiently large. Taking a finite $ \epsilon $-ball cover of $ \cl S_t[u_\infty] $ with centers in $ S_t[u_\infty] $ and taking an $ N_\epsilon $ such that all centers of the chosen balls satisfy $ \hat{u}(\point{q}) < t $ shows that when $ j > N_{\epsilon} $, there is $ \dist (\point{q}, S_t[\hat{u}]) \leqslant \epsilon $ for all $ \point{q} \in S_t[u_\infty] $.

    For the reverse inclusion, take a sequence of $ \point{q}_j \in \cl S_t[\hat{u}] $ converging to $ \point{q} $. At least we have $ \point{q} \in \cl K_\infty $. Choose $ \point{y} \in K_\infty $ with $ u_\infty(\point{y}) < t $. For a fixed $ 0 < \theta < 1 $, let
    \[
        \point{q}_{\theta} = (1 - \theta)\point{q} + \theta \point{y} \in K_\infty.
    \]
    By convexity of $ \hat{w} $, 
    \[
        \hat{w}((1 - \theta)\point{q}_j + \theta \point{y}) \leqslant (1 - \theta)\hat{w}(\point{q}_j) + \theta \hat{w}(\point{y}) \leqslant (1 - \theta) \frac{e^{k s_j t} - 1}{k s_j} + \theta \hat{w}(\point{y}).
    \]
    Taking limits
    \[
        u_\infty(\point{q}_{\theta}) \leqslant (1 - \theta)t + \theta u_{\infty}(\point{y}) < t.
    \]
    Therefore $ \point{q}_\theta \in \cl S_t[u_\infty] $, and taking $ \theta \to 0 $ we have $ \point{q} \in \cl S_t[u_\infty] $. So the Hausdorff convergence 
    \[
        \cl S_t[\hat{u}] \longrightarrow \cl S_t[u_\infty]
    \]
    holds.
\end{proof}

\subsection{Application of blow-up homogeneity}

\begin{lemma}\label{lemma::thesis::inclusion-with-moving-points}
    Take a sequence of blow-ups as defined in \cref{definition::thesis::definition-of-blow-up} and \ref{proposition::thesis::remaining-terms}. Suppose that
    \[
        \frac{\hat{\chi}(a)}{\hat{\chi}(b)} \longrightarrow 1
    \]
    for any fixed $ 0 < a < b $. Set $ h_j := L s_j $, and denote by $ \mathcal{E}_{j}(h) $ the John ellipsoid of the original section $ S_h[u_j] $ with its center translated to the origin.

    Then for every fixed $ 0 < \delta < 1 $, when $ j $ is large enough the ellipsoids satisfy
    \[
        (1 - \delta) \frac{1}{\sqrt{2}} \mathcal{E}_{j}(h_j) \subset \mathcal{E}_{j}(h_j / 2) \subset (1 + \delta) \frac{1}{\sqrt{2}} \mathcal{E}_{j}(h_j).
    \]
\end{lemma}
\begin{proof}
    The idea is the same as in \cite[Proposition 5.1]{collins2025boundary}. Suppose that, for some fixed $ 0<\delta<1 $, one of the inclusions fails along a subsequence. For this proof, apply the blow-up construction in \cref{definition::thesis::definition-of-blow-up} using the John ellipsoid of the larger section $ S_{4Ls_j}[u_j] $, while retaining the height scale $ s_j $. Changing the normalization matrix multiplies $ \hat{\chi}(r) $ by a factor independent of $ r $, so the ratio assumption is preserved.

    Passing to a further subsequence, the preceding lemmata give a two-homogeneous limit $ u_\infty $. With $ 4L $ in place of $ L $, \cref{lemma::thesis::convergence-of-sections-at-any-t} applies at every height below $ 2L $, in particular at $ L $ and $ L/2 $. Thus the normalized sections at these two heights converge in Hausdorff distance, and their limits satisfy
    \[
        S_{L/2}[u_\infty]=\frac{1}{\sqrt{2}}S_L[u_\infty].
    \]
    By continuity of John ellipsoids with respect to Hausdorff convergence, see \cite[Proposition 3.8]{artstein2022some}, their origin-centered translations converge to ellipsoids satisfying the same dilation relation. Consequently, for the fixed $ \delta $, both claimed inclusions hold for all sufficiently large $ j $. Applying the inverse linear normalization gives the inclusions for $ \mathcal{E}_j(h_j) $ and $ \mathcal{E}_j(h_j/2) $, contradicting the choice of the subsequence.
\end{proof}

\begin{corollary}\label{corollary::thesis::fixed-point-half-height-inclusion}
    Fix a base point and take a blow-up. For every $ 0 < \epsilon < \frac{1}{2} $, there exists $ h_{\epsilon} > 0 $ depending on the base point, such that for every $ 0 < h < h_{\epsilon} $,
    \[
        2^{-\frac{1}{2} - \epsilon} \mathcal{E}(h) \subset \mathcal{E}(h / 2) \subset 2^{-\frac{1}{2} + \epsilon} \mathcal{E}(h),
    \]
    where $ \mathcal{E}(h) $ still denotes the origin-centered translation of the John ellipsoid of $ S_h[u] $.
\end{corollary}
\begin{proof}
    At the fixed base point, the positive finite limit of $ \chi $ at zero implies that the ratio hypothesis of \cref{lemma::thesis::inclusion-with-moving-points} holds along every blow-up sequence. 
    
    Choose $ 0 < \delta < 1 $ such that $ 1 - \delta \geqslant 2^{-\epsilon} $ and $ 1 + \delta \leqslant 2^{\epsilon} $. If the inclusions failed at arbitrarily small heights, we could choose a sequence of such heights tending to zero. Applying \cref{lemma::thesis::inclusion-with-moving-points} to this sequence gives both inclusions for all sufficiently large indices, leading to a contradiction.
\end{proof}

Now we pick a universal $ h_{\epsilon} $ for all base points.
\begin{lemma}\label{lemma::thesis::uniform-aepsilon-scale-inclusion}
    For every $ 0 < \epsilon < \frac{1}{2} $, there exist $ 0 < a_{\epsilon} < 1 $ and $ h_{\epsilon} > 0 $, depending only on $ \epsilon $ and the problem $ (\regs, \regt, \rho, \overline{\rho}, b) $, such that for every base point $ \point{q} \in \cl Q $, after moving its transport pair to $ (\point{0}, \point{0}) $, and every $ 0 < h < h_{\epsilon} $,
    \[
        a_{\epsilon}^{\frac{1}{2} + \epsilon} \mathcal{E}_{\point{q}}(h) \subset \mathcal{E}_{\point{q}}(a_{\epsilon} h) \subset a_{\epsilon}^{\frac{1}{2} - \epsilon} \mathcal{E}_{\point{q}}(h),
    \]
    where $ \mathcal{E}_{\point{q}}(h) $ denotes the John ellipsoid of the height-$ h $ section at $ \point{q} $, with its center translated to the origin, in the corresponding cost coordinates.
\end{lemma}
\begin{proof}
    For a fixed $ \epsilon $, we call a height $ h > 0 $ \textbf{good} at $ \point{q} $ if
    \[
        2^{-\frac{1}{2} - \frac{\epsilon}{2}} \mathcal{E}_{\point{q}}(h) \subset \mathcal{E}_{\point{q}}(h / 2) \subset 2^{-\frac{1}{2} + \frac{\epsilon}{2}} \mathcal{E}_{\point{q}}(h),
    \]
    and \textbf{bad} otherwise. The idea is to bound the number of bad dyadic heights at each base point uniformly.

    Define
    \[
        \mathcal{M}_{\point{q}}(r) := \exp\left(-\int_{0}^{r} \frac{\beta_{\point{q}}(t)}{t} \dif t \right) \chi_{\point{q}}(r),
    \]
    which is nonincreasing in $ r $ by the monotonicity formula. By \cref{lemma::thesis::estimate-of-sections-and-limit-of-chi}, the inclusion $ D_r \subset S_{2r} $ and the uniform upper bound for the density give $ \chi_{\point{q}}(r) \leqslant C_1 $, uniformly in the base point. Since the exponential correction is at most one, we also have $ \mathcal{M}_{\point{q}}(r) \leqslant C_1 $.

    Fix a small common $ r_\ast > 0 $ within the range of the monotonicity formula and the section-radius inclusions. We prove a positive lower bound for $ \mathcal{M}_{\point{q}}(r_\ast) $ independent of the base point. Denote the normalized potential and coordinate domain at $ \point{q} $ by $ u_{\point{q}} $ and $ Q_{\point{q}} $. Let $ w_{\point{q}} = \Phi(u_{\point{q}}) $ be the common convexification, where $ \Phi(t) := (e^{kt} - 1) / k $ for $ k > 0 $ and $ \Phi(t) := t $ for $ k = 0 $. Since the original potential and cost are bounded, subtracting the supporting mountain gives a common bound for $ u_{\point{q}} $, and hence
    \[
        0 \leqslant w_{\point{q}} \leqslant M, \qquad w_{\point{q}}(\point{0}) = 0,
    \]
    with $ M > 0 $ independent of the base point. Choose $ 0 < \lambda < 1 $ such that $ \lambda M < \Phi(r_\ast) $. For every $ \point{x} \in Q_{\point{q}} $, convexity of the domain and the convexification gives
    \[
        \lambda \point{x} \in Q_{\point{q}}, \qquad w_{\point{q}}(\lambda \point{x}) \leqslant \lambda w_{\point{q}}(\point{x}) \leqslant \lambda M < \Phi(r_\ast).
    \]
    Thus $ \lambda Q_{\point{q}} \subset S_{r_\ast}[u_{\point{q}}] $. Using \cref{equation::thesis::section-radius-inclusions}, $ \frac{1}{4} S_{r_\ast}[u_{\point{q}}] \subset D_{\point{q},r_\ast} $, we obtain
    \[
        \frac{\lambda}{4} Q_{\point{q}} \subset D_{\point{q},r_\ast}.
    \]
    Nondegeneracy of the cost on the compact product of the domain closures gives uniform positive lower and upper bounds for the Jacobian determinants of the cost-coordinate maps. Consequently, there are constants $ c_Q, \rho_\ast > 0 $ such that $ |Q_{\point{q}}| \geqslant c_Q $ and the transformed density satisfies $ \rho_{\point{q}} \geqslant \rho_\ast $, uniformly in the base point. Therefore
    \[
        \chi_{\point{q}}(r_\ast) = r_\ast^{-\frac{n}{2}} \int_{D_{\point{q},r_\ast}} \rho_{\point{q}} \dif \point{x} \geqslant r_\ast^{-\frac{n}{2}} \rho_\ast \left( \frac{\lambda}{4} \right)^n c_Q =: c_\ast > 0.
    \]
    {Set $ \vartheta:=\alpha/(4n) $. The uniform error estimate $ 0 \leqslant \beta_{\point{q}}(r) \leqslant C r^\vartheta $ gives}
    \[
        \mathcal{M}_{\point{q}}(r_\ast) \geqslant {\exp\left(-\frac{C r_\ast^\vartheta}{\vartheta}\right)} c_\ast =: c_1 > 0.
    \]
    By monotonicity, for every base point and every $ 0 < r \leqslant r_\ast $,
    \[
        c_1 \leqslant \mathcal{M}_{\point{q}}(r) \leqslant C_1.
    \]

    Now we claim a criterion for good heights.
    \begin{claim}
        There exist an integer $ m \geqslant 1 $, a number $ \eta > 0 $, and a height $ H_\epsilon > 0 $ with $ 2^m H_\epsilon \leqslant r_\ast $, all independent of the base point, such that
        \[
            0 < h < H_{\epsilon}, \quad \log \frac{\mathcal{M}_{\point{q}}(2^{-m} h)}{\mathcal{M}_{\point{q}}(2^{m} h)} < \eta \quad \Rightarrow \quad h \text{ is good at } \point{q}.
        \]
    \end{claim}
    \begin{proof}
        Suppose the claim did not hold. For each integer $ j \geqslant 1 $, we could choose a base point $ \point{q}_j $ and a bad height $ 0 < h_j < \min\{ 2^{-j} r_\ast, 1 / j \} $ such that
        \[
            0 \leqslant \log \frac{\mathcal{M}_{\point{q}_j}(2^{-j} h_j)}{\mathcal{M}_{\point{q}_j}(2^{j} h_j)} < \frac{1}{j}.
        \]
        For any fixed $ 0 < a < b $, the interval $ [a h_j, b h_j] $ lies in $ [2^{-j} h_j, 2^j h_j] $ for all sufficiently large $ j $. Monotonicity gives
        \[
            0 \leqslant \log \frac{\mathcal{M}_{\point{q}_j}(a h_j)}{\mathcal{M}_{\point{q}_j}(b h_j)} < \frac{1}{j}.
        \]
        Moreover,
        \[
            0 \leqslant \int_{a h_j}^{b h_j} \frac{\beta_{\point{q}_j}(r)}{r} \dif r \leqslant {\frac{C}{\vartheta} (b^\vartheta - a^\vartheta) h_j^\vartheta} \longrightarrow 0.
        \]
        Thus, removing the exponential correction, we obtain
        \[
            \frac{\chi_{\point{q}_j}(a h_j)}{\chi_{\point{q}_j}(b h_j)} \longrightarrow 1.
        \]
        This is the ratio condition in \cref{lemma::thesis::inclusion-with-moving-points}. Choose its parameter $ \delta > 0 $ sufficiently small that $ 1 - \delta \geqslant 2^{-\epsilon/2} $ and $ 1 + \delta \leqslant 2^{\epsilon/2} $. That lemma then implies that $ h_j $ is good for all sufficiently large $ j $, contradicting its choice.
    \end{proof}

    Fix a base point $ \point{q} $. Take any $ 0 < H < H_\epsilon $ and denote $ h_i := 2^{-i} H $. Define
    \[
        \ell_i := \log \frac{\mathcal{M}_{\point{q}}(h_{i+1})}{\mathcal{M}_{\point{q}}(h_i)} \geqslant 0.
    \]
    By the uniform bounds for $ \mathcal{M}_{\point{q}} $, for every integer $ N \geqslant 1 $,
    \[
        \sum_{i = 0}^{N - 1} \ell_i = \log \frac{\mathcal{M}_{\point{q}}(h_N)}{\mathcal{M}_{\point{q}}(H)} \leqslant \log \frac{C_1}{c_1} =: B.
    \]
    Therefore $ \sum_{i = 0}^{\infty} \ell_i \leqslant B $, independently of the base point and the starting height.

    If $ i \geqslant m $ and $ h_i $ is bad, the contrapositive of the claim gives
    \[
        \eta \leqslant \log \frac{\mathcal{M}_{\point{q}}(2^{-m} h_i)}{\mathcal{M}_{\point{q}}(2^m h_i)} = \log \frac{\mathcal{M}_{\point{q}}(h_{i+m})}{\mathcal{M}_{\point{q}}(h_{i-m})} = \sum_{l = i-m}^{i+m-1} \ell_l.
    \]
    Summing over any finite collection of bad indices, each $ \ell_l $ occurs at most $ 2m $ times. Thus the number of bad indices $ i \geqslant m $ is at most $ 2mB / \eta $. Allowing all the first $ m $ indices to be bad, the total number of bad indices is bounded by
    \[
        N_\epsilon := m + \left\lceil \frac{2mB}{\eta} \right\rceil,
    \]
    independently of the base point and the starting height.

    We next give a comparison valid also at bad heights. Decrease $ H_\epsilon $ if necessary so that $ t \leqslant \Phi(t) \leqslant 2t $ for $ 0 \leqslant t \leqslant H_\epsilon $. For $ \point{x} \in S_h[u_{\point{q}}] $, convexity gives
    \[
        w_{\point{q}}(\point{x} / 4) \leqslant \frac{1}{4} w_{\point{q}}(\point{x}) < \frac{1}{4} \Phi(h) \leqslant \frac{h}{2} \leqslant \Phi(h/2).
    \]
    Consequently,
    \[
        \frac{1}{4} S_h[u_{\point{q}}] \subset S_{h/2}[u_{\point{q}}] \subset S_h[u_{\point{q}}].
    \]
    For a compact convex body $ K $ with John ellipsoid $ z + E $, John's theorem gives $ z + E \subset K \subset z + nE $. Thus the width of $ K $ in every direction lies between the width of $ E $ and $ n $ times that width. Applying this to the closures of the preceding sections, and using that the width of an origin-centered ellipsoid is twice its support function, gives
    \[
        \frac{1}{4n} \mathcal{E}_{\point{q}}(h) \subset \mathcal{E}_{\point{q}}(h / 2) \subset n \mathcal{E}_{\point{q}}(h).
    \]

    At good heights we use the sharper comparison in the definition, and at bad heights we use this last comparison. Relative to the good-height factors, each bad step costs at most a factor $ 4n $ on either side. Since there are at most $ N_\epsilon $ bad indices, iteration gives, for every integer $ N \geqslant 1 $,
    \[
        (4n)^{-N_\epsilon} 2^{-N(\frac{1}{2} + \frac{\epsilon}{2})} \mathcal{E}_{\point{q}}(H) \subset \mathcal{E}_{\point{q}}(2^{-N} H) \subset (4n)^{N_\epsilon} 2^{-N(\frac{1}{2} - \frac{\epsilon}{2})} \mathcal{E}_{\point{q}}(H).
    \]
    Choose one integer $ N \geqslant 1 $ sufficiently large that
    \[
        2^{N\epsilon/2} \geqslant (4n)^{N_\epsilon},
    \]
    and set $ a_\epsilon := 2^{-N} $ and $ h_\epsilon := H_\epsilon $. Then
    \[
        a_\epsilon^{\frac{1}{2} + \epsilon} \mathcal{E}_{\point{q}}(H) \subset \mathcal{E}_{\point{q}}(a_\epsilon H) \subset a_\epsilon^{\frac{1}{2} - \epsilon} \mathcal{E}_{\point{q}}(H).
    \]
    Since the base point and $ 0 < H < h_\epsilon $ were arbitrary, this proves the lemma.
\end{proof}

By iterating the result of \cref{lemma::thesis::uniform-aepsilon-scale-inclusion}, we have
\begin{corollary}\label{corollary::thesis::1-e-inclusion-for-john-ellipsoid}
    For every sufficiently small $ \epsilon > 0 $, there exist $ c_{\epsilon}, C_{\epsilon} > 0 $ and $ h_{\epsilon} > 0 $, independent of the base point, such that for any $ 0 < h < h_{\epsilon} $
    \[
        B_{c_{\epsilon} h^{\frac{1}{2} + \epsilon}} \subset \mathcal{E}(h) \subset B_{C_{\epsilon} h^{\frac{1}{2} - \epsilon}}.
    \]
\end{corollary}
\begin{proof}
    The result follows by iterating \cref{lemma::thesis::uniform-aepsilon-scale-inclusion}, using the uniform diameter and volume bounds at a fixed small height to control the initial ellipsoids, and section nesting together with John's theorem to pass to intermediate heights.
\end{proof}

\begin{lemma}\label{lemma::thesis::1-e-inclusion-at-0}
    For every sufficiently small $ \epsilon > 0 $, there exist $ c_{\epsilon}, C_{\epsilon} > 0 $ and $ h_{\epsilon} > 0 $ such that for any $ 0 < h < h_{\epsilon} $
    \[
        Q \cap B_{c_{\epsilon} h^{\frac{1}{2} + \epsilon}}(0) \subset S_{h} \subset Q \cap B_{C_{\epsilon} h^{\frac{1}{2} - \epsilon}}(0).
    \]
\end{lemma}
\begin{proof}
    The idea is to transform inclusions about John ellipsoids from Corollary \ref{corollary::thesis::1-e-inclusion-for-john-ellipsoid} to sections. Once the inclusion is proved, then the $ C^{1, 1 - \epsilon} $ estimate is easy to obtain.

    For the exterior containment, since we have the John ellipsoid $ z_h+\mathcal{E}(h) $ for $ S_h $ is contained in a ball $ z_h + B_{C_{\epsilon} h^{\frac{1}{2} - \epsilon}} $ centered at $ z_h $, and by enlarging the John ellipsoid by $ n $, the enlarged ellipsoid $ z_h + n\mathcal{E}(h) $ contains $ S_h $, and $ 0 \in \cl S_h $ is also contained. Thus for any point $ q $ in $ S_h $, $ |q| \leqslant |q - z_h| + |z_h - 0| \leqslant 2n C_{\epsilon} h^{\frac{1}{2} - \epsilon} $, so we have for
    \[
        S_h[u] \subset B_{2n C_{\epsilon} h^{\frac{1}{2} - \epsilon}}.
    \]
    By abuse of notation we still denote $ 2n C_{\epsilon} $ by $ C_{\epsilon} $. The same inclusion holds for $ v $ by duality of the problem.

    For the interior inclusion, we use the exterior bound and the duality relation of the potentials. Take a boundary point $ \point{p} \in \partial S_h[v] \cap P $, then $ v(\point{p}) = h $ and $ |\point{p}| \leqslant C_{\epsilon} h^{\frac{1}{2} - \epsilon} $, so we have
    \[
        v(\point{p}) \geqslant C |\point{p}|^{\frac{1}{\frac{1}{2} - \epsilon}},
    \]
    take $ \eta $ such that $ 2 + \eta = \frac{1}{\frac{1}{2} - \epsilon} $. Since $ h $ was arbitrary, the lower bound holds for every sufficiently small $ \point{p} $: take $ h=v(\point{p}) $, which is small by continuity. Denote $ \point{q} := \point{T}^{-1}(\point{p}) $; by duality of $ u $ and $ v $ we have
    \[
        0 \leqslant v(\point{p}) \leqslant u(\point{q}) + v(\point{p}) \leqslant b(\point{q}, \point{p}) \leqslant C |\point{q}| |\point{p}|,
    \]
    thus 
    \[
        C |\point{p}|^{2 + \eta} \leqslant v(\point{p}) \leqslant C |\point{q}||\point{p}|,
    \]
    so we have
    \[
        |\point{p}| = |\point{T}(\point{q})| \leqslant C |\point{q}|^{\frac{1}{1 + \eta}}.
    \]
    Apply this bound to $ u $ again, we obtain
    \[
        0 \leqslant u(\point{q}) \leqslant b(\point{q}, \point{p}) \leqslant C|\point{q}||\point{p}| \leqslant C |\point{q}|^{\frac{2 + \eta}{1 + \eta}}.
    \]
    By the definition of $ \eta $ we have
    \[
        \frac{2 + \eta}{1 + \eta} = (1 - \frac{1}{2 + \eta})^{-1} = (\frac{1}{2} + \epsilon)^{-1},
    \]
    so if $ u(\point{q}) = h $ then 
    \[
        |\point{q}| \geqslant C h^{\frac{1}{2} + \epsilon},
    \]
    i.e. there exists a constant $ c_\epsilon $ such that
    \[
        Q\cap B_{c_{\epsilon} h^{\frac{1}{2} + \epsilon}}(0) \subset S_h[u].
    \]
    The same result holds for $ v $ by duality of the problem.
\end{proof}

And thus we obtain the global $ C^{1, 1 - \epsilon} $ regularity of the optimal potentials.

\begin{theorem}[Global $ C^{1, 1 - \epsilon} $ regularity of $ u $ and $ v $]\label{theorem::thesis::C11-e-of-u-v}
    For every sufficiently small $ \epsilon > 0 $, we have $ u \in C^{1, 1 - \epsilon}(\cl \regs) $ and $ v \in C^{1, 1 - \epsilon}(\cl \regt) $, and there exists a constant $ M $ such that
    \[
        \norm{Du}_{C^{0, 1 - \epsilon}(\cl \regs)} + \norm{Dv}_{C^{0, 1 - \epsilon}(\cl \regt)} \leqslant M.
    \]
    The constant $ M $ depends only on the dimension $ n $, $ \epsilon $, {the H\"older exponent $ \alpha $, the positive lower and upper bounds for the densities, and their $ C^{0,\alpha} $ norms}, the source and target regions $ \regs $ and $ \regt $, the norm of the cost function $ c $ and constants in the MTW conditions.
\end{theorem}
\begin{proof}
    Fix $ \epsilon>0 $ and choose $ 0<\epsilon'<1/2 $ such that
    \[
        \theta:=\frac{1-2\epsilon'}{1+2\epsilon'}\geqslant1-\epsilon.
    \]
    Applying \cref{lemma::thesis::1-e-inclusion-at-0} with $ \epsilon' $, its proof gives $ |\point{p}|\leqslant C|\point{q}|^\theta $ at transport pairs near the normalized base pair, with constants and a neighborhood independent of the base point. The cost-coordinate maps and their inverses have uniform derivative bounds. Thus, for every base pair $ (\point{x}_0,\vec{T}(\point{x}_0)) $, this estimate gives
    \[
        |\vec{T}(\point{x})-\vec{T}(\point{x}_0)|
        \leqslant C|\point{x}-\point{x}_0|^\theta
    \]
    whenever $ |\point{x}-\point{x}_0| $ is sufficiently small, uniformly in $ \point{x}_0 $. Continuity extends the estimate to boundary points, and boundedness of the target domain gives the same estimate for the remaining pairs after enlarging $ C $.

    For the original potential, the identity $ D u(\point{x})=D_{\point{x}}b(\point{x},\vec{T}(\point{x})) $ and smoothness of $ b $ imply
    \[
        \begin{aligned}
            |Du(\point{x})-Du(\point{x}_0)|
            &\leqslant C\bigl(|\point{x}-\point{x}_0|
            +|\vec{T}(\point{x})-\vec{T}(\point{x}_0)|\bigr)\\
            &\leqslant C|\point{x}-\point{x}_0|^{1-\epsilon}.
        \end{aligned}
    \]
    The gradient is also uniformly bounded by the first derivatives of $ b $. Interchanging source and target gives the same estimates for $ Dv $, proving the theorem.
\end{proof}

\subsection{Obliqueness from blow-up}

We now use the blow-up limit, following ideas from \cite[Lemma 3, 4]{collins2025boundary}. We first prove the blow-up limits of the source and target regions are half-spaces under $ C^{1, \beta} $ boundary regularity, then if obliqueness does not hold we can pick a point inside $ \cl Q_{\infty} $ perpendicular to $ \cl P_{\infty} $ and then prove by contradiction.

\begin{lemma}[Uniform obliqueness]\label{lemma::thesis::uniform-obliqueness}
    Under the assumptions of \cref{theorem::thesis::maintheorem} and with the notation of \cref{definition::thesis::setting-q-coordinate}, denote the inward unit normal of $ Q $ at $ \point{0} $ by $ \vec{\nu}_{Q} $ and that of $ P $ at $ \point{0} $ by $ \vec{\nu}_{P} $. 

    There is a constant $ \mu>0 $, independent of the base boundary transport pair, such that
    \[
        \vec{\nu}_{Q} \cdot \vec{\nu}_{P} \geqslant \mu.
    \]
\end{lemma}
\begin{proof}
    \noindent \textbf{Step 1: The blow-up regions.}  Fix the boundary transport pair and take $ s_j\downarrow0 $ at this pair as in \cref{definition::thesis::definition-of-blow-up}. The positive limit $ \chi(0) $ from \cref{lemma::thesis::estimate-of-sections-and-limit-of-chi} gives the ratio assumption in \cref{lemma::thesis::limit-chi-picture-and-zero-negative-terms}, so the homogeneity and section-convergence results apply. We have
    \[
        Q_j := \matr{A}_j^{-1} Q, \quad P_j := s_j^{-1} \matr{A}_j^{\top} P, \quad \point{T}_j(\point{q}) := s_j^{-1} \matr{A}_j^{\top} \point{T}(\matr{A}_j \point{q}).
    \]
    After the blow-up transformation, the inward unit normals of $ Q_j $ and $ P_j $ are 
    \[
        \vec{\nu}_{Q, j} = \frac{\matr{A}_j^{\top} \vec{\nu}_{Q}}{|\matr{A}_j^{\top} \vec{\nu}_{Q}|}, \quad \vec{\nu}_{P, j} = \frac{\matr{A}_j^{-1} \vec{\nu}_{P}}{|\matr{A}_j^{-1} \vec{\nu}_{P}|},
    \]
    and for the boundary defining function, suppose the function for $ Q $ near zero is $ \phi(\point{q}) = \point{q} \cdot \vec{\nu}_{Q} + R(\point{q}) $, where on a fixed $ B_{r_0} $ by $ C^{1, \beta} $ regularity we have $ |R(\point{q})| \leqslant C |\point{q}|^{1 + \beta} $. {Define $ \phi_j(\point{q}):=\phi(\matr{A}_j\point{q})/|\matr{A}_j^\top\vec{\nu}_Q| $, a defining function for $ Q_j $.}
    
    Now take $ R > 0 $, since $ R \norm{\matr{A}_j} \to 0 $, thus for all $ j $ large enough $ \phi_j $ is defined on $ B_R $ and the estimate in $ B_{r_0} $ applies. Inserting their coordinates we have
    \[
        \phi_j(\point{q}) = \point{q} \cdot \vec{\nu}_{Q, j} + R_j(\point{q}),
    \]
    {Here $ R_j(\point{q})=R(\matr{A}_j\point{q})/|\matr{A}_j^\top\vec{\nu}_Q| $, and $ \sup_{B_R}|R_j| \leqslant C_R s_j^{(1 + \beta)(1 / 2 - \epsilon) - (1 / 2 + \epsilon)} = C_R s_j^{\beta/2 - (2 + \beta) \epsilon} $.}

    Take $ \epsilon $ small to make the error exponent $ \beta/2-(2+\beta)\epsilon $ positive, then the error tends to zero. After passing to subsequences we have $ \vec{\nu}_{Q, j} \to \vec{\nu}_{Q, \infty} $, and $ Q_\infty $ is a half-space. The same result applies to $ \vec{\nu}_{P, \infty} $ and $ P_\infty $.

    \noindent \textbf{Step 2: The bound for transport map.} We want to show that $ \point{T}_j(B_r(\point{0})) $ is bounded, for the convergence argument used later in Step 3. The goal is
    \begin{claim}
        There exist $ r > 0 $ and $ M < \infty $ independent of $ j $ such that
        \[
            \diam (\point{T}_j(B_r(\point{0}) \cap \cl Q_j))     \leqslant M.
        \]
    \end{claim}
    \begin{proof}
        First recall estimates for normalized $ u, v $ and $ b $. We have
        \[
            b_{\point{q} \point{p}} (\point{0}, \point{p}) = b_{\point{q} \point{p}} (\point{q}, \point{0}) = \matr{I},
        \]
        and smoothness gives
        \[
            \norm{b_{\point{q}\point{p}} - \matr{I}} \leqslant C|\point{q}|, \quad \norm{b_{\point{q}\point{p}} - \matr{I}} \leqslant C|\point{p}|.
        \]

        After blow-up with $ s_j $, still denote the transformed functions by adding a hat for readability, we have for any fixed $ R $, for $ \point{q} \in B_R(\point{0}) $,
        \[
            \norm{\hat{b}_{\point{q}\point{p}} - \matr{I}} \leqslant C |\matr{A}_j| |\matr{A}_j^{-\top}| |\matr{A}_j \point{q}| \leqslant {C R s_j^{\frac{1}{2} - 3\epsilon}} = o(1)
        \]
        as $ j \to \infty $. So integrating back we obtain
        \begin{equation}\label{equation::thesis::estimate-used-for-bound-T}
            |D_{\point{q}} \hat{b}(\point{q}, \point{p}) - \point{p}| \leqslant o(1) |\point{p}|, \quad |\hat{b}(\point{q}, \point{p})-\point{q}\cdot\point{p}| \leqslant o(1)|\point{q}||\point{p}|,
        \end{equation}
        In particular, $ |\hat{b}(\point{q},\point{p})|\leqslant C|\point{q}||\point{p}| $.

        Now we want to show if some transport image of the normalized section $ K_j $ is too distant from $ \point{0} $, the value of $ \hat{v} $ will have a linear lower bound. First by the normalization, for the section $ K_j $, by the definition of $ \matr{A}_j $ and John's theorem, there is
        \[
            B_1(\matr{A}_j^{-1}\point{z}_j) \subset K_j \subset B_n(\matr{A}_j^{-1}\point{z}_j) \subset B_{2n}(\point{0}).
        \]
        So the volume of $ K_j $ is bounded by a constant $ C_1 $ depending only on $ n $, and the same is true of the volume of $ \point{T}_j(K_j) $ by the Monge-Amp\'ere equation. 
        
        By passing to a subsequence if needed, we have $ \matr{A}_j^{-1}\point{z}_j \to \point{z}_{\infty} $, and $ B_1(\point{z}_\infty) \subset K_{\infty} $. For $ j $ sufficiently large we have $ |\matr{A}_j^{-1}\point{z}_j - \point{z}_{\infty}| < \frac{1}{4} $ and therefore
        \[
            B := B_{\frac{1}{2}}(\point{z}_\infty) \subset K_j, K_\infty.
        \]
        Then we apply \cref{lemma::thesis::bound-for-Du-on-arbitrary-G}, there is a constant $ C_B $ such that
        \[
            \sup_{B} |{D\hat{u}}| \leqslant C_B,
        \]
        therefore for $ \point{q} \in B $, using $ {D\hat{u}(\point{q})} = D_{\point{q}} {\hat{b}}(\point{q}, {\point{T}_j(\point{q})}) $ and the estimates in \cref{equation::thesis::estimate-used-for-bound-T} we have
        \[
            |{\point{T}_j(\point{q})} - {D\hat{u}(\point{q})}| \leqslant o(1) |{\point{T}_j(\point{q})}|,
        \]
        thus
        \[
            |{\point{T}_j(\point{q})}| \leqslant |{\point{T}_j(\point{q})} - {D\hat{u}(\point{q})}| + |{D\hat{u}(\point{q})}| \leqslant o(1) |{\point{T}_j(\point{q})}| + C_B,
        \]
        so we have
        \[
            |{\point{T}_j(\point{q})}| \leqslant 2 C_B
        \]
        for $ j $ sufficiently large. This shows that the interior ball in section $ K_j $ is transported inside a bounded region. So on $ B $
        \[
            \hat{v}({\point{T}_j(\point{q})}) = \hat{b}(\point{q}, {\point{T}_j(\point{q})}) - \hat{u}(\point{q}) \leqslant \hat{b}(\point{q}, {\point{T}_j(\point{q})}) \leqslant C|\point{q}||{\point{T}_j(\point{q})}| \leqslant 4n C C_B.
        \]
        Take $ t = 4n C C_B + 1 $ then $ {\point{T}_j(\point{q})} \in S_t[\hat{v}] $. 

        We then show that $ S_t[\hat{v}] $ is bounded by some $ B_N $. First we already know that $ |S_t[\hat{v}]| \leqslant V_t $ by the dual volume estimate in \cref{lemma::thesis::estimate-of-sections-and-limit-of-chi} and the scaling in \cref{definition::thesis::definition-of-blow-up}. Since $ \point{T}_j(B)\subset S_t[\hat{v}]\cap B_{2C_B} $ has volume at least $ c|B| $, its width in any direction is at most $ 4C_B $, and Fubini's theorem gives a perpendicular slice of area at least $ c|B|/(4C_B) $. For $ \point{p}\in S_t[\hat{v}] $ with $ |\point{p}|>2C_B $, take slices perpendicular to $ \point{p} $. Convexity puts the cone joining $ \point{p} $ to the convex hull of such a slice inside $ S_t[\hat{v}] $. Its height is at least $ |\point{p}|-2C_B $, so
        \[
            \frac{c|B|}{4nC_B}(|\point{p}|-2C_B)\leqslant V_t,
        \]
        giving a bound for $ |\point{p}| $. Denote this bound by $ S_t[\hat{v}] \subset B_{C_t} $.

        Now we can prove the lower bound for $ \hat{v} $ desired. Take any $ \point{p}\in P_j $ with $ |\point{p}| > C_t $. Since $ S_t[\hat{v}]\subset B_{C_t} $, we have $ \hat{v}(C_t \point{p} / |\point{p}|) \geqslant t $, using the convexification $ \hat{z} $ of $ \hat{v} $, we have
        \[
            t \leqslant \hat{z}(C_t \point{p} / |\point{p}|) \leqslant \frac{C_t}{|\point{p}|}\hat{z}(\point{p}) \leqslant \frac{C C_t}{|\point{p}|}\hat{v}(\point{p})
        \]
        by using \cref{remark::thesis::equivalence-of-sections-of-uw}. Thus $ \hat{v}(\point{p}) \geqslant c|\point{p}| $ for some $ c $ when $ |\point{p}| > C_t $. These constants are independent of $ j $.

        Finally we use the lower bound to prove the claim. Now let $ \point{q} \in B_r(\point{0}) \cap \cl {Q_j} $ with $ r $ to be determined, and denote its transport image by $ \point{p} := {\point{T}_j(\point{q})} $. If $ |\point{p}| > C_t $ we apply the lower bound and upper bound in \cref{equation::thesis::estimate-used-for-bound-T}, having
        \[
            c|\point{p}| \leqslant \hat{v}(\point{p}) = \hat{b}(\point{q}, \point{p}) - {\hat{u}(\point{q})} \leqslant \hat{b}(\point{q}, \point{p}) \leqslant C|\point{q}||\point{p}| \leqslant Cr|\point{p}|,
        \]
        Taking $ r < c / (2C) $ gives a contradiction. Thus $ |{\point{T}_j(\point{q})}|\leqslant C_t $ on the indicated set, and we may take $ M=2C_t $ in the diameter bound.
    \end{proof}

    \noindent \textbf{Step 3: The main proof.} Now we prove the lemma by contradiction. Suppose obliqueness does not hold at the fixed boundary pair, represented by $ \point{q}=\point{0}\in\partial Q $, by assumption
    \[
        \vec{\nu}_{Q} \cdot \vec{\nu}_{P} \leqslant 0.
    \]
    We shall take $ \point{q}_j \in B_r \cap \cl Q_j$ approaching $ \point{q}_\infty := -s \vec{\nu}_{P, \infty} $ for some $ 0 < s < r $ with $ r $ defined in the claim, for example we can take $ \point{q}_j = -s \vec{\nu}_{P, j} + \delta_j \vec{\nu}_{Q, j} $, since the angle between the two normals is at least $ \pi / 2 $, so $ -s \vec{\nu}_{P, j} $ falls above the tangent plane of $ Q_j $ at $ \point{0} $ and we choose $ \delta_j := 2\sup_{B_r}|R_j|+j^{-1}\to0 $, so that $ \point{q}_j\in B_r\cap Q_j $ for all sufficiently large $ j $. Let $ \point{p}_j := {\point{T}_j(\point{q}_j)} $, by the claim $ |\point{p}_j| \leqslant M $, and 
    \[
        \point{q}_j \cdot \point{p}_j = -s \vec{\nu}_{P, j} \cdot \point{p}_j + \delta_j \vec{\nu}_{Q, j} \cdot \point{p}_j \leqslant M \delta_j \longrightarrow 0.
    \]

    But since $ \point{q}_j $ and $ \point{p}_j $ are bounded, using the defining equality for transport pairs, 
    \[
        0 \leqslant \hat{u}(\point{q}_j) = \hat{b}(\point{q}_j, \point{p}_j) - \hat{v}(\point{p}_j) \leqslant \hat{b}(\point{q}_j, \point{p}_j) = \point{q}_j \cdot \point{p}_j + o(1),
    \]
    taking $ \limsup $ with $ j \to \infty $ on both sides, we have the $ \limsup $ of RHS is zero, therefore $ \limsup_{j \to \infty} \hat{u}(\point{q}_j) = 0 $. 
    
    But since $ \point{q}_j \to \point{q}_\infty \not= \point{0} $, take a small $ t^{\prime} $ such that $ S_{t^{\prime}}[u_\infty] \subset B_{|\point{q}_{\infty}| / 2}(\point{0}) $. By the Hausdorff convergence $ \cl S_{t^{\prime}}[u^{(j)}] \to \cl S_{t^{\prime}}[u_\infty] $ we have for $ j $ sufficiently large $ \point{q}_j $ lies outside $ S_{t^{\prime}}[u^{(j)}] $, which means $ \hat{u}(\point{q}_j) \geqslant t^{\prime} > 0 $, contradicting the limit.

    Thus the scalar product is positive at every boundary transport pair. The transport map, the boundary unit normals, and the coordinate changes depend continuously on the base pair. Compactness of the boundary therefore gives a common positive lower bound $ \mu $.
\end{proof}

\section{$ C^{2, \gamma} $ regularity of the potentials}\label{section::thesis::C2gamma}

In this section we prove \cref{theorem::thesis::maintheorem}, the $ C^{2, \gamma} $ estimates for the Kantorovich potentials. Our proof follows the approach of \cite{chen2021global,collins2025boundary}, but since the cost function and thus the equations are not identical, here we give a direct proof.

\subsection{Preparation}

Let the densities be $ C^{0, \alpha} $ and the boundaries be $ C^{1, \beta} $. Define $ \gamma := \min\{ \alpha, \beta \} $. The goal is to prove
\[
    u \in C^{2, \gamma}(\cl \regs), \quad v \in C^{2, \gamma}(\cl \regt),
\]
with uniform $ C^{2, \gamma} $ estimates. Consequently, the transport map is $ \point{T} \in C^{1, \gamma}(\cl \regs) $.

By \cref{lemma::thesis::uniform-obliqueness}, a linear change of the form $ \point{q}_{\mathrm{old}}=\matr{L}\point{q} $, $ \point{p}_{\mathrm{old}}=\matr{L}^{-\top}\point{p} $ aligns both inward normals with $ \vec{e}_n $ and preserves the mixed-Hessian normalization. Locally, we may therefore write
\[
    Q=\{q_n>\phi(\point{q}^{\prime})\},\qquad
    P=\{p_n>\psi(\point{p}^{\prime})\},
\]
where
\[
    0\leqslant\phi(\point{q}^{\prime})\leqslant C|\point{q}^{\prime}|^{1+\beta},
    \qquad
    0\leqslant\psi(\point{p}^{\prime})\leqslant C|\point{p}^{\prime}|^{1+\beta}.
\]

To normalize the density ratio, set $ \lambda_0=(\rho(\point{0})/\overline{\rho}(\point{0}))^{1/n} $. Replace the target coordinate by $ \point{p}_{\mathrm{old}}=\lambda_0\point{p} $ and the functions by
\[
    u_{\mathrm{new}}(\point{q})=\lambda_0^{-1}u(\point{q}),\qquad
    v_{\mathrm{new}}(\point{p})=\lambda_0^{-1}v(\lambda_0\point{p}),\qquad
    b_{\mathrm{new}}(\point{q},\point{p})=\lambda_0^{-1}b(\point{q},\lambda_0\point{p}).
\]
The source density is unchanged and the target density becomes $ \lambda_0^n\overline{\rho}(\lambda_0\point{p}) $. We continue to use the same notation for the resulting functions and domains. This preserves
\[
    D_{\point{q}}D_{\point{p}}b(\point{q},\point{0})
    =D_{\point{q}}D_{\point{p}}b(\point{0},\point{p})=\matr{I},
\]
and makes the density ratio at the origin equal to one. The normalized benefit has the form
\[
    b(\point{q},\point{p})=\point{q}\cdot\point{p}-R(\point{q},\point{p}),
\]
with
\[
    |R(\point{q},\point{p})|\leqslant C|\point{q}|^2|\point{p}|^2,\qquad
    |D_{\point{q}}R(\point{q},\point{p})|\leqslant C|\point{q}||\point{p}|^2,\qquad
    \|D^2_{\point{q}\point{q}}R(\point{q},\point{p})\|\leqslant C|\point{p}|^2,
\]
and the corresponding estimates with source and target interchanged.

Let $ \point{p}(\point{q},\vec{\xi}) $ be determined by $ D_{\point{q}}b(\point{q},\point{p}(\point{q},\vec{\xi}))=\vec{\xi} $. Define
\[
    \matr{A}(\point{q},\vec{\xi})
    :=D^2_{\point{q}\point{q}}b(\point{q},\point{p}(\point{q},\vec{\xi})),\qquad
    B(\point{q},\vec{\xi})
    :=|\det D_{\point{q}}D_{\point{p}}b(\point{q},\point{p}(\point{q},\vec{\xi}))|
    \frac{\rho(\point{q})}{\overline{\rho}(\point{p}(\point{q},\vec{\xi}))}.
\]
The potential satisfies
\[
    \det[D^2u-\matr{A}(\point{q},Du)]=B(\point{q},Du),\qquad
    D^2u-\matr{A}(\point{q},Du)>0,\qquad B(\point{0},\point{0})=1.
\]

Choose the convexification constant $ k $ for the normalized benefit and put $ w=(e^{ku}-1)/k $, with $ w=u $ when $ k=0 $. Since
\[
    Du=\frac{Dw}{1+kw},\qquad
    D^2u=\frac{D^2w}{1+kw}-\frac{k\,Dw\otimes Dw}{(1+kw)^2},
\]
the equation for $ w $ is
\[
    \det[D^2w-\matr{A}_w(\point{q},w,Dw)]=B_w(\point{q},w,Dw),
\]
where
\[
    \begin{aligned}
        \matr{A}_w(\point{q},w,\vec{\xi})
        &=(1+kw)\matr{A}\left(\point{q},\frac{\vec{\xi}}{1+kw}\right)
        +\frac{k}{1+kw}\vec{\xi}\otimes\vec{\xi},\\
        B_w(\point{q},w,\vec{\xi})
        &=(1+kw)^n B\left(\point{q},\frac{\vec{\xi}}{1+kw}\right).
    \end{aligned}
\]
In particular, $ D^2w-\matr{A}_w=(1+kw)(D^2u-\matr{A})>0 $. Near the origin, $ |\point{p}(\point{q},\vec{\xi})|\leqslant C|\vec{\xi}| $. The benefit estimates and the H\"older continuity of the densities therefore give
\begin{equation}\label{equation::thesis::convexified-coefficient-bounds}
    \|\matr{A}_w(\point{q},w,\vec{\xi})\|\leqslant C|\vec{\xi}|^2,\qquad
    |B_w(\point{q},w,\vec{\xi})-1|
    \leqslant C\bigl(|\point{q}|^\alpha+|\vec{\xi}|^\alpha+|w|\bigr).
\end{equation}

We record the following estimate.
\begin{lemma}\label{lemma::thesis::difference-of-dw-and-t}
    Set $ \point{p} := \point{T}(\point{q}) $. There exists a constant $ C > 0 $ such that
    \[
        |\point{p} - D w(\point{q})| \leqslant C |\point{q}| |\point{p}|^2.
    \]
\end{lemma}
\begin{proof}
    By the triangle inequality,
    \[
        |\point{p} - D w(\point{q})| \leqslant |\point{p} - D u(\point{q})| + |Du(\point{q})(1 - e^{ku})|.
    \]

    Since we have normalized $ u $ and $ v $, we have
    \[
        0 \leqslant u(\point{q}) \leqslant u(\point{q}) + v(\point{p}) \leqslant b(\point{q}, \point{p}) \leqslant C_1 |\point{q}| |\point{p}|,
    \]
    so we have $ |1 - e^{ku}| \leqslant C_2 |u| \leqslant C_1 C_2 |\point{q}||\point{p}| $.

    Using the normalized form $ b = \point{q}\cdot \point{p} - R(\point{q}, \point{p}) $ and the identity $ Du(\point{q}) = D_{\point{q}}b(\point{q}, \point{p}) $, we obtain
    \[
        \point{p} - D u(\point{q}) = \point{p} - D_{\point{q}} b(\point{q}, \point{T}(\point{q})) = \point{p} - \point{p} + D_{\point{q}} R(\point{q}, \point{p}) = D_{\point{q}} R(\point{q}, \point{p}),
    \]
    so we have
    \[
        |\point{p} - D u(\point{q})| \leqslant C_3 |\point{q}| |\point{p}|^2, \quad |Du(\point{q})| = |\point{p} - D_{\point{q}} R(\point{q}, \point{p})| \leqslant C_4 |\point{p}|.
    \]

    Combining these estimates gives
    \[
        |\point{p} - D w(\point{q})| \leqslant C |\point{q}| |\point{p}|^2.
    \]
\end{proof}

\subsection{Reflection region and comparison with a solution of $ \det D^2 w_h = 1 $}

We use $ h $ to denote the section height.

Let $ a_h $ be the greatest height attained by the part of the boundary of the section $ S_h[w] $ lying on $ \partial Q $, namely
\[
    a_h:=\max\{q_n:\point{q}\in\partial Q\cap\cl S_h[w]\}.
\]
By the section bound $ S_h[w]\subset B_{C_\epsilon h^{\frac12-\epsilon}} $ and the boundary estimate $ \phi(\point{q}^{\prime})\leqslant C|\point{q}^{\prime}|^{1+\beta} $, a point $ \point{q}\in\partial Q\cap\cl S_h[w] $ attaining this maximum satisfies
\[
    a_h = \phi(\point{q}^\prime) \leqslant C |\point{q}|^{1 + \beta} \leqslant C h^{(1 + \beta)(\frac{1}{2} - \epsilon)}.
\]

The transport image lies above $ p_n=0 $, so $ T^n(\point{q})\geqslant0 $. However, 
\[
    D w(\point{q}) = e^{ku} D_{\point{q}} b(\point{q}, \point{T}(\point{q})) = e^{ku} \left[ \point{T}(\point{q}) - D_{\point{q}} R(\point{q}, \point{T}(\point{q})) \right],
\]
so
\[
    D_n w(\point{q}) = e^{ku} \left[ T^n(\point{q}) - D_{n} R(\point{q}, \point{T}(\point{q})) \right] \geqslant -C |\point{q}| |\point{T}(\point{q})|^2.
\]
We therefore add a correction term. By \cref{lemma::thesis::1-e-inclusion-at-0} and \cref{theorem::thesis::C11-e-of-u-v} we have for some $ \epsilon $
\[
    |\point{q}| \leqslant C h^{\frac{1}{2} - \epsilon}, \quad |\point{T}(\point{q})| \leqslant C |\point{q}|^{1 - \epsilon} \leqslant C h^{(1 - \epsilon)(\frac{1}{2} - \epsilon)},
\]
thus
\[
    D_n w(\point{q}) \geqslant -C h^{(3 - 2\epsilon)(\frac{1}{2} - \epsilon)}.
\]
Define $ \kappa_h := C h^{(3 - 2\epsilon)(\frac{1}{2} - \epsilon)} $ and then
\[
    \tilde{w}_h(\point{q}) := w(\point{q}) + \kappa_h (q^n - a_h)
\]
to make $ D_n \tilde{w}_h \geqslant 0 $ on $ S_h[w]\cap\{q_n=a_h\} $.


Now we define 
\[
    G_h^{+} := \{\point{q}\in\cl Q:\tilde{w}_h(\point{q})<h,\ q_n\geqslant a_h\}
\]
and denote its reflection with respect to $ q^n = a_h $ by $ G_h^{-} $. The set $ G_h := \interior (G_h^+ \cup G_h^-) $ is convex because we have proved the nonnegative normal derivative of $ \tilde{w}_h $. Above the plane, $ G_h \cap \{ q^n \geqslant a_h \} \subset S_h[w] $ since the fixing term $ \kappa_h (q^n - a_h) $ is nonnegative. So $ G_h $ shrinks to the origin as $ h \to 0 $. And there is a similar containment estimate to \cref{lemma::thesis::1-e-inclusion-at-0} for $ G_h $, namely:
\begin{lemma}\label{lemma::thesis::1-e-inclusion-for-Gh}
    Denote $ \point{q}_h := (0, \cdots, 0, a_h) $. For every sufficiently small $ \epsilon > 0 $, there exist constants $ c_{\epsilon}, C_{\epsilon} > 0 $ such that
    \[
        B_{c_{\epsilon} h^{\frac{1}{2} + \epsilon}}(\point{q}_h) \subset G_h \subset B_{C_{\epsilon} h^{\frac{1}{2} - \epsilon}}(\point{q}_h).
    \]
\end{lemma}
\begin{proof}
    The outer inclusion follows from $ G_h^+\subset\cl S_h[w] $, the section outer bound, and reflection about $ q_n=a_h $. We prove the inner inclusion.

    Choose $ \epsilon>0 $ small enough that
    \[
        (1+\beta)(\tfrac12-\epsilon)>\tfrac12+\epsilon,
        \qquad (4-2\epsilon)(\tfrac12-\epsilon)>1.
    \]
    Then $ a_h=o(h^{1/2+\epsilon}) $. By \cref{lemma::thesis::1-e-inclusion-at-0} and the equivalence of the sections of $ u $ and $ w $, choose $ c>0 $ such that
    \[
        Q\cap B_{4c h^{1/2+\epsilon}}(\point{0})\subset S_{h/4}[w].
    \]
    Set $ \point{q}_h=a_h\vec{e}_n $. For small $ h $, we have $ a_h<c h^{1/2+\epsilon} $. If $ |\point{q}-\point{q}_h|<c h^{1/2+\epsilon} $ and $ q_n\geqslant a_h $, then $ |\point{q}|<2c h^{1/2+\epsilon} $. The boundary point $ \point{q}_b=(\point{q}^{\prime},\phi(\point{q}^{\prime})) $ also lies in this latter ball for small $ h $. By continuity and the preceding section inclusion, $ w(\point{q}_b)\leqslant h/4<h $, so $ \point{q}_b\in\partial Q\cap\cl S_h[w] $. The definition of $ a_h $ gives $ \phi(\point{q}^{\prime})\leqslant a_h $. Thus $ \point{q}\in\cl Q $, and again the section inclusion gives $ w(\point{q})\leqslant h/4<h/2 $.

    On this upper half-ball,
    \[
        0\leqslant\kappa_h(q_n-a_h)
        \leqslant Ch^{(4-2\epsilon)(1/2-\epsilon)}=o(h).
    \]
    Consequently $ \tilde{w}_h(\point{q})<h/2+o(h)<h $, so the upper half-ball belongs to $ G_h^+ $. Reflection gives the lower half-ball, and taking the interior of their union proves the inner inclusion.
\end{proof}


We consider the solution to 
\[
    \begin{cases}
        \det D^2 w_h = 1, & \text{in} G_h, \\
        w_h = h, & \text{on} \partial G_h.
    \end{cases}
\]

We compare this solution with $ w $.
\begin{lemma}\label{lemma::thesis::comparison-of-w-wh}
    On $ G_h \cap Q $ we have
    \[
        |w - w_h| \leqslant C h^{1 + \gamma / 4}.
    \]
    Here $ \gamma = \min \{ \alpha, \beta \} $.
\end{lemma}
\begin{proof}
    We proceed in five steps. After establishing the required estimates, we compare $ w_h $ with $ \tilde{w}_h $ using the comparison principle on $ G_h^+ $. We then compare $ \tilde{w}_h $ with $ w $ by their definition and use reflection to obtain the estimate in the lower half. 

    \noindent \textbf{Step 1: Preparation.} We first establish some estimates needed for the comparison. Recall that
    \[
        \kappa_h = C h^{(3 - 2\epsilon)(\frac{1}{2} - \epsilon)}.
    \]

    Now we bound $ \partial_n \tilde{w}_h $ from above. Take a boundary point $ \point{q}_b = (\point{q}^{\prime}, \phi(\point{q}^{\prime})) $. Its transport image lies on the target boundary, so 
    \[
        [\point{T}(\point{q}_b)]^n = \psi([\point{T}(\point{q}_b)]^{\prime}) \leqslant C |\point{T}(\point{q}_b)|^{1 + \beta}, 
    \]
    and by \cref{lemma::thesis::difference-of-dw-and-t} we have
    \[
        |\point{T}(\point{q}) - D w(\point{q})| \leqslant C |\point{q}||\point{T}(\point{q})|^2,
    \]
    so
    \[
        \partial_n w(\point{q}_b) \leqslant C |\point{T}(\point{q}_b)|^{1 + \beta} + C |\point{q}_b||\point{T}(\point{q}_b)|^2.
    \]
    Now use the $ C^{0, 1 - \epsilon} $ estimate for $ \point{T} $ obtained in the proof of \cref{lemma::thesis::1-e-inclusion-at-0}, we have
    \[
        \partial_n w(\point{q}_b) \leqslant C h^{(1+\beta)(1-\epsilon)(\frac12-\epsilon)} + C h^{(3-2\epsilon)(\frac12-\epsilon)}.
    \]
    Now take a point in $ G_h $ and on the hyperplane, say $ \point{q} = (\point{q}^{\prime}, a_h) $. The $ C^{1, 1 - \epsilon} $ H\"older estimate from \cref{theorem::thesis::C11-e-of-u-v} gives
    \[
        |\partial_n w(\point{q}) - \partial_n w(\point{q}_b)| \leqslant C|\point{q} - \point{q}_b|^{1 - \epsilon} \leqslant  C |a_h|^{1 - \epsilon} \leqslant C h^{(1 + \beta)(\frac{1}{2} - \epsilon)(1 - \epsilon)},
    \]
    and for the additional shift of $ \tilde{w}_h $ from $ w $, we have
    \[
        \partial_n \tilde{w}_h(\point{q}) - \partial_n w(\point{q}) = \kappa_h = C h^{(3 - 2\epsilon)(\frac{1}{2} - \epsilon)}.
    \]

    Combining these estimates yields
    \[
        0 \leqslant \partial_n \tilde{w}_h \leqslant C h^{(1+\beta)(1-\epsilon)(\frac12-\epsilon)} + C h^{(3-2\epsilon)(\frac12-\epsilon)}.
    \]

    \noindent \textbf{Step 2: The equation.} Recall that
    \[
        \det[D^2w-\matr{A}_w(\point{q},w,Dw)]=B_w(\point{q},w,Dw).
    \]
    On $ G_h^+ $, the section outer bound and the $ C^{1,1-\epsilon} $ estimate give
    \[
        |\point{q}|\leqslant Ch^{\frac12-\epsilon},\qquad
        |Dw|\leqslant Ch^{(1-\epsilon)(\frac12-\epsilon)},\qquad
        0\leqslant w<h.
    \]
    By \cref{equation::thesis::convexified-coefficient-bounds},
    \[
        \|\matr{A}_w(\point{q},w,Dw)\|\leqslant Ch^{2(1-\epsilon)(\frac12-\epsilon)},\qquad
        |B_w(\point{q},w,Dw)-1|\leqslant C\bigl(h^{\alpha(1-\epsilon)(\frac12-\epsilon)}+h\bigr).
    \]
    Since $ D^2\tilde{w}_h=D^2w $ and $ Dw=D\tilde{w}_h-\kappa_h\vec{e}_n $, the equation becomes
    \[
        \det[D^2\tilde{w}_h-\matr{A}_w(\point{q},w,D\tilde{w}_h-\kappa_h\vec{e}_n)]
        =B_w(\point{q},w,D\tilde{w}_h-\kappa_h\vec{e}_n),
    \]
    where $ w=\tilde{w}_h-\kappa_h(q_n-a_h) $. Set
    \[
        \begin{aligned}
            e_A&:=\sup_{G_h^+}\|\matr{A}_w(\point{q},w(\point{q}),Dw(\point{q}))\|,\\
            e_B&:=\sup_{G_h^+}|B_w(\point{q},w(\point{q}),Dw(\point{q}))-1|,\\
            e_N&:=Ch^{(1+\beta)(1-\epsilon)(\frac12-\epsilon)}+Ch^{(3-2\epsilon)(\frac12-\epsilon)}.
        \end{aligned}
    \]
    In the comparison below, $ \matr{A}_w $ and $ B_w $ are evaluated at $ (\point{q},w(\point{q}),Dw(\point{q})) $.

    We will show
    \[
        |\tilde{w}_h - w_h| \leqslant C \left[ h^{1 - 2\epsilon}(e_A + e_B) + h^{\frac{1}{2} - \epsilon} e_N \right]
    \]
    on $ G_h^+ $.

    \noindent \textbf{Step 3: Construct comparison functions.} Write $ \partial G_h^+ = \mathcal{C}_1 \cup \mathcal{C}_2 $, where $ \mathcal{C}_1 \subset \{ q^n > a_h \} $ is the boundary of the section and $ \mathcal{C}_2 \subset \{ q^n = a_h \} $ is on the hyperplane. We know $ w_h = \tilde{w}_h = h $ on $ \mathcal{C}_1 $, and by symmetry we know $ \partial_n w_h = 0 $, while $ 0 \leqslant \partial_n \tilde{w}_h \leqslant e_N $ on $ \mathcal{C}_2 $.

    Take $ C_{\ast} $ as the constant $ C_{\epsilon} $ in \cref{lemma::thesis::1-e-inclusion-for-Gh}, then we have
    \[
        |\point{q} - a_h \vec{e}_n|^2 \leqslant C_{\ast}^2 h^{1 - 2 \epsilon},
    \]
    and 
    \[
        0 \geqslant q_n-a_h - C_{\ast} h^{\frac{1}{2} - \epsilon} \geqslant - C_{\ast} h^{\frac{1}{2} - \epsilon}.
    \]

    We define some auxiliary functions, 
    \[
        Q_h(\point{q}) := |\point{q} - a_h \vec{e}_n|^2 - C_{\ast}^2 h^{1 - 2 \epsilon}, \quad L_h(\point{q}) := q_n-a_h - C_{\ast} h^{\frac{1}{2} - \epsilon}, 
    \]
    then throughout $ G_h^+ $
    \[
        \begin{gathered}
            Q_h \leqslant 0, \quad |Q_h| \leqslant C_{\ast}^2 h^{1 - 2 \epsilon} \\
            \partial_n Q_h = 0 \text{ on } \mathcal{C}_2, \quad D^2 Q_h = 2 \matr{I}, \\
            L_h \leqslant 0, \quad |L_h| \leqslant C_{\ast} h^{\frac{1}{2} - \epsilon}, \\
            \partial_n L_h = 1.
        \end{gathered}
    \]

    \noindent \textit{1.} We first deal with the upper bound of $ \tilde{w}_h - w_h $. Consider the function
    \[
        \hat{w} := \tilde{w}_h + s Q_h + \delta L_h
    \]
    with $ \delta > 0 $ small and  $ s = C(e_A + e_B) + \delta $, here $ C $ is a constant to be determined. Since both $ Q_h $ and $ L_h $ are nonpositive over $ G_h^+ $, we have $ \hat{w} \leqslant \tilde{w}_h $ on the boundary $ \partial G_h^+ $. And along $ \mathcal{C}_2 $ we have
    \[
        \partial_n (\hat{w} - w_h) = \partial_n \tilde{w}_h + \delta > 0
    \]
    The parameter $ \delta $ ensures strict positivity and can be arbitrarily small. And we also have
    \[
        D^2 \hat{w} = D^2 \tilde{w}_h + 2s \matr{I}.
    \]
    Recall that $ D^2 \tilde{w}_h - \matr{A}_{w} = D^2 w - \matr{A}_{w} \geqslant 0 $, and $ \matr{A}_w \geqslant -e_{A} \matr{I} $, so
    \[
        D^2 \hat{w} = (D^2 \tilde{w}_h - \matr{A}_{w}) + (\matr{A}_{w} + e_{A} \matr{I}) + (2s - e_{A}) \matr{I} \geqslant (D^2 \tilde{w}_h - \matr{A}_{w}) + (2s - e_{A}) \matr{I},
    \]
    Choosing $ C $ so that $ 2s - e_A>0 $ and using the concavity of $ \det^{1/n} $ gives
    \[
        \det D^2 \hat{w} \geqslant \left[ \det(D^2 \tilde{w}_h - \matr{A}_{w})^{\frac{1}{n}} + (2s - e_{A}) \right]^n \geqslant \left[ (1 - e_B)^{\frac{1}{n}} + (2s - e_{A}) \right]^n,
    \]
    now we enlarge $ C $ so that $ (1 - e_B)^{\frac{1}{n}} + (2s - e_{A}) > 1 $, then $ \det D^2 \hat{w} \geqslant 1 = \det D^2 w_h $. By the comparison principle we have
    \[
        \hat{w} \leqslant w_h,
    \]
    which is
    \[
        \tilde{w}_h - w_h \leqslant s |Q_h| + \delta |L_h|.
    \]
    Letting $ \delta \to 0 $, by the upper bounds of $ |Q_h| $ and $ |L_h| $ we have
    \[
        \tilde{w}_h - w_h \leqslant C h^{1 - 2 \epsilon} (e_A + e_B).
    \]

    \noindent \textit{2.} We then deal with the lower bound of $ \tilde{w}_h - w_h $. Consider the function
    \[
        \check{w} := \tilde{w}_h - s Q_h - (e_N + \delta)L_h,
    \]
    with the same $ \delta $ and $ s $. Then $ \check{w} \geqslant \tilde{w}_h $ on the boundary $ \partial G_h^+ $. On $ \mathcal{C}_2 $ we have
    \[
        \partial_n (\check{w} - w_h) = \partial_n \tilde{w}_h - (e_N + \delta) \leqslant -\delta < 0.
    \]
    These inequalities rule out possible negative minima of $ \check{w} - w_h $ on the boundary $ \partial G_h^+ $. Suppose there exists a negative minimum inside, then at that point 
    \[
        D^2 \check{w} \geqslant D^2 w_h \geqslant 0,
    \]
    so we have
    \[
        D^2 \tilde{w}_h - \matr{A}_w \geqslant (-\matr{A}_w + e_A \matr{I}) + (2s - e_A) \matr{I} + D^2 w_h \geqslant D^2 w_h + (2s - e_A) \matr{I},
    \]
    the same determinant inequality then gives
    \[
        \det(D^2\tilde{w}_h - \matr{A}_w) \geqslant \left[ (\det D^2 w_h)^{\frac{1}{n}} + (2s - e_A) \right]^n = \left[ 1 + (2s - e_A) \right]^n,
    \]
    now we enlarge $ C $ again so that $ \left[ 1 + (2s - e_A) \right]^n > 1 + e_B $, then this is a contradiction since the LHS is
    \[
        |\det(D^2\tilde{w}_h - \matr{A}_w)| = |B_w| \leqslant 1 + e_B.
    \]
    Therefore a negative interior minimum cannot exist, and
    \[
        \tilde{w}_h - w_h \geqslant s Q_h + (e_N + \delta)L_h.
    \]
    Letting $ \delta \to 0 $ again and by the upper bounds of $ |Q_h| $ and $ |L_h| $ we have
    \[
        \tilde{w}_h - w_h \geqslant -C h^{1 - 2 \epsilon} (e_A + e_B) - C h^{\frac{1}{2} - \epsilon} e_N.
    \]

    Substituting the estimates for $ e_A,e_B,e_N $, and combining terms with the same power, gives
    \[
        \begin{aligned}
            |\tilde{w}_h-w_h|\leqslant C\bigl[{}
            &h^{(4-2\epsilon)(\frac12-\epsilon)}
            +h^{1-2\epsilon+\alpha(1-\epsilon)(\frac12-\epsilon)}\\
            &+h^{2-2\epsilon}
            +h^{[1+(1+\beta)(1-\epsilon)](\frac12-\epsilon)}\bigr].
        \end{aligned}
    \]
    At $ \epsilon=0 $, these exponents are $ 2 $, $ 1+\alpha/2 $, $ 2 $, and $ 1+\beta/2 $, respectively. Since $ \gamma=\min\{\alpha,\beta\}>0 $, we may choose $ \epsilon>0 $ sufficiently small that all four exponents exceed $ 1+\gamma/4 $, while retaining the smallness conditions in the preceding construction. Thus
    \[
        |\tilde{w}_h - w_h| \leqslant C h^{1 + \frac{\gamma}{4}}.
    \]

    \noindent \textbf{Step 4: Comparison with $ w $.} Using the estimate for $ \kappa_h $ and the outer bound for $ G_h $, we have
    \[
        0 \leqslant \tilde{w}_h - w = \kappa_h (q^n - a_h) \leqslant C h^{(4-2\epsilon)(\frac12-\epsilon)}.
    \]
    This is $ o(h^{1+\gamma/4}) $. Thus, on $ \cl G_h^+ $,
    \[
        |w - w_h| \leqslant C h^{1 + \frac{\gamma}{4}}.
    \]

    \noindent \textbf{Step 5: The lower half.} Take a point $ \point{q} \in G_h^- \cap Q $ and denote its reflection with respect to $ q^n = a_h $ by $ \point{q}^{+} $. We have proved that
    \[
        |w(\point{q}^+) - w_h(\point{q}^+)| \leqslant C h^{1 + \frac{\gamma}{4}}.
    \]
    Since $ \point{q}\in Q $ lies below the cutting plane and $ q_n\geqslant0 $, we have $ |\point{q}-\point{q}^+|\leqslant2a_h $. The segment joining these points lies in $ Q $, and the $ C^{1,1-\epsilon} $ estimate gives
    \[
        \begin{aligned}
            |w(\point{q})-w(\point{q}^+)|
            &\leqslant \sup_{[\point{q},\point{q}^+]}|Dw|\,|\point{q}-\point{q}^+|\\
            &\leqslant C h^{(1-\epsilon)(\frac12-\epsilon)}h^{(1+\beta)(\frac12-\epsilon)}
            =C h^{(2+\beta-\epsilon)(\frac12-\epsilon)}.
        \end{aligned}
    \]
    Choose $ \epsilon $ sufficiently small that $ (2+\beta-\epsilon)(\frac12-\epsilon)>1+\gamma/4 $; this is possible since the exponent at $ \epsilon=0 $ is $ 1+\beta/2 $.
    By symmetry $ w_h(\point{q}^+) = w_h(\point{q}) $, so
    \[
        |w(\point{q}) - w_h(\point{q})| \leqslant |w(\point{q}) - w(\point{q}^+)| + |w(\point{q}^+) - w_h(\point{q}^+)| \leqslant C h^{1 + \frac{\gamma}{4}}.
    \]
\end{proof}

\subsection{Construction of the approximating sequence}

We construct a sequence of solutions approximating $ w $, as in \cite{chen2021global}. To use the same normalization, we establish uniform shape estimates for the domains and bound the difference of adjacent Hessians using \cite[Lemma 6.4]{chen2021global}.

Fix a sufficiently small initial height $ h_0 $ and start with $ G_{h_0} $. Define
\[
    h_k := 4^{-k} h_0, \quad \mathcal{G}_{k} := G_{4^{-k} h_0},
\]
and let $ w_k := w_{h_k} $ solve
\[
    \begin{cases}
        \det D^2 w_k = 1, & \text{in } \mathcal{G}_k, \\
        w_k = h_k, & \text{on } \partial \mathcal{G}_k.
    \end{cases}
\]

We first prove that $ \mathcal{G}_k \sim \mathcal{G}_{k + 1} $, namely there is a constant $ K $ depending only on $ n $ such that
\[
    K^{-1} \mathcal{G}_{k} \subset \mathcal{G}_{k + 1} \subset K \mathcal{G}_{k}.
\]
\begin{proof}
    \noindent \textit{1.} Pick a point $ \point{q} \in G_h $, we want to show that $ \frac{1}{8} \point{q} \in G_{h / 4} $. If $ \point{q} $ lies in the lower half, denote its reflection into the upper half by $ \point{q}^+ $; otherwise, set $ \point{q}^+=\point{q} $.

    Now take an axis point $ \point{p} := 4a_{h / 4} \vec{e}_n $, and consider
    \[
        \point{z} := \frac{1}{8} \point{q}^+ + \frac{7}{8} \point{p},
    \]
    Since both $ \point{q}^+ $ and $ \point{p} $ belong to $ \cl Q $, convexity gives $ \point{z}\in\cl Q $. Now we show that $ \point{z} \in G_{h / 4} $. For this we must evaluate $ \tilde{w}_{h / 4} $ at this point. We first check its $ q^n $ coordinate, 
    \[
        z^n \geqslant \frac{1}{8} a_h + \frac{7}{2} a_{h / 4} \geqslant a_{h / 4}.
    \]
    Thus $ \point{z} $ lies on or above $ q^n=a_{h/4} $, and $ \tilde{w}_{h/4}(\point{z}) $ is evaluated on the upper half. By definition, 
    \[
        \tilde{w}_{h / 4}(\point{z}) = w(\point{z}) + \kappa_{h / 4} |z^n - a_{h / 4}|,
    \]
    the first term can be bounded by the convexity of $ w $, that is
    \[
        w(\point{z}) \leqslant \frac{1}{8} w(\point{q}^+) + \frac{7}{8} w(\point{p}).
    \]
    By definition of $ \point{q}^+ $, $ w(\point{q}^+)\leqslant h $. For the other term, the $ C^{1,1-\epsilon} $ estimate gives $ w(\point{p}) \leqslant C |a_{h / 4}|^{2 - \epsilon} \leqslant C h^{(1 + \beta)(\frac{1}{2} - \epsilon)(2 - \epsilon)} = o(h) $. Therefore
    \[
        w(\point{z}) \leqslant \frac{1}{8}h + o(h).
    \]
    For the correction term, we also use the outer bound
    $ |\point{q}^+| \leqslant C h^{\frac12-\epsilon} $. Since
    \[
        0 \leqslant z^n-a_{h/4}
        = \frac18 (\point{q}^+)^n+\frac52 a_{h/4}
        \leqslant C h^{\frac12-\epsilon},
    \]
    the bound for $\kappa_{h/4}$ gives
    \[
        \kappa_{h/4}|z^n-a_{h/4}| \leqslant C h^{(4-2\epsilon)(\frac12-\epsilon)} = o(h).
    \]
    Therefore
    \[
        \tilde{w}_{h / 4}(\point{z}) \leqslant \frac{1}{8}h + o(h) < \frac{1}{4}h.
    \]
    This shows that $ \point{z} \in G_{h / 4} $.

    We now show that $ \frac{1}{8} \point{q} \in G_{h / 4} $. Let $ \point{z}^- $ be the reflection of $ \point{z} $ with respect to $ q^n = a_{h / 4} $. Both points belong to $ G_{h / 4} $, so does the vertical segment joining them. Since the horizontal coordinate of the segment is $ \frac{1}{8} \point{q}^{\prime} $, we only need to check that the vertical coordinate of $ \frac{1}{8} \point{q} $ lies between $ \point{z} $ and $ \point{z}^- $. Using $ (\point{q}^+)^n \geqslant q^n $ and $ (\point{q}^+)^n + q^n \geqslant 2 a_h $,
    \[
        \begin{gathered}
            z_n - \frac{1}{8} q^n = \frac{1}{8} ((\point{q}^+)^n - q^n) + \frac{7}{2} a_{h / 4} \geqslant 0,  \\
            \frac{1}{8} q^n - z_n^- = \frac{1}{8} (q^n + (\point{q}^+)^n) + \frac{3}{2} a_{h / 4} \geqslant 0.
        \end{gathered}
    \]
    Therefore $ \frac{1}{8} \point{q} \in G_{h / 4} $, i.e. $ \frac{1}{8} G_{h} \subset G_{h / 4} $.

    \noindent \textit{2.} We now pick a $ \point{q} \in G_{h / 4} $ and show that $ \frac{1}{8} \point{q} \in G_{h} $. This time we pick $ \point{p} := 4 a_h \vec{e}_n $ and define $ \point{z} $ in the same way. The proof is the same as part 1, and the main difference is now we have
    \[
        \tilde{w}_{h}(\point{z}) \leqslant \frac{1}{32}h + o(h) < h,
    \]
    so the inclusion $ \point{z} \in G_h $ holds.
\end{proof}

With this inclusion result, we may take $ K^{-2} \mathcal{G}_{k + 2} \subset \mathcal{G}_{k} \cap \mathcal{G}_{k + 1} $ as a region where both $ w_k $ and $ w_{k + 1} $ are defined. By previous results, we have
\[
    \sup_{K^{-2} \mathcal{G}_{k + 2} \cap Q} |w - w_k| \leqslant C h_k^{1 + \gamma / 4}.
\]

We first show that the difference of adjacent $ w_k $ has the same bound.
\begin{lemma}\label{lemma::thesis::bound-of-wk-wk1}
    We have on the whole $ K^{-2} \mathcal{G}_{k + 2} $,
    \[
        \sup_{K^{-2} \mathcal{G}_{k + 2}} |w_k - w_{k + 1}| \leqslant C h_k^{1 + \gamma / 4}.
    \]
\end{lemma}
\begin{proof}
    Put $ h=h_k $. For points inside $ Q $, the triangle inequality and the comparison with $ w $ give the result. It remains to treat points outside $ Q $. Take such $ \point{q} \in K^{-2} \mathcal{G}_{k + 2} \setminus Q $, then $ \point{q} $ must be below $ q^n = a_h $ and $ q^n = a_{h / 4} $ since the upper half is contained in $ Q $. So we take reflections respectively, say $ \point{q}_h^+ $ and $ \point{q}_{h / 4}^+ $. By the reflection symmetries,
    \[
        \begin{aligned}
            |w_k(\point{q})-w_{k+1}(\point{q})|
            &\leqslant
            |w_k(\point{q}_h^+)-w(\point{q}_h^+)| \\
            &\quad+
            |w(\point{q}_h^+)-w(\point{q}_{h/4}^+)| \\
            &\quad+
            |w(\point{q}_{h/4}^+)-w_{k+1}(\point{q}_{h/4}^+)| \\
            &\leqslant
            C h^{1+\gamma/4}
            +|w(\point{q}_h^+)-w(\point{q}_{h/4}^+)|.
        \end{aligned}
    \]

    To estimate the remaining term, observe that
    \[
        |\point{q}_h^+-\point{q}_{h/4}^+|
        =2|a_h-a_{h/4}|
        \leqslant C h^{(1+\beta)(\frac12-\epsilon)}.
    \]
    Both upper reflected points satisfy the section outer bound
    \[
        |\point{q}_h^+|+|\point{q}_{h/4}^+|
        \leqslant C h^{\frac12-\epsilon}.
    \]
    Now use the $C^{1,1-\epsilon}$ estimate
    \[
        \sup_{[\point{q}_h^+,\point{q}_{h/4}^+]}|Dw|
        \leqslant C h^{(1-\epsilon)(\frac12-\epsilon)}.
    \]
    Thus
    \[
        \begin{aligned}
            |w(\point{q}_h^+)-w(\point{q}_{h/4}^+)|
            &\leqslant
            \sup_{[\point{q}_h^+,\point{q}_{h/4}^+]}|Dw|\,
            |\point{q}_h^+-\point{q}_{h/4}^+| \\
            &\leqslant
            C h^{(2+\beta-\epsilon)(\frac12-\epsilon)}
            \leqslant C h^{1+\gamma/4}.
        \end{aligned}
    \]
    Combining the estimates proves the lemma.
\end{proof}

We next estimate the difference of the Hessians of $ w_k $ and $ w_{k+1} $.
\begin{lemma}\label{lemma::thesis::bound-for-D2-difference-at-0}
    There holds
    \[
        |D^2 w_{k}(\point{0}) - D^2 w_{k + 1}(\point{0})| \leqslant {C h_k^{\frac{\gamma}{4} - 2\epsilon}}.
    \]
\end{lemma}
\begin{proof}
    The proof has two steps. We first prove a $ C^4 $ estimate for $ w_k $ after scaling and then apply \cite[Lemma 6.4]{chen2021global} to bound $ |D^2 w_{k} - D^2 w_{k + 1}| $.

    \noindent \textbf{Step 1: Normalization.} For each $ h < h_0 $, take the John ellipsoid of $ G_h $, say $ z_h + \matr{J}_h B_1 $, {with $ \matr{J}_h $ symmetric positive definite}. Take $ \matr{T}_h := \frac{\sqrt{h}}{(\det \matr{J}_h)^{\frac{1}{n}}} \matr{J}_h $, then normalize $ G_h $ and the function $ w_h $ as
    \[
        \overline{G}_h := \matr{T}_h^{-1} G_h, \quad \overline{w}_h(\point{q}) := h^{-1} w_h(\matr{T}_h \point{q}),
    \]
    The normalized equation is
    \[
        \det D^2 \overline{w}_h = h^{-n}(\det \matr{T}_h)^2 \det D^2 w_h = 1.
    \]

    Now we want to give a $ C^4 $ bound independent of $ h $ for $ \overline{w}_h $, so we need a uniform shape estimate: 
    \begin{claim}
        There are fixed $ c_\ast $ and $ C_\ast $ independent of $ h $ such that for all small $ h < h_0 $,
        \[
            B_{c_\ast}(\point{0}) \subset \overline{G}_h \subset B_{C_\ast}(\point{0}).
        \]
    \end{claim}
    \begin{proof}
        By John's theorem, 
        \[
            z_h + \matr{J}_h B_1 \subset G_h \subset z_h + n\matr{J}_h B_1.
        \]
        After applying $ \matr{T}_h^{-1} $, we have
        \[
            \overline{G}_h \subset \matr{T}_h^{-1} z_h + n \frac{(\det \matr{J}_h)^{\frac{1}{n}}}{\sqrt{h}} B_1, 
        \]
        since we have $ G_h^+ \subset S_h[w] $ and $ |S_h[w]| \leqslant C h^{\frac{n}{2}} $, $ \det \matr{J}_h \leqslant C h^{\frac{n}{2}} $ and therefore the outer radius $ n \frac{(\det \matr{J}_h)^{\frac{1}{n}}}{\sqrt{h}} \leqslant Cn $. So $ \diam \overline{G}_h \leqslant 2Cn $, giving the $ C_{\ast} $ ball.

        For the inner ball we use the function $ \overline{w}_h $. The comparison from \cref{lemma::thesis::comparison-of-w-wh} gives
        \[
            |\overline{w}_h(\point{0})| = h^{-1} |w_h(\point{0}) - w(\point{0})| \leqslant C h^{\gamma / 4},
        \]
        so by choosing $ h_0 $ sufficiently small, we have $ 1 - \overline{w}_h(\point{0}) \geqslant \frac{1}{2} $. Now we apply the Alexandrov maximum principle to $ \overline{w}_h - 1 $ over $ \overline{G}_h $, and we have 
        \[
            |1 - \overline{w}_h(\point{0})|^n \leqslant C \diam (\overline{G}_h)^{n - 1} \dist(\point{0}, \partial \overline{G}_h) \int_{\overline{G}_h} |\det D^2\overline{w}_h|,
        \]
        using the diameter we obtained above we have 
        \[
            {\dist(\point{0}, \partial \overline{G}_h) \geqslant \frac{2^{-n}}{C C_{\ast}^{n - 1} |\overline{G}_h|} \geqslant \frac{1}{2^n C\omega_n C_{\ast}^{2n-1}} =: c_{\ast}>0.}
        \]
        {Here $ \omega_n=|B_1| $, and we used $ |\overline{G}_h|\leqslant\omega_n C_\ast^n $, so $ c_\ast $ is independent of $ h $.} Therefore we have the inner ball $ B_{c_\ast}(\point{0}) $.
    \end{proof}
    
    {By the uniform inner and outer ball bounds, a further affine normalization has its linear part and inverse bounded in terms of $ c_\ast,C_\ast $ and puts the domain between $ B_{1/n} $ and $ B_1 $. Apply \cite[Lemma 5]{ONEILL2021178} to $ \overline{w}_h-1 $ after this change and the corresponding rescaling of the function that preserves $ \det D^2w=1 $. Scaling back, for any $ 0<r<c_\ast $ there is a constant $ C $, depending only on $ n,r,c_\ast,C_\ast $, such that}
    \[
        \norm{\overline{w}_h}_{C^4(B_{r}(\point{0}))} \leqslant C.
    \]
    {The upper Hessian bound and $ \det D^2\overline{w}_h=1 $ also give a constant $ \lambda>0 $ with the same dependence such that, on $ B_r(\point{0}) $,}
    \[
        \lambda^{-1} \matr{I} \leqslant D^2 \overline{w}_h \leqslant \lambda \matr{I}.
    \]

    \noindent \textbf{Step 2: Apply the estimates.} Normalize $ \mathcal{G}_{k} $ using the preceding construction, and denote $ \overline{w}_k := \overline{w}_{h_k} $. {We also write $ \matr{T}_k:=\matr{T}_{h_k} $.} By the estimate we have
    \[
        \norm{\overline{w}_k}_{C^4(B_{r}(\point{0}))} \leqslant C.
    \]
    
    We express adjacent solutions in the same coordinates. Put $\matr{T}:=\matr{T}_{h_k}$ and define
    \[
    v_{k+1}(\point{q}) := h_k^{-1}w_{k+1}(\matr{T}\point{q}).
    \]
    Notice that $v_{k+1}$ uses the normalization of $\mathcal{G}_k$, rather than its own normalization. Since $ \det\matr{T}=h_k^{n/2} $, it satisfies
    \[
    \det D^2v_{k+1}=1
    \quad\text{in }\matr{T}^{-1}\mathcal{G}_{k+1},
    \qquad
    v_{k+1}=\frac14
    \quad\text{on }\partial(\matr{T}^{-1}\mathcal{G}_{k+1}).
    \]

    The consecutive-domain inclusions imply
    \[
        K^{-1} \matr{T}^{-1}\mathcal{G}_k \subset \matr{T}^{-1}\mathcal{G}_{k+1}  \subset K \matr{T}^{-1}\mathcal{G}_k,
    \]
    Therefore,
    \[
    B_{c_\ast/K}(\point{0}) \subset \matr{T}^{-1}\mathcal{G}_{k+1} \subset B_{KC_\ast}(\point{0}).
    \]
    We can thus apply the same $ C^4 $ interior estimate to $ v_{k+1} - \frac{1}{4} $ on its own domain $\matr{T}^{-1}\mathcal{G}_{k+1}$.

    Choose $ r  > 0 $ independent of $ k $ such that $2r < c_\ast / K^4$. The preceding interior estimates then give
    \[
        \norm{\overline{w}_k}_{C^4(B_{2r}(\point{0}))}
        +\norm{v_{k+1}}_{C^4(B_{2r}(\point{0}))}
        \leqslant C.
    \]
    Here $ C $ is still independent of $ k $ since we can use $ B_{c_{\ast} / K} $ and $ B_{C_\ast K}(\point{0}) $ to control both regions.

    Since $ K^{-2} \mathcal{G}_{k} \subset \mathcal{G}_{k + 2} $, we have
    \[
        B_{2r} \subset K^{-4} \matr{T}^{-1} \mathcal{G}_k \subset K^{-2} \matr{T}^{-1} \mathcal{G}_{k + 2},
    \]
    and recall \cref{lemma::thesis::bound-of-wk-wk1} we have
    \[
        \sup_{B_{2r}}|\overline{w}_k - v_{k+1}| \leqslant h_k^{-1} \sup_{K^{-2}\mathcal{G}_{k+2}}|w_k-w_{k+1}| \leqslant C h_k^{\gamma/4}.
    \]

    Now we shall apply \cite[Lemma 6.4]{chen2021global}, after rescaling by the fixed radius $ r $, we have for $ m = 0, 1, 2, 3 $,
    \[
        \norm{D^m\overline{w}_k - D^m v_{k + 1}}_{L^\infty(B_r)} \leqslant C h_k^{\frac{\gamma}{4}}.
    \]

    Scaling back, by definition
    \[
        (D^2 w_k - D^2 w_{k + 1})(\point{q}) = h_k \matr{T}^{-\top} \left[ D^2\overline{w}_k - D^2 v_{k + 1} \right](\matr{T}^{-1} \point{q}) \matr{T}^{-1}.
    \]

    By \cref{lemma::thesis::1-e-inclusion-for-Gh} we have $ \norm{\matr{T}_{h_k}^{-1}} \leqslant C h_k^{-\frac{1}{2} - \epsilon} $, substituting this into the equation above, 
    \[
        |D^2 w_k(\point{0}) - D^2 w_{k + 1}(\point{0})| \leqslant C h_k^{\frac{\gamma}{4} - 2\epsilon}.
    \]
\end{proof}

Choose $ \epsilon $ so that $ \gamma/4-2\epsilon > 0 $. For every $ m > k $, summing the estimate at $ \point{0}  $ gives
\[
    \begin{aligned}
        |D^2w_m(\point{0})-D^2w_k(\point{0})|
        &\leqslant
        \sum_{j=k}^{m-1}
        |D^2w_{j+1}(\point{0})-D^2w_j(\point{0})|\\
        &\leqslant
        C\sum_{j=k}^{m-1}h_j^{\gamma/4-2\epsilon}\\
        &\leqslant C h_k^{\gamma/4-2\epsilon}.
    \end{aligned}
\]
So $ D^2w_k(\point{0}) $ converges and is uniformly bounded from above. Since each matrix is positive definite and $ \det D^2w_k(\point{0})=1 $, its eigenvalues are also uniformly bounded from below. Therefore, there exist constants $ c,C>0 $, independent of $ k $, such that
\[
    c\matr{I} \leqslant D^2w_k(\point{0}) \leqslant C\matr{I}.
\]

\subsection{Proof of $ C^{2, \gamma} $ regularity}

We can now sharpen the section bounds in \cref{lemma::thesis::1-e-inclusion-at-0}.
\begin{lemma}\label{lemma::thesis::sqrt-inclusion-for-sections}
    There exist constants $ c,C>0 $ such that, for all sufficiently small $ h>0 $,
    \[
        B_{c \sqrt{h}}(\point{0}) \cap Q \subset S_h[w] \subset B_{C \sqrt{h}}(\point{0}) \cap Q.
    \]
\end{lemma}
\begin{proof}
    \noindent \textbf{Step 1.} We first establish
    \begin{claim}
        For each $ k $ we have
        \[
            B_{c \sqrt{h_k}}(\point{0}) \subset \mathcal{G}_k \subset B_{C \sqrt{h_k}}(\point{0}).
        \]
    \end{claim}
    \begin{proof}
        Recall that we now have
        \[
            c\matr{I} \leqslant D^2w_k(\point{0}) \leqslant C\matr{I}.
        \]
        The normalized domain lies between two fixed concentric balls. By \cref{lemma::thesis::bound-for-D2-difference-at-0}, the normalized solution satisfies
        \[
            \lambda^{-1} \matr{I} \leqslant D^2 \overline{w}_{h_k}(\point{0}) \leqslant \lambda \matr{I}.
        \]
        For any vector $ \vec{\xi} $, the definition gives
        \[
            h_k \vec{\xi}^{\top} D^2 \overline{w}_{h_k}(\point{0}) \vec{\xi} = \vec{\xi}^\top \matr{T}_k^\top D^2 w_k(\point{0}) \matr{T}_k \vec{\xi}, 
        \]
        the $ \lambda $ bound gives
        \[
            \lambda^{-1}h_k |\vec{\xi}|^2 \leqslant \vec{\xi}^\top \matr{T}_k^\top D^2 w_k(\point{0}) \matr{T}_k \vec{\xi} \leqslant \lambda h_k |\vec{\xi}|^2.
        \]
        And by the bound for $ D^2 w_k $ we have
        \[
            C |\matr{T}_k \vec{\xi}|^2 \geqslant \lambda^{-1} h_k  |\vec{\xi}|^2, \quad c |\matr{T}_k \vec{\xi}|^2 \leqslant \lambda h_k |\vec{\xi}|^2,
        \]
        which implies
        \[
            \frac{\lambda^{-1}}{C} h_k |\vec{\xi}|^2 \leqslant |\matr{T}_k \vec{\xi}|^2 \leqslant \frac{\lambda}{c} h_k |\vec{\xi}|^2.
        \]
        Therefore the eigenvalues of $ \matr{T}_k $ are bounded above and below by fixed positive multiples of $ \sqrt{h_k} $.

        We now use the inner and outer balls of $ \matr{T}_k^{-1} \mathcal{G}_k $ proved in the claim of \cref{lemma::thesis::bound-for-D2-difference-at-0}, namely
        \[
            \matr{T}_k B_{c_\ast}(\point{0}) \subset \mathcal{G}_k \subset \matr{T}_k B_{C_\ast}(\point{0}),
        \]
        and the bounds for the eigenvalues of $ \matr{T}_k $ give
        \[
            B_{c \sqrt{h_k}}(\point{0}) \subset \mathcal{G}_k \subset B_{C \sqrt{h_k}}(\point{0})
        \]
        for some $ c $, $ C $ independent of $ k $.
    \end{proof}

    \noindent \textbf{Step 2.} Now we transfer this bound to sections of $ w $ at $ h_k $. 

    For the inner inclusion, $ w_k\leqslant h_k $ in $ \mathcal{G}_k $. Comparing $ w_k $ with $ w $ on $ B_{c\sqrt{h_k}}\cap Q\subset\mathcal{G}_k\cap Q $ gives 
    \[
        w \leqslant |w - w_k| + w_k \leqslant C h_k^{1 + \frac{\gamma}{4}} + h_k \leqslant 2h_k,
    \]
    so we have 
    \[
        B_{c \sqrt{h_k}} \cap Q \subset S_{2h_k}[w].
    \]

    For the exterior containment, take $ \point{q} \in S_{h_k / 2}[w] $, 
    \begin{itemize}
        \item If $ \point{q} $ is above or on $ q^n = a_{h_k} $, then
        \[
            \tilde{w}_{h_k}(\point{q}) = w(\point{q}) + \kappa_{h_k} |q^n - a_{h_k}| \leqslant \frac{1}{2} h_k + o(h_k),
        \]
        so $ \point{q} \in \mathcal{G}_k \subset B_{C \sqrt{h_k}}(\point{0}) $.
        \item If $ \point{q} $ is below $ q^n = a_{h_k} $, let $ \point{q}^{\sharp}:=(\point{q}^{\prime},a_{h_k}) $ be its vertical projection onto $ q^n=a_{h_k} $. Then
        \[
            w(\point{q}^{\sharp}) \leqslant w(\point{q}) + \sup_{[\point{q},\point{q}^{\sharp}]}|Dw|\,|q^n-a_{h_k}| \leqslant \frac{1}{2} h_k + o(h_k),
        \]
        and since the correction term in $ \tilde{w}_{h_k} $ vanishes on the plane $ q^n = a_{h_k} $, we have $ \point{q}^{\sharp} \in \mathcal{G}_k $. Therefore
        \[
            |\point{q}| \leqslant |\point{q} - \point{q}^{\sharp}| + |\point{q}^{\sharp}| \leqslant |a_{h_k}| + C\sqrt{h_k} \leqslant C\sqrt{h_k}
        \]
        since $ a_{h_k}=o(\sqrt{h_k}) $.
    \end{itemize}
    So we have
    \[
        S_{h_k / 2}[w] \subset B_{C \sqrt{h_k}} \cap Q.
    \]
    
    \noindent \textbf{Step 3.} Finally we fill the gaps. For an arbitrarily small $ h $, choose $ k $ with
    \[
        \frac{1}{8} h_k \leqslant h \leqslant \frac{1}{2} h_k,
    \]
    then
    \[
        S_h[w] \subset S_{\frac{1}{2} h_k}[w] \subset B_{C \sqrt{h_k}} \cap Q \subset B_{C \sqrt{8h}} \cap Q,
    \]
    and choose $ k $ with
    \[
        2 h_k \leqslant h \leqslant 8 h_k,
    \]
    then
    \[
        S_h[w] \supset S_{2 h_k}[w] \supset B_{c \sqrt{h_k}} \cap Q \supset B_{c \sqrt{h / 8}} \cap Q
    \]

    After adjusting $ c $ and $ C $, we obtain both inclusions for $ S_h[w] $.
\end{proof}

Having established the $ \sqrt{h} $ section bounds on both sides, with the target estimate obtained by duality, we repeat the preceding comparison argument with these improved estimates. We record the resulting estimates.

The source and target section bounds give quadratic lower bounds for the potentials normalized at any boundary transport pair. At a transport pair, $ u(\point{q})+v(\point{p})=b(\point{q},\point{p})\leqslant C|\point{q}|\,|\point{p}| $, so the target lower bound gives $ |\point{p}|\leqslant C|\point{q}| $. Thus, in any fixed coordinate neighborhood,
\[
    |Dw(\point{q})-Dw(\point{z})|\leqslant C|\point{q}-\point{z}|,
    \qquad \point{z}\in\partial Q,\quad \point{q}\in\cl Q,
\]
for nearby points, uniformly in the boundary point $ \point{z} $. Applying this estimate at the boundary point $ (\point{q}^{\prime},\phi(\point{q}^{\prime})) $ in the normal-derivative calculation gives, {for $ \point{q}\in S_h[w]\cap\{q_n=a_h\} $,}
\[
    \begin{gathered}
        a_h\leqslant Ch^{(1+\beta)/2},  \\
        -Ch^{3/2}\leqslant D_nw(\point{q}) \leqslant Ch^{(1+\beta)/2}, \\
    \end{gathered}
\]
We therefore repeat the construction of fixing $ \tilde{w}_h $ with $ \kappa_h = C h^{3/2} $, and now
\[
    \sup_{G_h \cap \cl Q} |w - w_h| \leqslant Ch^{1 + \gamma / 2}.
\]
Here $ G_h $ and $ w_h $ are defined using the new $ \kappa_h $. From now on, $ \mathcal{G}_k:=G_{h_k} $ and $ w_k:=w_{h_k} $ refer to this updated construction.

The resulting $ G_h $ then satisfies $ B_{c\sqrt{h}} \subset G_h \subset B_{C\sqrt{h}} $. Hence the preceding interior estimates using \cite[Lemma 6.4]{chen2021global} now yield
\begin{equation}\label{equation::thesis::bound-for-Dm-w-w}
    \norm{D^m(w_{k+1}-w_k)}_{L^\infty(B_{c\sqrt{h_k}}(\point{0}))}
    \leqslant Ch_k^{1 + \frac{\gamma}{2} - \frac{m}{2}},
    \qquad m=0,1,2,3.
\end{equation}

We now prove the quadratic expansion of $ w $, following \cite{chen2021global}.
\begin{lemma}
    There is a positive definite matrix $ \matr{H} $ with $ \det \matr{H} = 1 $ and $ c \matr{I} \leqslant \matr{H} \leqslant C \matr{I} $, such that for $ \point{q} \in \cl Q $ near $ \point{0} $,
    \[
        \left| w(\point{q}) - \frac{1}{2} \point{q}^{\top} \matr{H} \point{q} \right| \leqslant C |\point{q}|^{2 + \gamma}.
    \]
\end{lemma}
\begin{proof}
    {By \cref{equation::thesis::bound-for-Dm-w-w}, the sequences $ w_k(\point{0}) $, $ Dw_k(\point{0}) $, and $ D^2w_k(\point{0}) $ are Cauchy and therefore admit limits.} Denote their limits by
    \[
        {\ell_0 := \lim_{k \to \infty} w_k(\point{0}), \quad \vec{\ell} := \lim_{k \to \infty} Dw_k(\point{0}), \quad \matr{H} := \lim_{k \to \infty} D^2 w_k(\point{0}).}
    \]
    Then we have
    \[
        \begin{gathered}
            |w_k(\point{0})-\ell_0| \leqslant C h_k^{1 + \gamma/2},  \\
            |Dw_k(\point{0}) - \vec{\ell}| \leqslant Ch_k^{(1+\gamma)/2}, \\
            \norm{D^2w_k(\point{0})-\matr{H}} \leqslant Ch_k^{\gamma/2}.
        \end{gathered}
    \]
    Taking limits also gives $ c \matr{I} \leqslant \matr{H} \leqslant C \matr{I} $ as well as $ \det \matr{H} = 1 $.

    For $ |\point{q}|\leqslant c\sqrt{h_k} $, we have
    \[
        D^2 w_k(\point{q}) - D^2 w_k (\point{0}) = D^2 w_0 (\point{q}) - D^2 w_0(\point{0}) + \sum_{i=0}^{k-1} \left[ D^2(w_{i+1} - w_i)(\point{q}) - D^2(w_{i+1}-w_i)(\point{0}) \right].
    \]
    The $ m = 3 $ estimate gives
    \[
        \left| D^2(w_{i+1} - w_i)(\point{q}) - D^2(w_{i+1}-w_i)(\point{0}) \right| \leqslant C h_i^{-\frac{1}{2} + \frac{\gamma}{2}} |\point{q} - \point{0}|,
    \]
    and the $ C^4 $ bound for $ w_0 $ gives
    \[
        |D^2 w_0 (\point{q}) - D^2 w_0(\point{0})| \leqslant C |\point{q}|.
    \]
    Therefore
    \[
        |D^2 w_k(\point{q}) - D^2 w_k (\point{0})| \leqslant C |\point{q}| \left[ 1 + \sum_{i = 0}^{k - 1} h_i^{-\frac{1}{2} + \frac{\gamma}{2}} \right] \leqslant C h_{k}^{\frac{\gamma}{2}}.
    \]

    Define the polynomial
    \[
        p(\point{q}) := \ell_0 + \vec{\ell} \cdot \point{q} + \frac{1}{2} \point{q}^{\top} \matr{H} \point{q}.
    \]
    Then, in $ B_{c \sqrt{h_k}} $
    \[
        |w_k(\point{q}) - p(\point{q})| \leqslant C h_k^{1 + \frac{\gamma}{2}}.
    \]
    Since this ball lies in $ \mathcal{G}_k $ and the comparison error is bounded by $ C h_k^{1+\gamma/2} $, we have in $ B_{c \sqrt{h_k}} \cap \cl Q $
    \[
        |w(\point{q}) - p(\point{q})| \leqslant C h_k^{1 + \frac{\gamma}{2}}.
    \]

    For any nonzero $ \point{q}\in\cl Q $ sufficiently near zero, choose $ k $ such that $ |\point{q}| \leqslant c \sqrt{h_k} < 2|\point{q}| $ and adjust the coefficient if necessary, 
    \[
        |w(\point{q}) - p(\point{q})| \leqslant C |\point{q}|^{2 + \gamma}.
    \]
    Since $ w(\point{0}) = 0 $ and $ Dw(\point{0}) = \point{0} $, the lemma is proved.
\end{proof}

Since the remaining argument is similar to that in \cite[Proof of Theorem 1.1]{chen2021global}, we present it briefly.

Repeating the preceding construction at each boundary point $ \point{z} $ gives quadratic approximations $ p_{\point{z}} $ with uniform error $ C|\point{q}-\point{z}|^{2+\gamma} $. In a fixed coordinate neighborhood, comparison on a common interior ball yields $ |\matr{H}_{\point{z}}-\matr{H}_{\point{z}^{\prime}}|\leqslant C|\point{z}-\point{z}^{\prime}|^\gamma $ for their Hessians.

We next obtain uniform interior Hessian bounds by the rescaling argument in \cite[Proof of Theorem 2.2, Step 3]{chen2014boundaryepsilon}. Fix an interior point $ \point{q}_1 $, set $ d=\dist(\point{q}_1,\partial Q) $, and choose a nearest boundary point $ \point{z} $. Normalize at the transport pair based at $ \point{z} $. Since $ u-w=O(|\point{q}-\point{z}|^4) $, the preceding quadratic expansion also holds for $ u $. A reciprocal linear change, with uniformly bounded matrix and inverse, makes its quadratic part $ |\point{q}|^2/2 $ and preserves the mixed-Hessian normalization. Choose a fixed small $ \kappa>0 $ and rescale source and target by $ r=\kappa d $ around $ (\point{q}_1,\point{T}(\point{q}_1)) $, dividing the potential and benefit by $ r^2 $ and subtracting the supporting functions. The source ball of radius $ 4r $ lies in the interior. The quadratic expansion and semiconvexity give, uniformly as $ d\to0 $,
\[
    \left\|u_r-\frac12|\cdot|^2\right\|_{C^0(B_3)}=o(1),
    \qquad
    \|b_r-\point{q}\cdot\point{p}\|_{C^{2,\alpha}(B_3\times B_3)}=o(1).
\]
Here the gradient convergence follows by bounding directional difference quotients using the quadratic expansion and semiconvexity. Hence the rescaled transport map converges uniformly to the identity on $ B_2 $. Restrict the transported densities to the closed source ball $ \cl B_2 $ and its image, retaining the potential on $ B_3 $. Since the map is a homeomorphism, both sets contain $ B_{1/3} $ and lie in $ B_3 $ for small $ d $. Restriction preserves the transport relation and $ c $-convexity. The rescaled densities have uniformly bounded $ C^{0,\alpha} $ norms and converge uniformly to the same positive constant, by the density normalization and $ \det\matr{H}_{\point{z}}=1 $. Dividing both by that constant verifies the hypotheses of \cite[Theorem 5.3]{dephilippis2015partial}, with cost $ -b_r $. Its interior estimate gives $ |D^2u_r(\point{0})|\leqslant C $, independently of $ d $. Scaling back gives the same uniform bound for $ D^2u(\point{q}_1) $. The compact interior is covered by \cref{lemma::thesis::integrable-derivatives-for-monotonicity}; applying the same argument to $ v $ proves global $ C^{1,1} $ bounds for both potentials, and hence for $ w $.

These $ C^{1,1} $ bounds give uniform ellipticity: the source bound controls $ h_{ij} $, while the target bound and $ h_{\overline{k}\overline{l}}=b_{i\overline{k}}h^{ij}b_{j\overline{l}} $ control $ h^{ij} $. For an interior point $ \point{q} $, put $ d=\dist(\point{q},\partial Q) $ and choose a nearest boundary point $ \point{z} $. Normalize at $ \point{z} $ and rescale $ B_{d/2}(\point{q}) $ by $ d $. The boundary approximation has rescaled error $ O(d^\gamma) $, and \cref{equation::thesis::convexified-coefficient-bounds} gives rescaled $ C^{0,\alpha} $ coefficient errors $ O(d^\alpha+d^2) $. The interior estimates for uniformly elliptic concave equations, followed by the linear Schauder estimate for the difference from $ p_{\point{z}} $, therefore give
\[
    |D^2w(\point{q})-\matr{H}_{\point{z}}|\leqslant Cd^\gamma,
    \qquad
    [D^2w]_{C^{0,\alpha}(B_{\theta d}(\point{q}))}
    \leqslant Cd^{\gamma-\alpha}
\]
for a fixed $ \theta>0 $; see \cite[Chapters 6 and 17]{gilbarg1998elliptic}. These estimates are preserved under the smooth coordinate changes.

We now conclude as in \cite[Proof of Theorem 1.1, following (6.27)]{chen2021global}: use the interior estimate for pairs of nearby interior points, and compare through nearest boundary points for the remaining pairs. Thus $ w\in C^{2,\gamma}(\cl Q) $. Inverting the convexification and the coordinate changes gives $ u\in C^{2,\gamma}(\cl\regs) $, and the same conclusion holds for $ v $ by duality.

\newpage
\printbibliography

\end{document}